\documentclass[a4paper,final, reqno, 10pt]{amsart}

\usepackage[x11names]{xcolor} 

\usepackage{ifdraft}

\usepackage{etoolbox}

\makeatletter
\newcommand{\myitem}[1]{%
	\item[#1]\protected@edef\@currentlabel{#1}%
}
\makeatother

\usepackage{xfrac}

\usepackage[margin=2.5cm]{geometry} 

\usepackage{amsmath, bm, mathtools}
\usepackage{amsthm}
\usepackage{amssymb}

\numberwithin{equation}{section}

\usepackage{enumitem}

\usepackage{tikz-cd}
\usetikzlibrary{calc}

\usepackage[final]{hyperref}
\hypersetup{
	unicode,
	colorlinks=true,
	linkcolor=blue,
	citecolor=Green4, 
	filecolor=magenta,      
	urlcolor=blue,
	pdfauthor={Author One, Author Two, Author Three}, 
	pdfkeywords={article, template, simple},
	pdfproducer={LaTeX}
}

\theoremstyle{plain}
\newtheorem{theorem}{Theorem}[section]
\newtheorem{corollary}[theorem]{Corollary}
\newtheorem{lemma}[theorem]{Lemma}

\newtheorem{proposition}[theorem]{Proposition}

\newtheorem{definition}[theorem]{Definition}
\theoremstyle{definition}

\newtheorem{remark}[theorem]{Remark}

\usepackage{graphicx, color}
\graphicspath{{fig/}}

\usepackage{mathrsfs} 
\usepackage{xfrac}

\usepackage{lipsum}

\usepackage{scalerel}
\newlength\bshft
\def\fakebold#1{\ThisStyle{\ooalign{$\SavedStyle#1$\cr%
			\kern-\bshft$\SavedStyle#1$\cr%
			\kern\bshft$\SavedStyle#1$}}}

\title{Local well-posedness of the Cahn--Hilliard--Biot System}

\author{Helmut Abels} 
\address{Fakultät für Mathematik, Universität Regensburg, 93040 Regensburg, Germany}
\email{\href{mailto:helmut.abels@ur.de}{helmut.abels@ur.de}}

\author{Jonas Haselböck} 
\address{Fakultät für Mathematik, Universität Regensburg, 93040 Regensburg, Germany}
\email{\href{mailto:jonas.haselboeck@ur.de}{jonas.haselboeck@ur.de}}

\date{}

\newcommand\numberthis{\addtocounter{equation}{1}\tag{\theequation}}
\newcommand{\norm}[1]{ \left\|#1 \right\| }
\newcommand{\bnorm}[1]{\big\| #1 \big\|}

\newcommand{\abs}[1]{ |#1 |}
\newcommand{\N}{\mathbb{N}}
\newcommand{\inN}{\in \N}
\newcommand{\X}{\bm{W}^{1, q}_{\Gamma_D}}
\newcommand{\Xp}{\bm{W}^{1, q'} }
\newcommand{\Xd}{ {\bm{W}_{\Gamma_D}^{-1, q}} }
\newcommand{\OT}{\Omega_T}

\newcommand{\vphi}{\varphi}
\newcommand{\eps}{\varepsilon}
\newcommand{\Tau}{\mathcal{T}}
\newcommand{\pt}{\partial_t}

\newcommand{\Wp}{W_{,\vphi}}
\newcommand{\E}{\mathcal{E}}
\newcommand{\A}{\mathcal{A}}
\newcommand{\TA}{\tilde{\mathcal{A}}}
\newcommand{\cL}{\mathcal{L}}
\newcommand{\WE}{W_{,\E}}
\newcommand{\F}{\mathcal{F}}

\newcommand{\R}{\mathbb{R}}
\newcommand{\bu}{\bm{u}}

\newcommand{\dx}{\, d\mathbf{x}}

\newcommand{\vr}{\varrho}

\newcommand{\C}{\fakebold{\mathbb{C}}}

\newcommand{\new}[1]{{#1}}

\def\XXint#1#2#3{{\setbox0=\hbox{$#1{#2#3}{\int}$ }
		\vcenter{\hbox{$#2#3$ }}\kern-.6\wd0}}

\DeclareRobustCommand{\SkipTocEntry}[5]{}

\newcommand{\fakepart}[1]{%
	\par\refstepcounter{part}
	\addtocontents{toc}{\protect\contentsline{section}{\protect\numberline{}#1}{}{}}
}

\begin{document}
	
		\begin{abstract} 
		We show short-time well-posedness of a diffuse interface model describing the flow of a fluid through a deformable porous medium consisting of two phases.
		The system non-linearly couples Biot's equations for poroelasticity, including phase-field dependent material properties, with the Cahn--Hilliard equation to model the evolution of the solid, where we further distinguish between the absence and presence of a visco-elastic term of Kelvin-Voigt type. 
		While both problems will be reduced to a fixed-point equation that can be solved using maximal regularity theory along with a contraction argument, the first case relies on a semigroup approach over suitable Hilbert spaces, whereas treating the second case under minimal assumptions with respect to spatial regularity necessitates the application of Banach scales.
	\end{abstract}

	\maketitle

\noindent
\textbf{Key words:} Cahn--Hilliard equation, Biot’s equations, poroelasticity, well-posedness, maximal regularity\\ 
\textbf{AMS-Classification:} 35A01, 35A02, 35K41, 74B10, 76S05

\tableofcontents

\thispagestyle{empty}

\section{Introduction and main results} \label{sec:derivation}

Uniting two well-established, empirically verified models into one system, the Cahn--Hilliard--Biot equations combine poroelasticity and phase-separation dynamics into a three-way coupled system describing fluid flow through a saturated, heterogeneous medium, that in particular accounts for mutual interactions between pore pressure, mechanical deformation, and varying material properties. This exposition aims to establish well-posedness in both the purely elastic and the visco-elastic case.\par
Among a myriad of possible applications for the Cahn--Hilliard--Biot model, the most exciting is certainly tumor modeling, where associating the values $\pm 1$ of the order parameter $\vphi$ with malignant and healthy tissue, respectively, allows for the inclusion of dynamically changing characteristics such as {(visco-\nolinebreak)}elastic responsiveness or porosity.
The coupling to a deformation field, which is also influenced by the pore pressure of interstitial fluid, further permits the introduction of growth-inhibiting effects due to mechanical stress. \\ 
In this regard, various models based on the Cahn--Hilliard system coupled to different flow fields and designed for application to tumor growth were studied in past years \cite{ebenbeck2019analysis, colli2017asymptotic, garcke16darcy, lowengrub2013analysis, MR4126782}, with more recent contributions in e.g. \cite{trautwein22visco, cavalleri2024phasefield}.
On the other hand, we refer to \cite{bociu2023mathematical, bociu2016analysis,  mow1980biphasic, sacco2019comprehensive} and the references cited therein for studies of poro-(visco-)elastic systems, and point out that results on nonlinear couplings of the Biot model can be found in \cite{MR1876882, MR2102316, MR3794344, MR4649995}. \par
The Cahn--Hilliard--Biot model was derived by Storvik et al. \cite{STORVIK2022107799} as the gradient flow of a generalized Ginzburg-Landau energy, and consist of the fourth order, parabolic Cahn--Hilliard system which is nonlinearily coupled to the elliptic momentum equation for the elastic displacement, and a parabolic fluid flow system satisfying Darcy's law. Analytic results, including existence of weak solutions and uniqueness thereof in certain situation\new{s}, can be found in \cite{haselboeck2024existence, riethmüller2023wellposedness, fritz2023wellposedness}, all of which rely on a visco-elastic regularization for the most general case. With respect to numerical results, we reference \cite{MR4832711, fritz2023wellposedness, storvik2024sequential}, and point to \cite{MR4873132} for formal sharp-interface asymptotics in the absence of visco-elasticity. \par 
\medskip 
Before we introduce the coupled system, let us review the two relevant subsystems individually. 

\addtocontents{toc}{\SkipTocEntry}
\subsection{Cahn--Hilliard dynamics}
Originally, the Cahn--Hilliard system 
\begin{alignat*}{2}
	\partial_t \vphi  &= \nabla \cdot \left( m(\vphi ) \nabla \mu \right)  \quad & & \textnormal{in }  (0, T) \times \Omega, \\
	\mu &= -\eps \Delta \vphi + \tfrac{1}{\eps} \psi'(\vphi) \quad & & \textnormal{in }  (0, T) \times \Omega, 
\end{alignat*}
was proposed in 1958 \cite{cahn1958free} by John W. Cahn  and John E. Hilliard as a model for phase-separation in binary alloys that arises as the $H^{-1}$-gradient flow of the Ginzburg-Landau energy 
\begin{equation*}
	 \int_\Omega \frac{\eps}{2} |\nabla \vphi|^2 + \frac{1}{\eps} \psi(\vphi) \dx , 
\end{equation*}
\new{which describes the surface energy of the interface between the two phases and can be regarded as the free energy of the mixture.}
Here, $\vphi : [0, T] \times \Omega \rightarrow \R$ is the phase-field variable which attains values in the interval $[-1, 1]$, thus indicating an appropriate mixture of the pure components that correspond to $\pm 1$, such that $ |\nabla \vphi|^2$ penalizes rapid changes in the concentration. \new{Associated to $\vphi$ is the so-called \textit{chemical potential} $\mu$, which is given by the first variation of the energy}. To enforce the desired phase-separation dynamics, the free-energy density $\psi$ typically features a double-well shape, e.g., 
\begin{equation*}
	\psi(\vphi) = \alpha (1 - \vphi^2)^2, \quad \alpha \in \R_{>0}, 
\end{equation*} 
penalizing $\vphi$ for deviations from the pure phases. Such regular potentials, however, \new{cannot} guarantee the order parameter $\vphi$ to stay within the bounds of physically meaningful values, i.e, $[-1, 1]$, and are therefore often replaced by a logarithmic potential, cf. \cite{cahn1958free}, or an obstacle potential, cf.  \cite{PhysRevA.38.434,blowey1991cahn}. Combined, this results in a diffuse interface model which is characterized by the presence of a thin layer, the thickness of which being proportional to the parameter $\varepsilon$, between all areas of opposing pure phases where the concentration experiences a continuous transition between the values $\pm 1$.  For analytic results on the Cahn--Hilliard system we refer to, e.g., \cite{blowey1991cahn, cahn1996cahn, elliot_garcke00}. \par 
Observing that phase-separation processes are not only driven by interface mechanics, but rather affected by diverse phenomena, Cahn and Larché \cite{larcht1982effect}, and later Onuki \cite{onuki1989ginzburg}, proposed the inclusion of elastic effect by adding the \new{linearized elastic energy}
\begin{equation*}
	 \int_\Omega \new{\frac{1}{2}} \C(\vphi) (\E(\bm{u})- \Tau(\vphi)):  (\E(\bm{u}) - \Tau(\vphi)) \dx 
\end{equation*}
to the Ginzburg-Landu energy, where $\bm{u} : [0, T] \times \Omega \rightarrow \R^n$ denotes displacement in the lattice with respect to a reference configuration. Under the assumption of infinitesimally small deformations, the linearized strain tensor can thus be approximated by the symmetrized gradient $\E(\bm{u}) = \frac{1}{2} (\nabla \bm{u} + \nabla \bm{u}^{T})$, from which we subtract the phase-field dependent eigenstrain $\Tau(\vphi)$, i.e., the strain the material would attain if it was homogeneous and unstressed. Lastly, the elasticity tensor $\C$ models the stiffness of the material as a function of the order parameter $\vphi$. \\ 
Assuming that the mechanical equilibrium is attained at a much faster time-scale than the diffusion processes take place, and therefore that the displacement $\bm{u}$ is quasi-static, the $H^{-1}$-gradient flow of the combined energies yields the parabolic-elliptic Cahn-Larchè system,
which has been extensively studied over the last decades, cf. e.g. \cite{bonetti2002model,MR1807441, garcke_2003, GARCKE2005165}. \par 
\medskip 

\addtocontents{toc}{\SkipTocEntry}
\subsection{Biot dynamics}
Conversely, Biot's \cite{biot1941general} model, describing the interactions between a deformable, porous, solid matrix and the fluid that saturates it, combines a mass conserving diffusion equation based on Darcy's law with linear elasticity to ensure the balance of forces. Following the exposition in \cite{SHOWALTER2000310, gurvich24weak, bociu2016analysis} and denoting the fluid pore pressure by $p$, the poroelasticity system in the quasi-static case reads
\begin{equation}\label{eq:biot}
	\begin{alignedat}{2}
		-\nabla \cdot  \big( \C \E(\bm{u}) \big) + \alpha \nabla p &= \bm{f}  \quad & & \textnormal{in }  (0, T) \times \Omega, \\ 
		\pt \big( \tfrac{1}{M} p +  \alpha \nabla \cdot \bm{u}  \big) + \nabla (\kappa \nabla p) &= S_f  \quad & & \textnormal{in }  (0, T) \times \Omega,
	\end{alignedat}
\end{equation} 
where $\bm{f}$ relates to external forces and $S_f$ is a source term. The coupling coefficient $\alpha$, also known as Biot-Willis coefficient \cite{coussy2004poromechanics, biot1957elastic}, indicates the sensitivity to mutual effects between pressure and deformation. More precisely, the term $\alpha \nabla p$ introduces pore pressure as an additional cause for stress to the momentum equation, whereas $ \alpha \nabla \cdot \bm{u}$ corresponds to changes in the fluid content related to local volume fluctuation. Moreover, the variable $\kappa$ combines material properties of the solid and the fluid, more precisely permeability and viscosity, respectively, to create the diffusion parameter of the Darcy flow which is driven by the pressure gradient. Similarly, $\frac{1}{M}$ unites the fluid's compressibility with the solid's porosity to an effective parameter that provides a relation between pressure increments and the increase of fluid within the medium. Fundamental analytical results pertaining this system include \cite{auriault1980dynamic,vzenivsek1984existence, SHOWALTER2000310}. \\ 
The porosity $\phi = \phi(\bm{x}, t)$ at any given point $\bm{x} \in \Omega$ and time $t$ is given as the fraction
\begin{equation*}
	 \phi(\bm{x}, t) = \frac{V_f (\bm{x}, t)}{\new{ V_f(\bm{x}, t)} + V_s(\bm{x}, t)}, 
\end{equation*}
where $V_f$ and $V_s$ denote the relative volume occupied by the solid and the fluid component, respectively, in an infinitesimally small neighborhood of $\bm{x}$. Hence, we can define the function for the volumetric fluid content $\theta : [0, T] \times \Omega \rightarrow \R$ as the deviation of the porosity from a reference value $\phi_0$
\begin{equation*}
	\theta(\bm{x}, t) \coloneqq  \phi(\bm{x}, t) - \phi_0(\bm{x}), 
\end{equation*}
and assuming the constitutive relation $\theta = \frac{1}{M}p +  \alpha \nabla \cdot \bm{u}$, which implies $p =  M( \theta - \alpha \nabla \cdot \bm{u})$, we obtain that \eqref{eq:biot} is equivalent to the elliptic-parabolic system  
\begin{equation}\label{eq:biot_2}
	\begin{alignedat}{2}
		-\nabla \cdot  \big( \C\, \E(\bm{u}) \big) + \alpha M \nabla ( \theta + \alpha \nabla \cdot \bm{u}) &= \bm{f}  \quad & & \textnormal{in }  (0, T) \times \Omega, \\ 
		\pt  \theta + \nabla (\kappa \nabla p) &= S_f  \quad & & \textnormal{in }  (0, T) \times \Omega,\\ 
		p & =  M( \theta - \alpha \nabla \cdot \bm{u})  \quad & & \textnormal{in }  (0, T) \times \Omega. 
	\end{alignedat}
\end{equation} 
\new{
	In this context, we can define the fluid energy 
	\begin{equation*}
		\int_\Omega \frac{M}{2} (\theta - \alpha \nabla \cdot \bm{u})^2 \dx, 
	\end{equation*}
	and observe that the integrand is proportional to the square of the pressure $p$. Taking a gradient flow of this energy then leads to the already mentioned Darcy law for the evolution of the fluid. 
}
\\  
While there are various possibilities to include visco-elastic effects of the solid matrix into the model, see \cite{ biot1956theory, coussy2004poromechanics, SHOWALTER2000310,  bociu2016analysis, both21global}, this exposition focuses on visco-elasticity of Kelvin-Voigt type, which is obtained by augmenting the momentum equation with the divergence of the strain rate tensor $\nabla \cdot ( \C_{\nu} \E(\pt \bm{u})$ such that 
\begin{equation*}
	- \nabla \cdot ( \C_{\nu} \E(\pt \bm{u}) -\nabla \cdot  \big( \C \E(\bm{u}) \big) + \alpha \nabla p = \bm{f}  \quad \textnormal{in }  (0, T) \times \Omega,
\end{equation*}
where $\C_{\nu}$ denotes the visco-elasticity tensor.

\addtocontents{toc}{\SkipTocEntry}
\subsection{The Cahn--Hilliard-Biot system}
We reiterate that $\vphi, \mu , \theta, p : [0, T] \times \Omega \rightarrow \R$ are scalar functions in time and space, and that $\bm{u} : [0, T] \times \Omega \rightarrow \R^n $ is vector valued. 
As mentioned above, the following system was derived by Storvik et al.~in \cite{STORVIK2022107799} as a gradient flow of the energy 
	\begin{align}
		 \int_\Omega \frac{\eps}{2} |\nabla \vphi|^2 + \frac{1}{\eps} \psi(\vphi) \dx 
		+ \int_\Omega W(\vphi, \E(\bm{u})) \dx 
		+ \int_\Omega  \frac{M(\vphi)}{2} (\theta - \alpha(\vphi) \nabla \cdot \bm{u})^2 \dx, 
	\end{align}
	and augmented by source terms and external forces
	\begin{subequations} \label{eq:strong_formulation}
	\begin{align}
		\partial_t \vphi  &= \nabla \cdot \left( m(\vphi ) \nabla \mu \right) + S_s(\vphi, \E(\bu), \theta)  \label{eq:strong_formulation_phi} & & \textrm{in }  (0, T) \times \Omega ,\\
		&\begin{aligned}[b]
			\! \! \! \mu = -\eps \Delta \vphi &+ \frac{1}{\eps} \psi'(\vphi) + W_{,\vphi} (\vphi, \E(\bu)) \\&-M(\vphi) (\theta - \alpha(\vphi) \nabla \cdot \bm{u}) \alpha'(\vphi) \nabla \cdot \bm{u}  
			 +  \frac{M'(\vphi)}{2} (\theta - \alpha(\vphi) \nabla \cdot	 \bm{u})^2
		\end{aligned} & & \textrm{in } (0, T) \times \Omega , \label{eq:strong_formulation_mu} \\ 
		 - \nabla \cdot \bm{\sigma} &= \bm{f} & & \textrm{in } (0, T) \times \Omega , \label{eq:strong_formulation_u} \\ 
		\bm{\sigma} &= W_{,\E} (\vphi, \E(\bu))  + \vr\, \C_{\nu}(\vphi)  \E( \partial_t \bm{u})  - \alpha (\vphi) M(\vphi) (\theta - \alpha(\vphi) \nabla \cdot \bm{u})\, \bm{I} & & \textrm{in } (0, T) \times \Omega ,\label{eq:strong_formulation_u_2} \\ 
		\partial_t \theta &=  \nabla \cdot (\kappa(\vphi) \nabla p) + S_f(\vphi, \E(\bu), \theta) 	& & \textrm{in } (0, T) \times \Omega , \label{eq:strong_formulation_theta}\\   
		p& = M(\vphi) (\theta - \alpha (\vphi) \nabla \cdot \bm{u}) & & \textrm{in }  (0, T) \times \Omega. \label{eq:strong_formulation_p}
	\end{align}
	Here, $\vr \in \{0, 1 \}$ is a parameter indicating the presence of the visco-elastic regularization.
	We would like to point out that the original derivation did not include the linear visco-elastic term ${\nabla \cdot ( \C_\nu \E(\pt \bm{u}))}$, which is purely dissipative therefore does not interfere with the gradient flow structure the system exhibits in the absence of source terms and applied forces.\\
	\new{Moreover, we denote by $\Wp, \WE$ the derivatives of the elastic energy density with respect to $\vphi$ and $\E(\bm{u})$, respectively.} \par 
	As for initial conditions, we impose 
	\begin{align}
		\vphi (0) = \vphi_0  \quad  \textrm{and} \quad
		\varrho\bm{u}(0) = \varrho\bm{u}_0 \quad  \textrm{and} \quad 
		\theta(0) = \theta_0  \quad  \textrm{in } \Omega . 
	\end{align}
	For the boundary conditions, we shall always assume 
		\begin{align}
		\nabla \vphi \cdot \bm{n}  = 0  \quad  &\textrm{on }  (0, T) \times \partial \Omega,  \label{bc:contac_angle}\\ 
		\nabla \mu \cdot \bm{n}  = 0  \quad  \textrm{and} \quad 
		\nabla p \cdot \bm{n} = 0  \quad  &\textrm{on }  (0, T) \times \partial \Omega. \label{bc:mass_conservation}
	\end{align}
	Moreover, from a physical point of view, it is desirable to impose mixed boundary conditions of the form 
	\begin{align}
		\bu &= \bm{0}  
		 &\textrm{on }  (0, T) \times \Gamma_D,\\
		\bm{\sigma}\,  \bm{n} &= \bm{g}  \quad  &\textrm{on }  (0, T) \times \Gamma_N,  \label{bc:discplacement_neumann}
	\end{align}  
	where we split the boundary into a Dirichlet part $\Gamma_D$ and a Neumann part $\Gamma_N$.
	However, depending on the geometry of the domain, this would preclude us from deducing necessary regularity for the system.
	A detailed discussion of possible geometric conditions that allow for more general boundary conditions can be found in Remark \ref{rem:bc}. There, we also mention other types of boundary conditions along with the respective regularity results. 
\end{subequations}

\addtocontents{toc}{\SkipTocEntry}
\subsection{Main results}
The goal of this paper is to establish local-in-time well posedness of the original Cahn--Hilliard--Biot system and its visco-elastic extension. In both cases our strategy will involve rewriting the problem as a fixed-point equation by eliminating dispensable variables and linearizing highest order operators, such that the application of a contraction principle will provide a unique solution. Central to this approach is showing that the linearized system admits maximal regularity, along with Lipschitz estimates for the corresponding right-hand sides in suitable spaces. \\ 
A novelty in this exposition is that unlike other treatments of the Biot model, c.f. e.g. \cite{SHOWALTER2000310}, we utilize the equivalent formulation \eqref{eq:biot_2}, as opposed to \eqref{eq:biot}, for the poro-(visco-)elastic subsystem, thereby avoiding difficulties arising from the phase-field dependent material parameters $M$ and $\alpha$ by application of chain rule. 
\par 
\medskip 
In absence of the visco-elastic regularization $\nabla \cdot ( \C_\nu \E(\pt \bm{u}))$, we will show the following theorem.

	\begin{theorem}[\new{Elastic case}]\label{thm:existence_strong_nve}
		Suppose $\vr =0$ and let $s = 4r$ with $r> 4$ if $n = 3$ and $r > \frac{12}{5}$ if $n \leq 2$. Under the assumptions \ref{II.A:domain}-\ref{II.A:source_terms}, and for initital conditions 
		\begin{align*}
			\vphi_0 &\in  ((W^{1, 3}(\Omega))', W^{3, \frac{3}{2}}_{N}(\Omega))_{1- \frac{1}{r}, r}
			\cap W^{2, p^*}(\Omega) \textnormal{ for some } p^* > \tfrac{3}{2} \textnormal{ if } n = 3 \textnormal{ and } p^* \geq \tfrac{3}{2} \textnormal{ if } n \leq 2, \\
			\theta_0 &\in ((H^1(\Omega))', H^1(\Omega))_{1- \frac{1}{s}, s}, 
		\end{align*}
		 there exists a (possibly small) $T > 0$ such that the system \eqref{eq:strong_formulation} has a unique solution 
		\begin{subequations}
			\begin{align}
				\vphi &\in W^{1, r} (0, T; (W^{1, 3}\new{(\Omega)})') \cap L^r (0, T; W^{3,\sfrac{3}{2}}_N (\Omega)),  \\ 
				\bm{u} & \in L^s (0, T; \bm{H}^2_{\Gamma_D} (\Omega)), \\ 
				\theta & \in W^{1, s} (0, T; (H^1 (\Omega))') \cap L^{s} (0, T;H^1(\Omega)). 
			\end{align}
		\end{subequations}
	\end{theorem}

After using elliptic theory on \eqref{eq:strong_formulation_u} to find an identity for the displacement $\bm{u}$ and inserting this expression, along with \eqref{eq:strong_formulation_p} for the pressure $p$, into \eqref{eq:strong_formulation_theta}, we discover that, to the highest order, the evolution of the fluid is governed by the sum of two operators. Since establishing maximal regularity for each operator separately is not sufficient, we will show that the aggregate corresponds to a dissipative operator on an appropriate Hilbert space and deduce that it generates an analytic semigroup, thus obtaining the desired property. 
\par 
\medskip

Moreover, we will show the following theorem for the visco-elastic Cahn--Hilliard--Biot system.
	
	\begin{theorem}[Visco-elastic case]\label{thm:existence_strong}
		Suppose $\vr =1$ and let $q = 2p$ with $q > 3$ if $ n = 3$ and $q \geq 3$ if $n = 2$. 
		\new{
			Under the assumptions \ref{A:domain}-\ref{A:source_terms} there exists some $r > 1$ such that for initial conditions
			\begin{align*}
				&\vphi_0 \in W^{1, q}(\Omega) \cap W^{2, p}(\Omega), \\
				& \theta_0 \in  \big(W^{-2, q}_N(\Omega), L^q(\Omega)\big)_{1- \frac{1}{3r}, 3r}, \\ 
				& \bm{u}_0 \in \bm{W}^{1, q}_{\Gamma_D}(\Omega),
			\end{align*}
		}
		and a (possibly small) $T > 0$, the system \eqref{eq:strong_formulation} has a unique solution 
		\begin{subequations}
		\begin{align}
				\vphi &\in W^{1, r} (0, T; W^{-2, p}_N (\Omega)) \cap L^r (0, T; W^{2, p}_N (\Omega)),  \\ 
				\bm{u} & \in W^{1, r} (0, T; \X (\Omega)), \\ 
				\theta & \in W^{1, 3r} (0, T; W^{-2, q}_N (\Omega)) \cap L^{3r} (0, T; L^q(\Omega)). 
		\end{align}
		\end{subequations}
	\end{theorem}
For details on the necessary conditions on $r$ we refer to Section~\ref{sec:contraction}.\par
Note that our proof relies on $\vphi(t, \cdot) \in C^\beta(\overline{\Omega})$, $t \in [0, T]$,  for some $\beta >0$, thus guaranteeing invertibility of the involved elliptic operators in suitable spaces, see \ref{A:iso} and Remark \ref{rem:bc}, and that the spatial regularity postulated above is minimal in this regard. \\ 
Consequently, maximal regularity in the appropriate spaces does not follow immediately from classical theory, and will be shown by establishing $\mathcal{BIP}$ \new{("bounded imaginary powers", see Section \ref{sec:preliminaries} for details)} of the respective operators, such that interpolation and the application of Banach scales yields the desired property.

\begin{remark}
		 \textit{(i)} Under more restrictive assumptions, one could obtain better regularity for the displacement $\bm{u}$. Then, by similar arguments to those below, strong and even classical solutions can be retrieved. However, due to a counterexample by Shamir \cite{shamir1968regularization}, this is certainly not possible if we want to impose mixed boundary conditions where the segments with different conditions are allowed to meet arbitrarily. \\ 
		 \new{
		 	\textit{(ii)} The local-in-time well-posedness results presented in this work naturally raise the question about global-in-time existence of solutions. Such results exist for the Navier-Stokes equation, where the energy is convex and the data is chosen sufficiently close to the global minimizer, which is zero, i.e., the data is "small". However, since the free energy of the Cahn--Hilliard--Biot system is non-convex, there is in general not a unique global minimizer and one would expect that global-in-time solutions exist for data which are sufficiently close to a global or local minimum. A possible proof of this conjecture might be based on the \L ojasiewicz-Simon inequality, but is beyond the scope of this exposition. We refer to \cite{MR2571478} where this strategy is employed on a Navier--Stokes/Cahn--Hilliard system.  
		}
\end{remark}

\addtocontents{toc}{\SkipTocEntry}
\subsection{Structure of this paper}

We start in Section \ref{sec:preliminaries} by introducing the notation and reviewing several prerequisites on pointwise multiplication and interpolation along with results related to maximal $L^r$-regularity and semigroup theory. \\
Corresponding to the distinction between the visco-elastic and purely elastic case, the remaining portion of this paper is split into two parts, where Part~\ref{part:nve} is focused on the proof of Theorem~\ref{thm:existence_strong_nve} and Part~\ref{part:ve} deals with Theorem~\ref{thm:existence_strong}. \\ 
In this regard, we begin Section~\ref{sec:assumptions_non_ve} with the assumptions for the fully elastic case and proceed in Section~\ref{sec:proof_nve} with the proof of Theorem~\ref{thm:existence_strong_nve}, moving all arguments pertaining maximal regularity and the contraction estimates to Section \ref{sec:proof_nve} and Section~\ref{sec:contraction_non_ve}, respectively. \\ 
Similarly, the second part contains the assumptions for the visco-elastic case (Section~\ref{sec:assumptions}) and the proof of Theorem~\ref{thm:existence_strong} (Section~\ref{sec:proof}) with the argumentation for maximal regularity and Lipschitz estimates being moved to Section \ref{section:inverse} and Section \ref{sec:contraction}, respectively.

\medspace

\section{Notation and preliminaries} \label{sec:preliminaries}

This section begins with notation of relevant function spaces and a theorem on the compositions of classically differentiable functions with Sobolev functions. We proceed by recalling important results on pointwise multiplication, interpolation theory and maximal $L^p$-regularity. \\ 
Thereafter, we recall the notions of $\mathcal{H}^\infty$-calculus and $\mathcal{BIP}$ along with $\mathcal{R}$-sectorial operators in greater detail. The tools to establish these properties for given operators are provided by Theorems \ref{thm:H_infty_calculus}-\ref{thm:r_sec_max_reg}, which also connect the concepts to maximal regularity theory. Related to these ideas are fractional power scales and interpolation-extrapolations scales, that, combined with Theorems \ref{thm:H_infty_calculus}-\ref{thm:r_sec_max_reg}, will allow us to establish maximal regularity for an even wider range of operators. \par
\medskip 
For brevity we write $\OT$ instead of $(0, T) \times \Omega$ for any $T > 0$. 

\addtocontents{toc}{\SkipTocEntry}
\subsection*{Function spaces}

For $p \in [1, \infty]$, we denote by $W^{1, p}(\Omega)$ the usual Sobolev spaces and use the abbreviations
\begin{equation*}
	H^1(\Omega) \coloneqq W^{1, 2}(\Omega), \quad H^2_{N}(\Omega) \coloneqq \{ f \in W^{2, 2}(\Omega):  \nabla f \cdot \bm{n}  = 0 \quad \textrm{on } \partial \Omega\}, 
\end{equation*}
where $\bm{n}$ is the outer normal vector on $\partial \Omega$. Vector-valued functions and their corresponding function spaces are indicated by bold fonts, i.e., $\bm{f} \in \bm{W}^{1, p}(\Omega)$. Moreover, we use 
\begin{equation*}
	W^{-1, p}(\Omega) = (W^{1, p'}_0(\Omega))', 
\end{equation*}
where $W^{1, p'}_0(\Omega) \coloneqq \{ f \in W^{1, p'}(\Omega) \colon f_{|\partial \Omega} = 0 \}$. For $\Gamma_D \subseteq \partial \Omega$, let $\X (\Omega)$ be the subspace of $\bm{W}^{1, q}(\Omega)$, such that the trace on $\Gamma_D$ vanishes
\begin{equation*}
	\X(\Omega) \coloneqq \{\bm{v} \in \X \colon \bm{v}_{| \Gamma_D} = 0 \}, 
\end{equation*}
and 
\begin{equation*}
	\Xd(\Omega) \coloneqq  (\Xp_{\Gamma_D}(\Omega))'.
\end{equation*}
\par

The following theorem is a fundamental result on the composition of classically differentiable functions and Sobolov functions of certain order. The formulation below is taken from \cite[Thm. 3]{abels2021local}. 

\begin{theorem}[Composition with Sobolev functions]\label{thm:compositon}
	Let $\Omega \subset \R^n$ be a bounded domain with boundary of class $C^1$, $m, N \in \N$ and let $1 \leq p < \infty$ such that $m - n/p > 0$. Then, for every $f \in C^m(\R^N)$ and every $R > 0$ there exists a constant $C > 0$ such that for all $\bm{\vphi} \in \bm{W}^{m, p} (\Omega)$ with $\norm{\bm{\vphi}}_{\bm{W}^{m, p}}  \leq R$, we have $f(\bm{\vphi}) \in W^{m, p} (\Omega)$ and $\norm{f(\bm{\vphi})}_{\bm{W}^{m, p}} \leq C$. Moreover, if $f \in C^{m+1}(\R^N)$, then for all $R > 0$ there exists a constant $L > 0$ such that 
	\begin{equation*}
		\norm{f(\bm{\vphi}_1) - f(\bm{\vphi}_2)}_{W^{m, p}} \leq L \norm{\bm{\vphi}_1 - \bm{\vphi}_2}_{\bm{W}^{m, p}}
	\end{equation*} 
	for all $\bm{\vphi}_i \in \bm{W}^{m, p}(\Omega)$ with $\norm{\bm{\vphi}_i}_{W^{m, p}} \leq R$, $i = 1, 2$.  
  
  \end{theorem}

\begin{proof}
	The first assertion of the theorem is a consequence of \cite[Chapter 5, Theorem 1 and Lemma]{runst2011sobolev}. The second assertion follows easily from the first in conjunction with the mean value theorem. 
\end{proof}

\addtocontents{toc}{\SkipTocEntry}
\subsection*{Multiplication results}

We proceed by presenting results concerning the pointwise multiplication of some distributions with a Sobolev function. Note that, though the references given below only show the claims on the whole space, the boundary of our domain $\Omega$ is regular enough that the results can easily be transferred by means of a standard extension operator, e.g. Stein's extension operator. 
\\ 
For $2p = q > n$, we find by computation along with \cite[4.~4.~Theorem 2, Lemma p. 175]{runst2011sobolev} that the following multiplication operations are continuous
\begin{gather}
	W^{-1, p}(\Omega) \times W^{1, q}(\Omega) \rightarrow  W^{-1, p}(\Omega) \label{mult:-1p}, 
	\quad (f, g) \rightarrow \langle f, g\, \cdot \rangle \\ 
	W^{-1, q}(\Omega) \times W^{1, q}(\Omega) \rightarrow  W^{-1, q}(\Omega), 
	\quad (f, g) \rightarrow \langle f, g\, \cdot \rangle\label{mult:-1q}. 
\end{gather}
Moreover, since it holds for any $\delta  > 0$ that $1+ \delta > \frac{n}{q}$ if $q > n$, we deduce 
\begin{equation} \label{mult:triple}
		W^{2\delta, q}(\Omega) \times W^{- \delta, q}(\Omega) \rightarrow  W^{-1, q}(\Omega), 
		\quad (g, f) \rightarrow \langle f, g\, \cdot \rangle. 
\end{equation}
from \cite[Thm. A.1]{MR4339668}. By the assumptions $q > n$ and $2p= q \geq 3$ along with \cite[Thm. A.1]{MR4339668}, we can find some sufficiently small $\tilde{\delta} > 0$ such that 
\begin{equation}\label{eq:mutl_delta_p}
	W^{1, p}(\Omega) \times W^{- \tilde{\delta}, p}(\Omega) \rightarrow  W^{-1, p}(\Omega), 
	\quad (g, f) \rightarrow \langle f, g\, \cdot \rangle. 
\end{equation}

\addtocontents{toc}{\SkipTocEntry}
\subsection*{Interpolation results}

Denote by $(X_0, X_1)_{\vartheta, p}$ the real interpolation space between the Banach spaces $X_0$ and $X_1$ with exponent $\vartheta$ and summation index $p$ and let $X_1 \hookrightarrow X_0$ be a densely injected Banach couple, $r \in (1, \infty)$ and $T > 0$. Then, we have the following continuous embedding, cf. \cite[Ch. III, Thm. 4.10.2]{amann1995linear}, 
\begin{equation}\label{int:BUC}
		L^r(0, T; X_1) \cap W^{1, r} (0, T; X_0) \hookrightarrow BUC([0, T], (X_0, X_1)_{1- \frac{1}{r}, r})
\end{equation}

Furthermore, we use the following well-known result, \cite[Lem.~1]{abels2021local}: 

\begin{lemma}\label{lemma:interpolation_hoelder}
	Let $X_0 \subseteq Y \subseteq X_1$ be a Banach space such that 
	\begin{equation*}
			\norm{x}_{Y} \leq C \norm{x}_{X_0}^{1-\vartheta} \norm{x}_{X_1}^\vartheta
	\end{equation*}
	for every $x \in X_0$ and a constant $C > 0$, where $\vartheta \in (0, 1)$. Then, 
	\begin{equation*}
		C^{0, \beta} ([0, T]; X_1) \cap L^\infty (0, T; X_0) \hookrightarrow C^{0, \beta \vartheta} ([0, T]; Y)
	\end{equation*}
	continuously. 
\end{lemma}

\addtocontents{toc}{\SkipTocEntry}
\subsection*{Maximal $L^p$-regularity}

The following definition is taken from \cite{ARENDT20071}. Here, we assume $D, Y$ to be Banach spaces such that the embedding $D \xhookrightarrow{d} Y$ is dense. 

\begin{definition}
	An operator $A \in \mathcal{L}(D, Y)$ is said to have $L^p$-\textit{maximal regularity} for $p \in (1, \infty)$ if for some interval $(a, b)$ and all $f \in L^p(a, b; Y)$ there exists a unique $u \in W^{1, p}(a, b; Y) \cap L^p(a, b; D)$ satisfying
	\begin{equation*}
		\pt u + Au = f \quad \textit{a.e. on } (a, b), \quad u(a) = 0. 
	\end{equation*}
\end{definition}

\begin{remark}
		If $A$ has $L^p$-maximal regularity for some $p \in (1, \infty)$ on a given interval $(a, b)$, then $A$ already has $L^p$-maximal regularity for \textit{all} bounded intervals and \textit{all} $p \in (1, \infty)$, cf. \cite{ARENDT20071}, and we write $A \in \mathcal{MR}$. 
\end{remark}

\addtocontents{toc}{\SkipTocEntry}
\subsection*{$\mathcal{H}^\infty$-calculus, $\mathcal{BIP}$ and analytic semigroups}\label{subsec:BIP}

In connection with maximal regularity, we will also encounter the notions of sectorial operators, bounded imaginary powers ($\mathcal{BIP}$) and analytic semigroups, which are the subject of this section. This introduction in based on  \cite{MR2047641, Weis2006, PrussSimonett2016}. \\ 
In the following, we denote the domain and range of some operator $A$ by $D(A)$, $R(A)$, respectively, and use $\sigma(A), \rho(A)$ as an abbreviation for the spectrum and the resolvent set, respectively. \\ 
A closed linear operator $A \colon D(A) \rightarrow X$ in some Banach space $X$ is \textit{sectorial}, if 
\begin{enumerate}[label*=(P.\roman*), noitemsep]
	\item $\overline{D(A)} = X = \overline{R(A)}$ and $(0, \infty) \subset \rho(-A)$; 
	\item \label{eq:resolvent_estimates} $\norm{t (t + A)^{-1}}_{\mathcal{L}(X)} \leq M $ for all $t > 0$ and some $M < \infty $. 
\end{enumerate}
Then, there exists some $\theta > 0$ such that $\sup\, \{ \norm{\lambda (\lambda + A)^{-1}}_{\mathcal{L}(X)}  \colon \abs{\arg \lambda} < \theta \} < \infty$ and $\rho(- A) \supset \Sigma_{\theta}$, where the sector $ \Sigma_{\theta}$ with vertex at the origin and opening angle $2 \theta$ is defined as  
\begin{equation*}
	\Sigma_{\theta} \coloneqq \{  \lambda \in \mathbb{C} \setminus \{ 0\}  \colon \abs{\arg \lambda} < \theta \}. 
\end{equation*}
The \textit{spectral angle} $\phi_A$ of the sectorial operator $A$ is the infimum 
\begin{equation}
	\phi_A \coloneqq  \textnormal{inf}\, \Big\{ \phi \colon \rho(- A) \supset \Sigma_{\pi - \phi}, \sup_{\lambda \in  \Sigma_{\pi - \phi}} \norm{\lambda (\lambda + A)^{-1}}_{\mathcal{L}(X)}   < \infty \Big\}. 
\end{equation}

In particular, it holds that $\phi_A \in [0, \pi )$ and $\phi_A \geq \sup\, \{  \abs{\arg \lambda} \colon \lambda \in \sigma(A) \}$. \\ 
We proceed by defining 
\begin{align*}
	\mathcal{H}^{\infty} (\Sigma_\phi) &\coloneqq \{  f \colon \Sigma_\phi \rightarrow \mathbb{C} \textnormal{ holomorphic and bounded}\}
\end{align*}
for any $\phi \in (0, \pi]$ and further set 
\begin{equation*}
	\mathcal{H}_0 (\Sigma_\phi) \coloneqq \Big\{ f \in  \mathcal{H}^{\infty} (\Sigma_\phi)  
	\colon \abs{f(z)} \leq C \frac{\abs{z}^s}{(1+\abs{z})^{2s} }   \textnormal{ for all } z \in \Sigma_\phi \textnormal{ and some }  C, s > 0
	\Big\} . 
\end{equation*}
For a fixed sectorial operator $A$ and a given angle $\phi \in (\phi_A, \pi)$ let 
\begin{equation} \label{def:dunford_calculus}
	f (A) \coloneqq \frac{1}{2\pi i } 
	\int_{\gamma(t)} f(\lambda) (\lambda - A)^{-1} \, d\lambda 
	\quad \textnormal{for any } f \in \mathcal{H}_0(\Sigma_\phi), 
\end{equation}
where
\begin{equation*}
	\gamma(t) = \begin{cases}
		 t e^{- i \theta} \quad &\textnormal{if } t < 0,\\ 
		 t e^{ i \theta} \quad &\textnormal{if } t \geq  0, 
	\end{cases}
\end{equation*}
for some $\phi_A < \theta < \pi $ and note that this definition is independent of $\theta$. The integral is absolutely convergent and the map $f \mapsto f(A)$ describes an algebra homomorphism $\mathcal{H}_0 (\Sigma_\phi) \rightarrow \mathcal{L}(X)$. 

\begin{definition}[Operator semigroup, {\cite[Def. 3.3.1]{PrussSimonett2016}, \cite[Def.~5.1]{pazy2012semigroups}}]
	A family of operators $\{  T (t) \}_{t \geq 0} \subset{\mathcal{L}(X)}$ is called a \emph{semigroup} if 
	\begin{equation}\label{semigroup_property}
		T(t + s) = T(t) T(s), \quad s,t > 0, \quad T(0) = I, 
	\end{equation}
	is satisfied. The semigroup is \emph{strongly continuous}, if in addition
	\begin{equation}\label{strongly_continuous_SG}
		\lim_{t \searrow 0} T(t) x = x, \quad x \in X, 
	\end{equation}
	holds. \\ 
	If there exists a sector $\Sigma_\theta$ such that the map $z \mapsto T(z)$ is analytic in $\Sigma_\theta$ and \eqref{semigroup_property}, \eqref{strongly_continuous_SG} still hold for arbitrary $z_1, z_2 \in \Sigma_\theta$ and $z \rightarrow 0$ in $\Sigma_\theta$, respectively, then $(T(z))_{z \in \Sigma_\theta}$ is called \emph{strongly continuous analytic semigroup}. 
\end{definition}

Given a sectorial operator $A$ with spectral angle $\phi_A < \frac{\pi}{2}$, the calculus defined in \eqref{def:dunford_calculus} gives rise to the operator semigroup $\{  e^{-t A} : t > 0 \}$ with \emph{generator} $-A$. Moreover, we have the following perturbation theorem.

\begin{theorem}[{\cite[Thm. 1.3.1 ]{amann1995linear}}] \label{thm:perturbation_semigroup}
	Suppose $-A : D(A) \subset X \rightarrow X$ generates a strongly continuous analytic semigroup and $B \in \mathcal{L}(D(A), X)$ such that 
	\begin{equation*}
		\norm{B x}_{X} \leq C \norm{x}_{X} \quad x \in D(A). 
	\end{equation*} 
	Then $-(A + B)$ also generators a strongly continuous analytic semigroup.  
\end{theorem}

In particular, due to the inclusion $D(A) \hookrightarrow X$, it follows that shifting the generator of a strongly continuous analytic semigroup by a scalar multiples of the identity retains this property.  \par 
\medskip
De Simon \cite{deSimon} was able to show that if $-A$ generates an analytic semigroup $T$ of negative exponential type, i.e. $\norm{T (t)} \leq M e^{- \mu t }$, on a Hilbert space, then $A$ already has maximal $L^p$ regularity on $\R_+$. The latter condition can be neglected if we restrict ourselves to bounded intervals, leading to the following theorem.

\begin{theorem}[Maximal regularity of generators of analytic semigroups] \label{thm:max_reg_analytic_hilbert}
	Let $- A : D(A) \subset H \rightarrow H$ be a densely defined operator on the Hilbert space H, generating a strongly continuous analytic semigroup. Then $A$ has maximal $L^p$ regularity, $p \in (1, \infty)$, on bounded intervals. 
\end{theorem}

\begin{definition}[$\mathcal{H}^\infty$-calculus, McIntosh '86 \cite{mcintosh1986operators}]
	A sectorial operator $A$ with spectral angle $\phi_A$ admits a bounded $\mathcal{H}^\infty$-calculus if there is some $\phi > \phi_A$ and a constant $K_{\phi} < \infty$ such that 
	\begin{equation}\label{eq:h_infty-calculus}
		\norm{f(A)}_{\mathcal{L}(X)} \leq K_\phi \norm{f}_{C^0(\Sigma_\phi)} \textnormal{ for all } f \in \mathcal{H}_0 (\Sigma_\phi)
	\end{equation}
	In this case, we write $A \in \mathcal{H}^\infty(X)$ and call 
	\begin{equation*}
		 \phi_A^\infty \coloneqq \textnormal{inf}\, \{ \phi > \phi_A  \colon \eqref{eq:h_infty-calculus} \textnormal{ is valid}\}
	\end{equation*}
	the $\mathcal{H}^\infty$-angle of $A$. 
\end{definition}

\begin{remark}
	If \eqref{eq:h_infty-calculus} holds for some $\phi > \phi_A$, then there is a unique extension of the functional calculus to $\mathcal{H}^\infty (\Sigma_\phi)$, i.e., $f (A)$ is well-defined for all $f \in \mathcal{H}^\infty(\Sigma_\phi)$ and the estimate holds. 
\end{remark}

Furthermore, we say that a sectorial operator $A$ admits \textit{bounded imaginary powers} in $X$, $A \in \mathcal{BIP}(X)$, if $A^{i s} \in \mathcal{L}(X)$ for all $s \in \R$ and there exists some constant $C > 0$ such that $\norm{A^{i s} }_{ \mathcal{L}(X)} \leq C$ for all $\abs{s} \leq 1$. We observe that the functions $f_s$ with $z \mapsto z^{is}$ belong to $\mathcal{H}^{\infty} (\Sigma_\phi)$ for any $s \in \R$ and $\phi \in (0, \pi)$, and deduce immediately 
\begin{equation*}
	\mathcal{H}^\infty(X) \subset \mathcal{BIP}(X). 
\end{equation*} 
For a given sectorial operator $A \in \mathcal{BIP}(X)$ we define the \textit{power angle} $\theta_A$ as 
\begin{equation*}
	\theta_A  \coloneq \limsup_{\abs{s} \rightarrow \infty} \frac{1}{\abs{s}} \log \norm{A^{is}}_{\mathcal{L}(X)}.
\end{equation*}

	It is a well-known fact, cf. \cite[p. 177]{amann1995linear}, that for any $A \in \mathcal{BIP}(X)$ with power angle $\theta_A < \frac{\pi}{2}$, the operator $-A$ generates a strongly continuous analytic semigroup. These notions give rise to a sufficient condition for $A \in \mathcal{MR}$.  
	
	\begin{theorem}[{\cite[Thm. 4.8.10]{amann1995linear}}] \label{thm:max_reg_semigroups}
		Let $ 0 < T < \infty$, $X$ be a Banach space and $A : D(A) \subset X \rightarrow X$. Suppose that there exists some $\omega > 0$ such that $A+ \omega \in \mathcal{BIP}(X)$ with power angle $\theta_{A+ \omega}  < \frac{\pi}{2}$. Then the Cauchy problem
		\begin{equation*}
			\pt \vphi + A \vphi = f(t), \quad t \in (0, T), \quad \vphi(0) = 0, 
		\end{equation*}
		where $f \in L^r(0, T;  X)$, has maximal regularity of type $L^r$ on $(0, T)$. 
	\end{theorem}

Let us remark that Amann further requires that $-A$ generates an analytic semigroup. However, this requirement is actually redundant, as it is already follows from $A+ \omega \in \mathcal{BIP}(X)$ being the generator of an analytic semigproup and Theorem \ref{thm:perturbation_semigroup}.

\par 
\medskip 
To see some examples for operators which are in these classes, let $m \in \N$ and consider differential operators of the form $\mathcal{A}(x, D) = \sum_{\abs{\alpha} \leq 2m} a_{\alpha}(x) D^{\alpha}$ of order $2m$ with $\mathcal{L}(E)$-valued coefficients $a_{\alpha}(x)$ where $D = i (\partial_1, \ldots, \partial_n)$. While, in general, $E$ may be an arbitrary Banach space, we are only interested in $E = \new{\mathbb{C}}$ and therefore restrict the following exposition to this case, i.e., we assume that the coefficient functions $a_\alpha$ are \new{complex}-valued. To account for general boundary conditions, let $\mathcal{B}_j (x, D) = \sum_{\abs{\beta} \leq m_j} b_{j \beta(x)} D^{\beta}$, $m_j < 2m$, $j = 1, \ldots, m$, be boundary operators with \new{complex}-valued coefficients $ b_{j \beta(x)}$. Keep in mind that to obtain a well-posed problem, $\mathcal{A}$ and $\mathcal{B}_j$ need to be compatible, which is the case if, e.g., the Lopatinskii-Shapiro condition, cf. \cite[Sec. 8.1]{MR2006641}, is satisfied. \\ 
The principal symbol of our differential operator is defined as 
\begin{equation*}
	\mathcal{A}_{\#} (x, \xi) \coloneqq  \sum_{\abs{\alpha} = 2m} a_{\alpha}(x) \xi^{\alpha}, \quad \xi \in \R^n,
\end{equation*} 
and, in the case of \new{complex}-valued coefficient, we say that $\mathcal{A}_{\#}$ is \textit{parameter elliptic} in $x$, if there is an angle $\phi \in [0, \pi)$ such that 
\begin{equation}\label{def:parameter_elliptic}
		\sigma (\mathcal{A}_{\#} (x, \xi))  \subset \Sigma_\phi \textnormal{ for all } \xi \in \R^n, \quad \abs{\xi} = 1. 
\end{equation}
The infimum over all angles such that this condition is satisfied is called \textit{angle of ellipticity} $\phi_{\mathcal{A}_{\#}(x)}$
\begin{equation*}
	\phi_{\mathcal{A}_{\#}(x)} \coloneqq \inf \{ \phi \in [0, \pi)  \colon \eqref{def:parameter_elliptic} \textnormal{ is valid}\} = 
	\sup_{\abs{\xi} = 1}  \{ \arg \mathcal{A}(x, \xi)\}. 
\end{equation*}

\begin{theorem}[{\cite[[Thm. 2.3]{MR2047641}}] \label{thm:H_infty_calculus}
	Let $1 < p < \infty$, and $n \in \N$ and assume that $\Omega$ is a bounded domain in $\R^{n+1}$ with compact $C^{2m}$-boundary $\partial \Omega$. Let $(\mathcal{A}, \mathcal{B}_1,\ldots, \mathcal{B}_m)$ be as above, satisfying the Lopatinskii-Shapiro condition, and assume that there exists some $\phi_\mathcal{A} \in [0, \pi)$ such that for all $x \in \overline{\Omega}$ the principal symbol $\mathcal{A}_{\#} (x, \xi)$ is parameter elliptic with $\phi_{\mathcal{A}_{\#}(x)} \leq \phi_\mathcal{A}$.  Moreover, suppose that 
	\begin{enumerate}[label= (\roman*), noitemsep]
		\item $a_\alpha \in C^\rho(\overline{\Omega})$ for some $\rho> 0$ and each $\abs{\alpha} = 2m$; 
		\item $a_\alpha \in [L^\infty + L^{r_k}](\Omega)$ for each $\abs{\alpha} = k < 2m$ with $r_k \geq p $ and $2m - k > \sfrac{n}{r_k}$; 
		\item $b_{j \beta} \in C^{2m-mj}(\partial \Omega)$ for each $j, \beta$. 
	\end{enumerate}
	Let $A_B$ denote the realization of $\mathcal{A}(x, D)$ in $X = L^p(\Omega)$ with domain
	\begin{equation*}
		D(A_B) = \{u \in H^{2m}_p (\Omega; E) : \mathcal{B}_j (x, D)u = 0 , \ j = 1, \ldots, m\}. 
	\end{equation*}
	Then, for each $\phi > \phi_{\mathcal{A}}$ there is $\mu_{\phi} \geq 0$ such that $\mu_\phi + A_B \in \mathcal{H}^\infty(L^p(\Omega))$ with $\phi^\infty_{\mu_\phi + A_B} \leq \phi_{\mathcal{A}}$.  
\end{theorem}

\new{As we will only use $\mathcal{R}$-sectoriality as a criterion for maximal regularity, we refrain from providing a precise definition, but refer the interested reader to \cite{MR2006641}.} In general, verifying that a given operator is $\mathcal{R}$-sectorial is a formidable task and requires quite subtle arguments. Fortunately, Clémont and Prüss \cite{MR1816437}, established the following fundamental relation between $\mathcal{R}$-bounded operators and the class $\mathcal{BIP}(X)$. 

\begin{theorem}\label{thm:bip->r_sectorial}
	Let $A \in \mathcal{BIP}(X)$ with power angle $\theta_A$ on a Banach space $X$ of class $\mathcal{HT}$. Then $A$ is $\mathcal{R}$-sectorial and $\phi^{\mathcal{R}}_A  < \theta_A$. 
\end{theorem}

Note that a Banach space $X$ is of class $\mathcal{HT}$, if the Hilbert transform is bounded in $L^p(\R, X)$ for some $p \in (1, \infty)$.\\ 
Finally, we can state the most important theorem of this section - a criterion for maximal regularity based on the notion of $\mathcal{R}$-sectoriality.  

\begin{theorem}[{\cite[Thm. 4.4]{MR2006641}}] \label{thm:r_sec_max_reg}
	Suppose $A$ is a sectorial operator with spectral angle $\phi_A <  \tfrac{\pi}{2}$ in a Banach space $X$ of class $\mathcal{HT}$. Then the Cauchy-problem
	\begin{equation*}
		\pt \vphi + A \vphi = f(t), \quad t \geq 0, \quad \vphi(0) = 0, 
	\end{equation*}
	where $f \in L^r(\R_+;  \new{X}))$, has maximal regularity of type $L^r$ on $\R_+$ if and only if $A$ is $\mathcal{R}$-sectorial with $	\phi^{\mathcal{R}}_A < \tfrac{\pi}{2}$. 
\end{theorem}

\new{A similar result characterizing maximal regularity in terms of $\mathcal{R}$-boundedness of the resolvents in UMD-spaces can be found in \cite[Thm.\ 4.2]{MR1825406}}

\addtocontents{toc}{\SkipTocEntry}
\subsection*{Fractional power scales and interpolation-extrapolation scales} \label{subsec:scales}

\begin{definition}[Banach scale]
Let $\mathsf{A} \in \{\N, \mathbb{Z}, \R^+, \R\}$ be an index set and suppose $E_\alpha \coloneqq (E_\alpha , \norm{\cdot}_{\alpha})$ is a Banach space for each $\alpha \in \mathsf{A}$ 
such that, $E_\alpha \hookrightarrow E_\beta$ for each pair $\alpha ,\beta \in \mathsf{A}$ with $\alpha > \beta $. Moreover, assume that there are isomorphic operators $A_\alpha \in \mathcal{L}_{is} (E_{\alpha+1}, E_\alpha)$ and that the diagram
\begin{center}
	\begin{tikzcd}
		\cdots \arrow[r, hook ]
		&E_{\alpha + 1} \arrow[r, hook] \arrow[d, "A_\alpha" right  , "\cong" left] 
		& E_{\beta + 1} \arrow[d, "A_\beta" right, "\cong" left ] \arrow[r, hook ]
		& \cdots 
		\\
		\cdots \arrow[r, hook ]
		&E_{\alpha} \arrow[r, hook] 
		&E_{\beta}\arrow[r, hook ] 
		& \cdots  
		\end{tikzcd}\\
\end{center}
	is commutative. Then $[(E_\alpha, A_\alpha)\  ; \ \alpha \in \mathsf{A} ]$ is a Banach scale over $\mathsf{A}$. 
\end{definition}

While the following construction can be done under less restrictive assumptions, these suffice for our purposes. The more general theory can be found in \cite[Chapter V]{amann1995linear}. \par 
Suppose that $A \in \mathcal{BIP}(E)$ and note that the fractional powers, $A^\alpha : D(A^\alpha ) \rightarrow E$, for $A$ are well-defined for $\alpha \in \R$. Therefore, we can define 
\begin{equation*}
	E_\alpha \coloneqq (E_\alpha, \norm{\cdot}_{\alpha}) \coloneqq (D(A^\alpha), \norm{A^\alpha \cdot}), \quad \alpha \in \R^+
\end{equation*} 
and set $A_\alpha$ as the $E_\alpha$-realization of $A$, i.e., 
\begin{equation*}
	D(A^\alpha) \coloneqq \{ x \in E_\alpha \cap D(A) \colon A x \in E_\alpha \} , \quad A_{\alpha} x \coloneqq  A x, 
	\quad \alpha \in \R^+. 
\end{equation*}
We call $[(E_\alpha, A_\alpha)\ ; \ \alpha \in \R^+]$ the \textit{one-sided fractional power scale} generated by $(E, A)$ and \cite[Thm. V 1.2.4]{amann1995linear} tells us that it is in fact a Banach scale. \\ 
Since $A^{-1}$ is well-defined, we can also endow $E$ with the norm $\norm{A^{-1} \cdot}$, which is weaker then $\norm{\cdot}$, and obtain the \textit{extrapolation space of $E$ generated by $A$}
\begin{equation*}
	(E_{-1}, \norm{\cdot }_{-1}) \coloneqq \overline{E}^{\norm{A^{-1} (\cdot )}}
\end{equation*}
by completion. Since $A$ is densely defined, it holds that the closure of $E_1$ in $E$ is again $E$, and we note that there exists a unique continuous extension $A_{-1}$ of $A$ such that $A_{-1}$ is an isometric isomorphism from $E$ onto $E_{-1}$, which we call the \textit{$E_{-1}$-realization of $A$}. \\ 
By \cite[Thm. V 1.3.7]{amann1995linear} the linear operator $A_{-1}$ is of class $\mathcal{BIP}(E_{-1})$, such that we can inductively define
\begin{align*}
	E_{-k} &\coloneqq (E_{_k}, \norm{\cdot }_{-k}) \coloneqq \overline{E_{-k+1}}^{\norm{(A_{-k+1} )^{-1}(\cdot)} }, \\ 
	A_{-k} &\coloneqq \textnormal{closure of } A_{-k+1} \textnormal{ in } E_{-k}, 
\end{align*}
for all $1 \leq k \leq m$. By the considerations above, for the pair $(F, B) \coloneqq (E_{-m}, A_{-m})$ the one-sided fractional 
power scale $[F_\alpha, B_\alpha \ ; \ \alpha \in \R^+]$ exists and we obtain the \textit{extrapolated power scale of order $m \in \N \setminus \{ 0\}$ generated by $(E, A)$} as 
\begin{equation*}
	(E_\alpha, A_\alpha) \coloneqq (F_{\alpha + m, B_{\alpha +m }}) \quad \textnormal{for }  - m \leq \alpha < \infty, 
\end{equation*}
which is also a Banach-scale, cf. \cite[Thm.~V 1.3.8]{amann1995linear}. \\ 
Note that, following the notation from Amann, we set
\begin{equation*}
	E^\# \coloneqq E', \quad A^\# \coloneqq A'
\end{equation*}
to avoid confusion with too many dashes.
Suppose that the Banach space $E$ is reflexive. Then, the dual operator $A^\#$ is of class $\mathcal{BIP}(E^\#)$ and the extrapolated fractional power scale $[(E^\#_\alpha, A^\#_\alpha) \ ; \ \alpha \in [-m, \infty )]$ is therefore well-defined for all $m \in \N$. \\ 
The arguments leading to \cite[Thm.~1.4.9]{amann1995linear} further show that the pairs $(E_\alpha, A_\alpha)$, $ \alpha \in (-\infty, -m]$, are independent of the order $m \in \N$ of the extrapolated fractional power scale, such that our inductive construction already yields the the two-sided power scales 
\begin{equation*}
		[(E_\alpha, A_\alpha) \ ; \ \alpha \in  \R ], \quad [(E^\#_\alpha, A^\#_\alpha) \ ; \ \alpha \in \R ]. 
\end{equation*}
Moreover, by \cite[Thm.~1.4.12]{amann1995linear}, these are Banach scales and it holds that 
\begin{equation*}
	(E_\alpha)'  = E^\#_{-\alpha}, \quad (A_\alpha)'  = A^\#_{-\alpha}, \quad \alpha \in \R
\end{equation*}
with respect to the duality-pairing induced by $\langle \cdot, \cdot \rangle_{E, E^\#}$.  \\ 
Lastly we consider interpolation-extrapolation scales and their relation to extrapolated fractional power scales. Given a densely injected extrapolated discrete power scale 
\begin{equation*}
	[(E_j, A_j) \, \, j \in \mathbb{Z} \cap [-m , \infty)] 
\end{equation*}
of order $m$ generated by $(E, A)$ and fixing for each $\vartheta \in (0, 1)$ an admissible interpolation functor $(\cdot, \cdot)_{\vartheta}$, we put 
\begin{equation*}
	E_\alpha \coloneqq  (E_j, E_{j+1})_{\alpha-j}, \quad A_\alpha \coloneqq E_\alpha\textnormal{-realization of } A_j
\end{equation*}
for $j < \alpha < j+1$ and $j \in \mathbb{Z} \cap [-m, \infty)$. Then the \textit{interpolation-extrapolation scale of order $m$ generated by $(E, A)$ and $(\cdot, \cdot)_{\vartheta}$} is a densely injected Banach-scale. If $A \in \mathcal{BIP}(E)$, \cite[Thm. V 1.5.4]{amann1995linear} states that the interpolation-extrapolation scale generated by $(E, A)$ and the complex interpolation functor $[\cdot, \cdot]_{\vartheta}$ is equivalent to the extrapolated fractional power scale of order $m$. Moreover, in this case it also follows, cf. \cite[Thm. V 1.5.5]{amann1995linear}, that if $A \in \mathcal{BIP}(E)$ with power angle $\theta_A$, then $A_{\alpha} \in \mathcal{BIP} (E_\alpha)$ with power angle $\theta_A$.

\part{The elastic case} \label{part:nve}

The aim of this first part is to prove Theorem \ref{thm:existence_strong_nve}, which establishes the existence of a unique solution to the original Cahn--Hilliard--Biot system on a short time interval via Banach's fixed-point theorem. \\ 
After introducing all necessary assumptions in Section \ref{sec:assumptions_non_ve}, we proceed in Section \ref{sec:proof_nve} by reducing the system through the elimination of dispensable variables to two evolution equations that exhibit maximal regularity. Using Lipschitz estimates for the corresponding right-hand side that we shall prove in Section \ref{sec:contraction_non_ve}, we can apply the aforementioned contraction principle on a suitable linearization, concluding the proof. \\ 
We emphasize that this strategy crucially relies on establishing maximal regularity in Section \ref{section:inverse_non_ve}, and point out that this requires showing that the sum of two highest-order operators generates an analytic semigroup over an appropriate Hilbert space.  

\section{Assumptions: elastic case }\label{sec:assumptions_non_ve} 

The Cahn--Hilliard--Biot system \eqref{eq:strong_formulation} contains many fixed functions to include source terms, describe applied forces, and model the intrinsic energy that is crucial to the gradient flow structure. Since the existence of solutions, and therefore the validity of our results, does necessarily depend on their properties, we now state the precise assumptions required in the proof of Theorem \ref{thm:existence_strong_nve}. 

\begin{enumerate}[label = {(I.A\arabic*)} ]
	\item  \label{II.A:domain}
	Let $\Omega \subset \mathbb{R}^n$ be a bounded $C^{4}$-domain in dimension $n \leq 3$. 
	 Additionally, assume that the subset $\Gamma_D \subset \partial \Omega$ is relatively closed such that $\mathcal{H}^{n-1} (\Gamma_D) > 0$ and $\overline{\Gamma}_D \cap \overline{\Gamma}_N = \emptyset$, where $\Gamma_N\coloneqq \partial \Omega \setminus \Gamma_D$. 
	\item  \label{II.A:psi} 
	The potential $\psi : \R \rightarrow \R$ is of class $\psi \in C^3(\R)$.  
	\item \label{II.A:W}
	The elastic free energy density $W \in C^1(\mathbb{R} \times \mathbb{R}^{n \times n}_{sym})$ is of the form 
	\begin{equation*}
		W(\vphi', \E') = \C(\vphi') (\E' - \Tau(\vphi')):  (\E' - \Tau(\vphi')), 
	\end{equation*} 
	where $\C: \R \rightarrow \mathcal{L}(\R^{n \times n}_{sym})$ is a tensor-valued function of class $\bm{C}^2(\R, \R^{4n})$ such that all derivatives are bounded. We require it to fulfill the standard assumptions of linear elasticity, i.e., $\C(\vphi')$ is  symmetric and uniformly positive definite on $\R^{n \times n}_{sym}$, mapping symmetric matrices to symmetric matrices such that 
	\begin{align*}
		\E: \C(\vphi') \E &\geq c |\E|^2,\\ 
		\mathcal{D} : \C(\vphi') \E &= \C(\vphi')\mathcal{D} : \E
	\end{align*}
	for all symmetric matrices $\E, \mathcal{D} \in \R^{n\times n}_{sym}$ and all $\vphi' \in \R$.\\ 
	The eigenstrain $\Tau : \R \rightarrow \R^{n \times n}_{sym}$ is a matrix-valued function of class $\bm{C}^2(\R, \R^{n \times n}_{sym})$ with bounded derivatives.
	\item \label{II.A:g}
	The function $\bm{g} : \Gamma_N \rightarrow \mathbb{R}^n$, modeling applied outer forces, fulfills $\bm{g} \in L^s(0, T; \bm{H}^{\sfrac{3}{2}}(\Gamma_N))$.
	\item \label{II.A:f}
	The function $\bm{f} \colon \Omega \rightarrow \R^n$, modeling body forces, satisfies $\bm{f} \in \bm{L}^s(0, T; \bm{H}^1_{\Gamma_D}(\Omega))$.
	\item \label{II.A:m}
	There exists a constant $\underline{m} > 0$ such that the mobility $m \in C^{1}(\mathbb{R})$ fulfills 
	\begin{equation*}
		\underline{m} \leq m(z) \quad \textrm{for all} \quad z \in \mathbb{R}. 
	\end{equation*}
	\item \label{II.A:kappa}
	There exist a constant $\underline{\kappa} > 0$ such that the function $\kappa \in C^{1}(\mathbb{R})$ fulfills 
	\begin{equation*}
		\underline{\kappa} \leq \kappa(z) \quad \textrm{for all} \quad z \in \mathbb{R}. 
	\end{equation*}
	\item \label{II.A:phase_coefficients}
	The maps $\alpha, M : \R \rightarrow \R$ are of class $C^{2}(\R)$. Additionally, suppose that $M$ is uniformly positive, i.e., there exists some $\underline{M} >0$ with $M(z) \geq \underline{M}$ for all $z \in \R$. 
\end{enumerate}

\new{
	\begin{figure}[h!]
		\centering
		\begin{tikzpicture}[scale=2, \ifdraft{blue}{}]
			
			\draw[thick]
			(0,0)
			.. controls (1.2,-0.4) and (2.2,-0.2) .. (2.6,0.6)
			.. controls (3.0,1.4) and (2.2,2.2) .. (1.4,2.4)
			.. controls (0.6,2.6) and (-0.2,1.8) .. (-0.3,1.0)
			.. controls (-0.4,0.4) and (-0.2,0.1) .. (0,0);
			
			\draw[ultra thick]
			(1.2,1.0)
			.. controls (1.6,0.7) and (2.1,0.9) .. (2.0,1.3)
			.. controls (1.9,1.8) and (1.4,1.9) .. (1.1,1.6)
			.. controls (0.9,1.3) and (1.0,1.1) .. (1.2,1.0);
			
			\node at (0.5,0.8) {$\Omega$};
			
			\node[right] at (2.05,1.25) {$\Gamma_D$};
			\node[left] at (-0.45,1.4) {$\Gamma_N$};
			
		\end{tikzpicture}
		\caption{Sketch of an admissible domain in 2D}
	\end{figure}
}

\begin{remark} 
	Note that in the following we will impose an additional assumption \ref{II.A:source_terms} on the source terms $S_s$ and $S_f$. \\ 
	In view of the explicit formula for the elastic energy density \ref{II.A:W}, tedious calculations lead to the following growth conditions: 
	\begin{enumerate}[label = (\roman*)]
		\myitem{(I.A3.2)} \label{II.A:W_growth_conditions}
		There exists a constant $C_2 > 0$ such that for all $\vphi' \in \mathbb{R}$ and all symmetric $\E' \in  \mathbb{R}^{n \times n}$,
		\begin{align*}
			\abs{W_{,\vphi}(\vphi', \E') }  + \abs{  W_{,\vphi \vphi}(\vphi', \E') } + \abs{ W_{,\vphi \vphi \vphi}(\vphi', \E') }  & \leq C_2 \left( \abs{\E'}^2 + \abs{\vphi'}^2 + 1  \right), \\ 
			\abs{  W_{,\vphi \E}(\vphi', \E') }  + \abs{  W_{,\vphi \vphi \E}(\vphi', \E') } & \leq C_2 \left( \abs{\E'} + \abs{\vphi'}  + 1\right), \\ 
			\abs{  W_{,\vphi \E \E}(\vphi', \E') }  & \leq C_2. \\ 
		\end{align*}
	\end{enumerate}
\end{remark}

\addtocontents{toc}{\SkipTocEntry}
\subsection*{Higher regularity for elliptic systems}\label{subsectio:higher_regularity_elliptic}

We consider elliptic problems of the form
	\begin{align}\label{prob:elliptic_mixed_bc}
	\left\{ 
	\begin{alignedat}{2}
		- \nabla \cdot \big( \C(\phi) \E(\bm{v})   \big)  &= \bm{f} & & \quad \textnormal{in } \Omega,\\ 
		\bm{v}  &= \bm{0} & &\quad \textnormal{on }   \Gamma_D, \\ 
		\C(\phi) \E(\bm{v})  \bm{n} &= \bm{g} & &\quad \textnormal{on }  \Gamma_N. 
	\end{alignedat} 
	\right. 
\end{align}
Assuming sufficient regularity on the domain along with $\mathcal{H}^{n-1}(\Gamma_D)  \neq \emptyset$, we deduce with the help of the Lax-Milgram theorem that the corresponding operator describes a topological isomorphism between the spaces $\bm{W}^{1, 2}_{\Gamma_D}(\Omega)$ and $\bm{W}^{-1, 2}_{\Gamma_D}(\Omega)$. 
In particular, a function $\bm{v} \in \bm{W}^{1, 2}_{\Gamma_D}\new{(\Omega)}$ is mapped to 
 \begin{equation*}
 	 \Big(  \bm{w} \mapsto \int_\Omega \C(\phi) \E(\bm{v}) \colon \E(\bm{w}) \dx \Big) \in \bm{W}^{-1, 2}_{\Gamma_D}(\Omega),
 \end{equation*}
 and for $\bm{f} \in \bm{W}^{-1, 2}_{\Gamma_\new{D}}\new{(\Omega)}$, $\bm{g} \in W^{-\frac{1}{2},2}(\Gamma_N)$, we can define the tuple 
 \begin{equation*}
 	(\bm{f}, \bm{g}) \coloneqq  
 	{}_{\bm{W^{-1, 2}}}\langle \bm{f}, \cdot \rangle_{\bm{W}^{1, 2}} 
 	+ {}_{\bm{W}^{ -\frac{1}{2}, 2} (\Gamma_N )} \langle \bm{g}, \cdot \rangle_{\bm{W}^{\frac{1}{2}, 2}(\Gamma_N)}  \in \bm{W}^{-1, 2}_{\Gamma_D}(\Omega). 
 \end{equation*}
Classical texts, see \cite{MR125307}, \cite[Thm. 9.31]{renardyrogers} and the references cited therein, derive higher regularity in the case of pure (Dirichlet or Neumann) boundary conditions. Since we only consider the case where the segments on which different boundary conditions are imposed do not meet, we could divide our domain, derive higher regularity on the different parts, and eventually obtain the desired result upon combining our estimates. We particularly refer to \cite[Thm. 4.18]{mclean2000}, which includes all necessary assumptions for this approach.\\  However, the proof of the following theorem, which can be found in Appendix \ref{sec:w^2^p}, relies on results from maximal regularity theory. 
\par 
\medskip 
\noindent

\begin{theorem}\label{thm:higher_elliptic}
	Suppose that $\Omega \subset \mathbb{R}^n$ is a bounded $C^{2}$-domain. 
	Additionally, assume that the subset $\Gamma_D \subset \partial \Omega$ is relatively closed such that $\mathcal{H}^{n-1} (\Gamma_D) > 0$ and $\overline{\Gamma}_D \cap \overline{\Gamma}_N = \emptyset$, where $\Gamma_N\coloneqq \partial \Omega \setminus \Gamma_D$.\\ 
	Let $\phi \in W^{1, q}(\Omega)$ for some $q > n$, and $\C: \R \rightarrow \mathcal{L}(\R^{n \times n}_{sym})$ be a tensor-valued function of class $C^{2}(\R, \R^{4n})$ that is uniformly positive on symmetric matrices. Then there exists a constant $C > 0$, which only depends on the norm of $\C(\phi) \in  W^{1, q}(\overline{\Omega})$ and its ellipticity parameter, such that the unique solution $\bm{v}$ to the elliptic problem \eqref{prob:elliptic_mixed_bc}, where $\bm{f} \in \bm{L}^2(\Omega)$ and $\bm{g} \in \bm{H}^{\sfrac{1}{2}}(\Gamma_N)$, is in $\bm{H}^2(\Omega)$, and can be estimated by 
	\begin{equation*}
		\norm{\bm{v}}_{\bm{H}^2} \leq C \big( \norm{\bm{f}}_{ \bm{L}^2} +  \norm{\bm{g}}_{ \bm{H}^{\sfrac{1}{2}}(\Gamma_N) } \big). 
	\end{equation*}
	In particular, the operator 
	\begin{equation*}
		\mathcal{C}(\phi) : \bm{H}^2(\Omega) \rightarrow (\bm{L}^2(\Omega), \bm{H}^{\sfrac{1}{2}} (\Gamma_N)), 
		\quad \bm{v} \mapsto (- \nabla \cdot \big( \C(\phi) \E(\bm{v}), \C(\phi) \E(\bm{v}) \bm{n})
	\end{equation*}
	is invertible for all $\phi \in  W^{1, q}(\Omega)$ and the inverse $\mathcal{C}^{-1}(\phi)$ is uniformly bounded for all $\norm{\phi }_{ W^{1, q}} \leq R$. 
\end{theorem}

\begin{proof}
	This is an immediate consequence of Theorem \ref{thm:elastic_non_symmetric} and \eqref{eq:symmetry_C}. 
\end{proof}

\medspace

\section{Proof of Theorem \ref{thm:existence_strong_nve}} \label{sec:proof_nve}

The central ideas of the proof of Theorem \ref{thm:existence_strong_nve} can be summarized as follows: in order to reduce the question of existence of solutions to a fixed-point problem, we eliminate the variables $\mu$ and $p$, inserting the respective expression into the equations for $\vphi$ and $\theta$. By means of elliptic theory, in particular Theorem~\ref{thm:higher_elliptic}, the displacement $\bm{u}$ is uniquely determined by $\vphi$ and $\theta$, allowing us to further reduce the system into two evolution equations for the latter two variables. \par
At this point, it is crucial to establish maximal regularity and therefore the invertibility of the corresponding operator, see Proposition \ref{prop:inverse_boundedness_non_ve}, such that the evolution equations mentioned above can be rewritten as an equivalent fixed-point problem. \par  
After linearization, we will eventually employ a contraction principle and deduce unique solvability of the Cahn--Hilliard--Biot system on a short time scale, as asserted in the theorem. 

\addtocontents{toc}{\SkipTocEntry}
\subsection*{Boundary conditions} Before considering the equations in the bulk, we need to derive the boundary conditions for the higher-order derivatives of $\vphi$ and $\theta$ which are compatible with \eqref{bc:contac_angle}-\eqref{bc:discplacement_neumann}. From \eqref{eq:strong_formulation_p} we deduce by applying \eqref{bc:contac_angle} and \eqref{bc:mass_conservation}
\begin{equation*}
	0 = \nabla p \cdot \bm{n} = M(\vphi) \big( \nabla \theta \cdot \bm{n}  - \alpha(\vphi) \nabla (\nabla\cdot \bm{u}) \cdot \bm{n} \big) \quad \textnormal{on } (0, T) \times \partial \Omega , 
\end{equation*}
and recalling that $M$ is uniformly positive, we obtain the equivalent equation 
\begin{equation} \label{eq:bc_theta_higher}
	\nabla \theta \cdot \bm{n}  = \alpha(\vphi) \nabla (\nabla \cdot \bm{u}) \cdot \bm{n} \quad \textnormal{on } (0, T) \times \partial \Omega. 
\end{equation}
Similarly, it follows that  
\begin{equation} \label{eq:bc_phi_higher}
	 \nabla\Delta \vphi \cdot \bm{n} = W_{,\vphi \E} (\vphi, \E(\bm{u})) \nabla \E(\bm{u}) \bm{n} + \alpha'(\vphi) p \nabla(\nabla \cdot  \bm{u}) \cdot \bm{n} \quad \textnormal{on } (0, T) \times \partial \Omega. 
\end{equation}
Instead of imposing these boundary conditions in the sense of traces, we will use a weaker formulation that incorporates them implicitly, avoiding the highly technical computations required when dealing with inhomogeneous boundary conditions in the framework of maximal regularity. 

\addtocontents{toc}{\SkipTocEntry}
\subsection*{Reformulation of the fluid equation}
Our first goal is to reformulate \eqref{eq:strong_formulation_theta}, replacing all instances of $p$ and $\bm{u}$ by terms that solely depend on $\theta$ and $\vphi$. To that end, recall that \eqref{eq:strong_formulation_p} provides the identities 
\begin{equation}\label{eq:p_theta_identities}
	 p = M(\vphi)(\theta - \alpha(\vphi) \nabla \cdot \bm{u})
	 \quad \textnormal{and} \quad 
	 \theta = \frac{1}{M(\vphi)} p + \alpha(\vphi) \nabla \cdot \bm{u}. 
\end{equation}
In the absence of the visco-elastic regularization $\C_{\nu} (\E(\pt \bm{u}))$, the equations \eqref{eq:strong_formulation_u}-\eqref{eq:strong_formulation_u_2} simplify to the elliptic problem
\begin{equation*}
	- \nabla \cdot  \big( \C (\vphi) \E(\bm{u}) -   \C (\vphi) \Tau(\vphi) - \alpha(\vphi) p \bm{I}  \big) = \bm{f}, 
\end{equation*}
or equivalently 
\begin{equation*}
	- \nabla \cdot ( \C(\vphi) \E(\bm{u} ) + \alpha^2(\vphi) M(\vphi) (\nabla \cdot \bm{u}) \bm{I})=  - \nabla \cdot (   \C (\vphi) \Tau(\vphi) + \alpha(\vphi) M(\vphi)\theta  \bm{I} ) + \bm{f}
\end{equation*}
which are, due to the boundary conditions imposed in \ref{II.A:domain}, uniquely solvable by Theorem \ref{thm:higher_elliptic} with solution
\begin{align*}
	\bm{u} &= \mathcal{C}^{-1}(\vphi) \big( - \nabla \cdot (   \C (\vphi) \Tau(\vphi) + \alpha(\vphi) p \bm{I} ) + \bm{f}, \bm{g} \big)\\ 
	 &= \tilde{\mathcal{C}}^{-1}(\vphi) \big(- \nabla \cdot (   \C (\vphi) \Tau(\vphi) + \alpha(\vphi) M(\vphi) \theta  \bm{I} ) + \bm{f}, \bm{g} \big).
	\numberthis
	\label{eq:u_solution_elliptic}
\end{align*} 
Here, $\tilde{\mathcal{C}}^{-1}(\vphi)$ is the inverse of the operator 
\new{
\begin{equation*}
	\bm{u} \mapsto - \Big( \nabla \cdot \big( \C(\vphi) \E(\bm{u} ) + \alpha^2(\vphi) M(\vphi) (\nabla \cdot \bm{u}) \bm{I} \big), \big( \C(\vphi) \E(\bm{u} ) + \alpha^2(\vphi) M(\vphi) (\nabla \cdot \bm{u}) \bm{I} \big) \bm{n} 
				\Big), 
\end{equation*}
}which also possesses all properties required in Theorem \ref{thm:higher_elliptic}. \\ 
Neglecting terms of lower order, we define  
\begin{equation*}
	\tilde{\mathcal{B}}(\phi) :  L^2(\Omega) \rightarrow L^2(\Omega), \quad 
	q  \mapsto \frac{1}{M(\phi)} q  + \alpha(\phi)\,  \nabla \cdot \big[ \mathcal{C}^{-1}(\phi) \big( - \nabla  (   \alpha(\phi) q ), \bm{0} \big )\big],
\end{equation*}
and, comparing the right-hand side to \eqref{eq:p_theta_identities}, observe that $\tilde{\mathcal{B}}(\phi)$ provides a formula for $\theta$ that merely depends on $\phi$ and $q$, eliminating instances of $\bm{u}$ up to terms of lower order. In view of our goal, we would like to obtain a similar expression for $p$, and therefore consider the inverse.\\ 
The following lemma provides a rigorous argument for invertibility of this operator. 

\begin{lemma}\label{lem:properties_B_L^2}
 For any $\phi \in W^{1, q}$, $q > n$, the operator $\tilde{\mathcal{B}} (\phi)$ is invertible, positive and self adjoint. Moreover, there exists a constant $C > 0$ such that $\bnorm{\tilde{\mathcal{B}}^{-1}(\phi)}_{\mathcal{L}(L^2)} \leq C$. 
\end{lemma}

\begin{proof}
	First of all, we observe with the help of Theorem \ref{thm:compositon}, \ref{II.A:W} and standard embedding results, that $\C(\phi)$ is of class $W^{1, q}(\Omega)$. Now, let $\lambda \in [0, 1]$ and set  
	\begin{equation*}
		\tilde{\mathcal{B}}_\lambda (\phi) : L^2(\Omega) \rightarrow L^2(\Omega), \quad 
		q \mapsto  \frac{1}{M(\phi)} q  +  \lambda \alpha(\phi)\,  \nabla \cdot \big[ \mathcal{C}^{-1}(\phi) \big( - \nabla \cdot (\alpha(\phi) q \bm{I} ), \bm{0}\big) \big]. 
	\end{equation*}
	Moreover, we introduce $\bm{v} \coloneqq \mathcal{C}^{-1}(\phi) \big( - \nabla \cdot (\alpha(\phi) q \bm{I} ), \bm{0} \big ) \in \bm{W}^{1, 2}_{\Gamma_D} $, and recalling the discussion in Section~\ref{subsectio:higher_regularity_elliptic}, find that $\bm{v}$ is the unique weak solution to the corresponding elliptic problem, satisfying the weak formulation 
	\begin{equation} \label{eq:C_weak_formulation}
		\int_\Omega \C(\phi) \E(\bm{v}) \colon \E(\bm{w}) \dx 
		= \int_\Omega \alpha (\phi) q\,  \nabla \cdot \bm{w} \dx \quad \textnormal{for all } \bm{w} \in \bm{W}^{1, 2}_{\Gamma_D} (\Omega). 
	\end{equation}
	 Thus, it is possible to use $\bm{v}$ as a test function, leading to 
	\begin{equation}\label{eq:invertible_est_1}
		\int_\Omega \alpha (\phi) q\,  \nabla \cdot \bm{v} \dx \geq c \norm{\bm{v}}_{\bm{H}^1}^2  \geq 0. 
	\end{equation}
	Exploiting this inequality, we will now show that the operators $\TA_\lambda$ are injective and have closed range for all $\lambda \in [0, 1]$, which implies that they are in particular semi-Fredholm. Pointing out that $M$ is uniformly positive and bounded, cf. \ref{II.A:phase_coefficients}, we observe that the inner product $(u, v ) \mapsto ( M(\phi)\,  u, v) _{L^2}$ gives rise to an equivalent norm on $L^2(\Omega)$ for any $\phi$ as in the assumptions. Thus, 
	\begin{align*}
		\norm{\tilde{\mathcal{B}}_\lambda (\phi) q}_{L^2}^2 
		&= \norm{\frac{1}{M(\phi)} q  +  \lambda \alpha(\phi)\,  \nabla \cdot \big[ \mathcal{C}^{-1}(\phi) \big( - \nabla \cdot ( \alpha(\phi) q \bm{I} ), \bm{0} \big )\big]}_{L^2}^2 \\ 
		& \geq c \bigg(  q  +  \lambda M(\phi)\alpha(\phi)\,  \nabla \cdot \bm{v} , \frac{1}{M(\phi)} q  +  \lambda \alpha(\phi)\,  \nabla \cdot \bm{v}  \bigg)_{L^2} \\ 
		& \geq c \big( \frac{1}{M(\phi)} q , q \big)_{L^2}
		+ 2c  \lambda \int_\Omega \alpha(\phi) q \, \nabla \cdot \bm{v} \dx 
		+ c \lambda\new{^2} \big( M(\phi) \alpha (\phi) \nabla\cdot \bm{v} , \alpha (\phi) \nabla\cdot \bm{v}\big)_{L^2} \\ 
		& \geq c \norm{q}_{L^2}^2 + c \lambda\big( \norm{\bm{v}}_{\bm{H}^1}^\new{2} +  \new{\lambda} \norm{\alpha(\phi) \nabla \cdot \bm{v} }_{L^2}^\new{2} \big), \numberthis \label{eq:tildeA_bounded_below}
 	\end{align*}
 	for all $\lambda \in [0, 1]$ and any $q \in L^2(\Omega)$, which already implies the assertion. \new{Note that for the last inequality we used \eqref{eq:invertible_est_1} along with the fact that $M$ is uniformly positive. Moreover, due to continuity and since $\phi \in W^{1, q}(\Omega) \hookrightarrow C^0(\overline{\Omega})$, the function $M(\phi)$ is bounded, which implies an estimate from below for the reciprocal.}\\ 
 	Obviously, the operator $\tfrac{1}{M(\phi)} \mathcal{I}$ is bijective on $L^2(\Omega)$ and thus has Fredholm index $0$.  Since the mapping $[0, 1] \ni \lambda \mapsto \tilde{\mathcal{B}}_\lambda \in \mathcal{L}(L^2(\Omega)) $ describes a homotopy between $\tfrac{1}{M(\phi)} \mathcal{I}$ and $\tilde{\mathcal{B}}(\phi)$, Fredholm theory further implies that the indices must already coincide, i.e., 
 	\begin{equation*}
 		\textnormal{index} (\tilde{\mathcal{B}}(\phi)) = \textnormal{index} (\tfrac{1}{M(\phi)} \mathcal{I}) = 0. 
 	\end{equation*}
 	 Emphasizing the fact that $\ker(\tilde{\mathcal{B}}(\phi))$ is trivial, we deduce that $\textnormal{coker}(\tilde{\mathcal{B}}(\phi)) = {0}$, and conclude that ${\mathcal{B}}(\phi)$ is surjective. \\ 
 	 The fact that this operator is also positive can be seen from the following computation, which utilizes \ref{II.A:phase_coefficients} \new{and once again \eqref{eq:invertible_est_1}}
 	 \begin{equation*}
 	 	\big( \tilde{\mathcal{B}}(\phi) q, q \big)_{L^2}
 	 	= \big( \frac{1}{M(\phi)}q  +   \alpha(\phi)\,  \nabla \cdot \bm{v}, q \big)_{L^2} 
 	 	=  \big( \frac{1}{M(\phi)}q , q \big) +  \big(  \alpha(\phi)\,  \nabla \cdot \bm{v}, q \big)_{L^2} 
 	 	\geq c \norm{q}_{L^2}^2 
 	 	\geq 0 . 
 	 \end{equation*}
 	 To show that $\tilde{\mathcal{B}}(\phi)$ is self adjoint,  let $ q, r \in L^2(\Omega)$, and set $\bm{v} = \mathcal{C}^{-1}(\phi) \big( - \nabla \cdot ( \alpha(\phi) q \bm{I} ), \bm{0} \big )$ along with $\bm{w} = \mathcal{C}^{-1}(\phi) \big( - \nabla \cdot ( \alpha(\phi) r \bm{I} ), \bm{0} \big) $. Using $\bm{v}$ and $\bm{w}$ as test functions in \eqref{eq:C_weak_formulation} for the respective problems, a simple computation shows that $\tilde{\mathcal{B}}(\phi)$ is symmetric. Since $\tilde{\mathcal{B}}(\phi)$ is defined on the whole Hilbert space $L^2(\Omega)$ and bijective, symmetry already implies that  $\tilde{\mathcal{B}}(\vphi_0)$ is self adjoint.\\ 
 	 For the last assertion, note that the constant $c$ in \eqref{eq:tildeA_bounded_below} does not depend on $\phi$ but solely on the bounds of $\C, M$ and $\alpha$. Hence, it holds for all $\phi$ satisfying the assumptions that
 	 \begin{equation*}
 	 	\norm{\tilde{\mathcal{B}}^{-1}(\phi) \vartheta}_{L^2} \leq \frac{1}{c} \norm{\vartheta}_{L^2} \quad \textnormal{for all } \vartheta \in L^2(\Omega), 
 	 \end{equation*}
 	 which implies the assertion. 
\end{proof}
We shall now prove that the inverse of $\tilde{\mathcal{B}}(\phi)$ is given by 
\begin{equation*}
	\tilde{\mathcal{A}}(\phi) :  L^2(\Omega) \rightarrow L^2(\Omega), \quad 
	\vartheta \mapsto M(\phi)\big( \vartheta - \alpha(\phi) \nabla \cdot \big[ \tilde{\mathcal{C}}^{-1}(\phi) \big( - \nabla  (   \alpha(\phi) M(\phi) \vartheta )  \big), \bm{0} \big]  \big), 
\end{equation*}
and that this operator exhibits the same properties as $\tilde{\mathcal{B}}(\phi)$.

\begin{corollary}\label{cor:properties_TA}
	For any $\phi \in W^{1, q}(\Omega)$, $q > n$, the operator $\TA(\phi)$ is invertible, positive and self adjoint. Moreover, for every $R > 0$ there exists a constant $C = C(R) > 0$ such that 
	\begin{align*}
		\norm{ \TA(\phi) \vartheta}_{H^1} \leq C\norm{\vartheta}_{H^1}
		\quad  &\textnormal{ for all } \theta \in H^1(\Omega), \\ 
		\norm{ \TA^{-1}(\phi) q}_{H^1} \leq C\norm{q}_{H^1}
		\quad  &\textnormal{ for all } q \in H^1(\Omega), 
	\end{align*}
	for all $\phi \in W^{1, q}(\Omega)$ satisfying $\norm{\phi}_{W^{1, q}} \leq R$. In particular, the restriction 
	\begin{equation*}
		 \TA(\phi)_{|H^1} : H^1(\Omega) \rightarrow H^1(\Omega) 
	\end{equation*}
	is bijective. 
\end{corollary}
\begin{proof}
	In the lemma above we have shown that for every $\vartheta \in L^2(\Omega)$ there exists a unique $q \in L^2(\Omega)$ such that $\vartheta = \tilde{\mathcal{B}}(\phi)q$. Hence, we find that 
	\begin{align*}
		\TA(\phi) \vartheta 
		&= \TA(\phi)\bigg( \frac{1}{M(\phi)} q  + \alpha(\phi)\,  \nabla \cdot \big[ \mathcal{C}^{-1}(\phi) \big( - \nabla  (   \alpha(\phi) q ), \bm{0} \big ) \big] \bigg)\\ 
		& = \begin{aligned}[t]
			q &+ M(\phi)  \alpha(\phi)\,  \nabla \cdot \big[ \mathcal{C}^{-1}(\phi) \big( - \nabla  (   \alpha(\phi) q ), \bm{0}\big)\big]\\ 
			&- M(\phi) \alpha(\phi) \nabla \cdot \big[ \tilde{\mathcal{C}}^{-1}(\phi) 
			\big( -\nabla ( \alpha(\phi)  q  + \alpha^2(\phi) M(\phi)\,  \nabla \cdot [ \mathcal{C}^{-1}(\phi) \big( - \nabla  (   \alpha(\phi) q ), \bm{0}\big )]  ) , \bm{0} \big)    \big]. 
		\end{aligned}
	\end{align*}
	As above, we abbreviate $\bm{v} \coloneqq \mathcal{C}^{-1}(\phi) \big( - \nabla (\alpha(\phi) q ), \new{\bm{0}}\big )$ and obtain by adding $- \nabla\cdot ( \alpha^2(\phi) M(\phi) (\nabla \cdot  \bm{v})  \bm{I})$
	\begin{equation*}
		- \nabla \cdot \big( \C(\phi) \E(\bm{v})+ \alpha^2(\phi) M(\phi) (\nabla \cdot  \bm{v} ) \bm{I} \big) = - \nabla \big( \alpha(\phi) q  + \alpha^2(\phi) M(\phi) (\nabla \cdot  \bm{v} ) \big). 
	\end{equation*}
	Exploiting the uniqueness of solutions to elliptic problems, we infer 
	\begin{equation*}
		\bm{v} = \tilde{\mathcal{C}}^{-1}(\phi) \big( - \nabla (\alpha(\phi) q   + \alpha^2(\vphi) M(\phi) \nabla \cdot  \bm{v}  ), \bm{0}\big )
	\end{equation*}
	and hence, 
	\begin{align*}
		 \TA(\phi)\tilde{\mathcal{B}}(\phi) q
		&= q + M(\phi)  \alpha(\phi)\,  \nabla \cdot \bm{v} 
		- M(\phi) \alpha(\phi) \nabla \cdot \big[ \tilde{\mathcal{C}}^{-1}(\phi) 
		\big( -\nabla ( \alpha(\phi)  q  + \alpha^2(\phi) M(\phi)\,  \nabla \cdot \bm{v}   ) ,\bm{0} \big)   \big] \\ 
		&= q + M(\phi)  \alpha(\phi)\,  \nabla \cdot \bm{v}  - M(\phi)  \alpha(\phi)\,  \nabla \cdot \bm{v}  = q. 
	\end{align*}
	Conversely, we can use this result and the fact that for every $\vartheta \in L^2(\Omega)$ there exists an unique $q \in L^2(\Omega)$ such that $\vartheta = \tilde{\mathcal{B}}(\vphi)q$ to similarly deduce  
	\begin{equation*}
		\tilde{\mathcal{B}}(\phi) \TA(\phi) \vartheta =\tilde{\mathcal{B}}(\phi) \big( \TA(\phi) \tilde{\mathcal{B}}(\phi)  q \big) =   \tilde{\mathcal{B}}(\phi) q = \vartheta. 
	\end{equation*}
	Hence, it follows that $\TA(\phi)$ is indeed the inverse to the bijective operator $\tilde{\mathcal{B}}(\phi)$. 
	Moreover, since we have shown above that $\tilde{\mathcal{B}} (\phi) = \TA^{-1}(\phi)$ is self adjoint, so is $\TA(\phi)$
	\begin{equation*}
		\big( \TA(\phi) \vartheta, \vartheta \big)_{L^2}
		= \big( \TA(\phi) \vartheta,  \TA^{-1}(\phi) \TA(\phi) \vartheta \big)_{L^2}
		= \big( \TA^{-1}(\phi)\TA(\phi) \vartheta,   \TA(\phi) \vartheta \big)_{L^2} 
		\big( \vartheta,  \TA(\phi) \vartheta \big)_{L^2} . 
	\end{equation*}
	Now, let $\vartheta \in L^2(\Omega)$ and set $q = \TA(\phi) \vartheta$. Then,  
 	\begin{equation*}
		\big( \TA(\phi) \vartheta, \vartheta \big)_{L^2}
		= \big( q , \TA^{-1}(\phi) q \big)_{L^2}
			= \big( q ,{\mathcal{B}} (\phi) q \big)_{L^2}
		\geq 0, 
	\end{equation*}
	which means that $\TA(\phi)$ is indeed positive. \\ 
	Using the Rellich–Kondrachov embedding theorem and Hölder's inequality, we observe that for two arbitrary functions $f \in W^{1, q}(\Omega)$, $g\in H^1(\Omega)$ 
		\begin{align}\label{eq:H^2_estimate}
			\norm{f\, g}_{H^1}\leq C \norm{f}_{W^{1,q}} \norm{g}_{H^1}.
		\end{align}
	Due to Theorem \ref{thm:compositon} and Assumption \ref{II.A:phase_coefficients} we infer $M(\phi), M(\phi)\alpha(\phi)  \in W^{1, q}(\Omega)$. Regardless of the dimension, a similar argument implies that the coefficients of the operator $\mathcal{C}(\phi)$ are of class $W^{1, q}(\Omega)$, such that we can we can invoke Theorem~\ref{thm:higher_elliptic} to obtain for all $\vartheta \in H^1(\Omega)$
	\begin{align*}
		\norm{\TA(\phi) \vartheta }_{H^1}
		&\leq \norm{M(\phi) \vartheta}_{H^1} + \norm{M(\phi) \alpha(\phi) \nabla \cdot \big[ \tilde{\mathcal{C}}^{-1}(\phi) \big( - \nabla  (   \alpha(\phi) M(\phi) \vartheta ), \bm{0}  \big) \big]}_{H^1}\\ 
		& \leq C \norm{\vartheta}_{H^1} + C \norm{\tilde{\mathcal{C}}^{-1}(\phi) \big[ - \nabla  (   \alpha(\phi) M(\phi) \vartheta, \bm{0} )  \big]}_{H^2}\\ 
		&\leq C \norm{\vartheta}_{H^1} + C \norm{{\mathcal{C}}^{-1}(\phi)}_{\mathcal{L}((\bm{L}^2, \bm{H}^{\sfrac{1}{2}}(\Gamma_N))   , \bm{H}^2) } \norm{ - \nabla  (   \alpha(\phi) M(\phi) \vartheta )  }_{L^2}\\
		&\leq C \norm{\vartheta}_{H^1} + C \norm{\alpha(\phi) M(\phi) }_{W^{1,q}} \norm{\vartheta}_{H^1}
		\leq C \norm{\vartheta}_{H^1}, 
		\numberthis \label{eq:TA_theta_H^1}
	\end{align*}
	and note that $\TA(\phi) (H^1(\Omega))  \subset H^1(\Omega)$. 
	Analogously, one can show 
	\begin{equation}\label{eq:TA^-1_theta_H^1_1}
		\norm{\TA^{-1}(\phi) q }_{H^1} \leq C \norm{q}_{H^1}  \quad \textnormal{for all } q \in H^1(\Omega), 
	\end{equation}
	which implies that for all $q \in H^1(\Omega)$ the unique preimage in the (larger) space $L^2(\Omega)$ is even of class $H^1(\Omega)$. Combined with the fact that $\TA(\phi)$ is bijective on $L^2(\Omega)$, this suffices to deduce invertibility for the restriction to $H^1(\Omega)$. 
\end{proof}

Now, we fix some $\vphi_0 \in W^{1,q}$, $q > n$. Since $\TA$ is self adjoint and positive, the mapping 
\begin{equation*}
	(u, v) \mapsto ( u, v )_{H} \coloneqq (u, \TA(\vphi_0)v)_{L^2},  \quad  u, v  \in L^2(\Omega)
\end{equation*}
defines an inner product on $L^2(\Omega)$, giving rise to an equivalent norm on this space.
Setting
\begin{equation*}
	(H, ( \cdot , \cdot )_{H}  ) \coloneqq (L^2(\Omega), (\cdot , \TA(\vphi_0) \cdot )_{L^2}), 
\end{equation*}
the considerations above show this to be a well-defined Hilbert space.

\begin{lemma}\label{lem:norm_equivalence_H}
The inner product $( \cdot, \cdot )_{H}$ induces a norm on the space $H$ that is equivalent to $\norm{\cdot}_{L^2}$. In particular $(H, ( \cdot , \cdot )_{H})$ is a Hilbert space and $H \cong L^2(\Omega)$.  
\end{lemma}
\begin{proof}
	As mentioned above, $( \cdot , \cdot )_{H}$ is indeed an inner product on $H$. Hence, it suffices to show that there exist $C > c> 0$ such that 
	\begin{equation}\label{eq:norm_equivalence_H}
		c \norm{u}_{L^2} \leq \norm{u}_{H} \leq C \norm{u}_{L^2} \quad \textnormal{for all } u \in H. 
	\end{equation}
	To this end, let  $v = \TA(\vphi_0) u $ and compute
	\begin{align*}
		\norm{u}_H^2 
		= (u , u )_{H} 
		&= \big(u,  \TA(\vphi_0) u \big)_{L^2}
		=  \big( \TA^{-1}(\vphi_0) v ,  v \big)_{L^2}
		\leq \| \TA^{-1}(\vphi_0) v \|_{L^2} \norm{v}_{L^2}\\
		&\leq \|\TA^{-1}(\vphi_0) \|_{\mathcal{L}(L^2)} \norm{v}_{L^2}^2
		\leq C \| \TA(\vphi_0) u \|_{L^2}^2
		\leq  C \|\TA(\vphi_0) \|_{\mathcal{L}(L^2)}^2 \norm{u}_{L^2}^2. 
		\numberthis \label{eq:norm_equivalence_H_1}
	\end{align*}
	Moreover, the Cauchy-Schwarz inequality implies 
	\begin{align*}
		\norm{v}_{L^2}^2 = \big(v, v \big)_{L^2} =  \big( v, \TA(\vphi_0) u \big)_{L^2}
		= ( v, u )_{H} 
		\leq \norm{v}_{H} \norm{u}_{H} , 
	\end{align*}
	and	after applying \eqref{eq:norm_equivalence_H_1} to find $ \norm{v}_{H} \leq C \norm{v}_{L^2}$, we obtain
	\begin{equation}\label{eq:norm_equivalence_H_2}
		 \norm{v}_{L^2} \leq C  \norm{u}_{H}. 
	\end{equation}
	Combing \eqref{eq:norm_equivalence_H_1} with \eqref{eq:norm_equivalence_H_2} therefore entails 
	\begin{align*}
		  \norm{u}_{L^2} = \|\TA^{-1}(\vphi_0) v\|_{L^2} \leq C \norm{v}_{L^2} \leq C  \norm{u}_{H} \leq \tilde{C} \norm{u}_{L^2},  
 	\end{align*}
	 and \eqref{eq:norm_equivalence_H} follows by normalization. This shows that $\norm{\cdot}_{L^2}$ and $\norm{\cdot}_{H}$ are equivalent, which also implies the second assertion. 
\end{proof}
\par 
\medskip  

We note that due to $H \cong L^2(\Omega)$, we obtain the Gelfand triplet 
\begin{equation*}
	V \coloneqq H^1(\Omega)  \hookrightarrow H \cong H' \hookrightarrow    V' \eqqcolon (H^1(\Omega))', 
\end{equation*}
where we shall denote by $V$ the space $H^1(\Omega)$ to emphasize the deviation from the canonical choices. 
The natural inclusion $V \hookrightarrow V'$ in this framework utilizes the Riesz isomorphism $\mathcal{R}_H$ of the Hilbert space $H$, i.e., for every $f \in V$ (or even $f \in H$), we have 
\begin{equation} \label{eq:dual_space_H}
	f  \mapsto \big( V \ni g \mapsto {}_{V'}\langle \mathcal{R}_H f, g \rangle_{V}  =  (f, g)_H \big) \in V'. 
\end{equation}
We stress that while $V'$ and $(H^1(\Omega))'$ coincide as sets, it would be wrong to simply identify two elements by the identity map. To see this, we first consider the Hilbert space $H$ and note that
\begin{equation*}
	\mathcal{R}_H(f) = (f, \cdot)_{H} = (f, \TA(\vphi_0) \, \cdot)_{L^2} = (\TA(\vphi_0)f,  \cdot)_{L^2} =  \mathcal{R}_{L^2}( \TA(\vphi_0) f),  
\end{equation*}
 where $\mathcal{R}_{L^2}$ and $\mathcal{R}_{H}$ are the Riesz isomorphisms associated with the Hilbert spaces $L^2(\Omega)$ and $H$, respectively. Thus, the proper identification between $(L^2(\Omega))'$ and $H'$ is given by the commutative diagram 

 \begin{center}
 	\begin{tikzcd}
 		&L^2(\Omega) \arrow[r, "\cong" above, "\textnormal{id}" below] \arrow[d, "\mathcal{R}_{L^2}(f) = (f{,}\cdot)_{L^2} " left, hook] 
 		&H \arrow[d,  "\mathcal{R}_H(f) = (f{,}\cdot)_H " right, hook]
 		\\
 		&(L^2(\Omega))'  \cong L^2(\Omega )\arrow[r, "\tilde{\mathcal{A} }(\vphi_0)" below, "\cong" above] 
 		& H \cong H' 
 	\end{tikzcd}, \\
 \end{center}
 or more explicitly 
 \begin{equation*}
 	\mathcal{I}_{H'} : (L^2(\Omega))'   \rightarrow H', \quad \zeta \mapsto \big( H \ni  g \mapsto  {}_{(L^2)'}\langle \zeta,   \TA(\vphi_0) g \rangle_{L^2}  = (\TA(\vphi_0) \zeta, g)_{L^2} \big). 
 \end{equation*}
Now consider some $\xi \in (H^1(\Omega))'$ and note that the embedding $(L^2(\Omega)) \xhookrightarrow{d} (H^1(\Omega))'$ is dense, such that there exists at least one sequence $(\zeta_n)_{n \inN} \subset L^2(\Omega)$ that converges to $\xi$ in $(H^1(\Omega))'$. To obtain an extension of $\mathcal{I}_{H'}$ it remains to show that for any such sequence the limit $\lim_{n \rightarrow\infty} \mathcal{I}_{H'} \zeta_n$ in $V'$ is uniquely determined. Consider the linear map
 \begin{equation} \label{eq:id_H'_V'}
 	\mathcal{I}_{V'} : (H^1(\Omega))' \rightarrow V',  \quad \xi \mapsto \big( V \ni  g \mapsto   {}_{(H^1)'}\langle \xi, \TA(\vphi_0) g \rangle_{H^1}  \big),
 \end{equation}
 and note that it is well-defined due to $\TA(\vphi_0) \in \mathcal{L}(H^1(\Omega))$.
Recall that $\TA(\vphi_0)$ is a bijective, bounded, linear operator on $H^1(\Omega)$, cf. Corollary \ref{cor:properties_TA}, and therefore admits  a $C > 0$ such that $\|\TA(\vphi_0) g\|_{H^1} \leq C \norm{g}_{H^1}$ for all $g \in H^1(\Omega)$. Along with the norm identity $\norm{\cdot}_{V} = \norm{\cdot}_{H^1}$ and since $H' \hookrightarrow V'$, we compute
\begin{align*}
	\norm{\mathcal{I}_{H'} \zeta_n - \mathcal{I}_{V'} \xi }_{V'}
	&= \sup_{g \in V\setminus \{ 0 \}} \frac{ \abs{ {}_{(L^2)'}\langle  \zeta_n , \TA(\vphi_0) g \rangle_{L^2}  - {}_{V'}\langle  \zeta_n - \xi, \TA(\vphi_0) g \rangle_{V} }}{\norm{g}_{V}}\\ 
	&= \sup_{g \in H^1(\Omega)\setminus \{ 0 \}} \frac{ \abs{{}_{(H^1)'}\langle  \zeta_n - \xi, \TA(\vphi_0) g \rangle_{H^1} }}{\norm{g}_{H^1}}
	\leq \sup_{g \in H^1(\Omega)\setminus \{ 0 \}}C\, \frac{ \abs{{}_{(H^1)'}\langle  \zeta_n - \xi, \TA(\vphi_0) g \rangle_{H^1} }}{  \| \TA(\vphi_0)g\|_{H^1}}\\ 
	& =  \sup_{g \in H^1(\Omega)\setminus \{ 0 \}} C\, \frac{ \abs{{}_{(H^1)'}\langle  \zeta_n - \xi, g \rangle_{H^1} }}{  \| g\|_{H^1}}\\ 
	& = C \norm{\zeta_n - \xi}_{(H^1)'} \xrightarrow{n \rightarrow \infty} 0, 
\end{align*}
emphasizing the importance of $\TA(\vphi)$ being bijective for the second-to-last equality.  
Thus, the linear map $\mathcal{I}_{V'}$ is the unique extension of $\mathcal{I}_{H'}$. Moreover, exploiting $\TA(\vphi_0)$ being bijective, it is easy to deduce the same for $\mathcal{I}_{V'}$, giving rise to the commutative diagram
\begin{center}
	\begin{tikzcd}
		&H^1(\Omega) \arrow[r, "\cong" above, "\textnormal{id}" below] \arrow[d, "\mathcal{R}_{L^2}(f) = (f{,}\cdot)_{L^2} " left, hook] 
		&V \arrow[d,  "\mathcal{R}_H(f) = (f{,}\cdot)_H " right, hook]
		\\
		&(H^1(\Omega))' \arrow[r, "\mathcal{I}_{V'}" below, "\cong" above] 
		& V' 
	\end{tikzcd}.\\
\end{center}
\par 
\medskip 
Remembering our goal to replace all instances of $p$ and $\bm{u}$ in \eqref{eq:strong_formulation_u}, we compute 
\begin{align*}
	&\pt \theta - \nabla \cdot (\kappa(\vphi) \nabla p ) 	\numberthis \label{eq:strong_w} \\ 
	& \quad = \pt \theta  - \nabla \cdot \big( \kappa(\vphi) \nabla (\TA(\vphi_0) \theta ) \big) 
	- \nabla \cdot \Big( \kappa(\vphi) \nabla \Big( \alpha(\vphi) M(\vphi) \nabla \cdot  \big[ \tilde{\mathcal{C}}^{-1}(\vphi) \big( -\nabla(\C(\vphi) \Tau(\vphi)) + \bm{f }, \bm{g}  \big) \big]   \Big)\Big), 
\end{align*} 
and recall that we are interested in a solution $\theta$ to this equation satisfying the regularity condition $ \theta \in L^s(0, T; H^1(\Omega)) \cap W^{1,s}(0, T; (H^1(\Omega)'))$. By virtue of the isomorphism $\mathcal{I}_{V'}$, it is equivalent to consider the corresponding equation over $V'$ instead of $(H^1(\Omega))'$, such that we may apply \eqref{eq:id_H'_V'} to obtain 
\begin{align*}
	 \mathcal{I}_{V'} \Big({}_{(H^1)'}\langle - \nabla \big( \kappa(\vphi) \nabla (\TA(\vphi) \theta ) \big)  , \cdot  \rangle_{H^1} \Big) 
	 &= {}_{(H^1)'}\langle - \nabla \big( \kappa(\vphi) \nabla (\TA(\vphi) \theta ) \big)  ,  \TA(\vphi_0)\, \cdot  \rangle_{H^1} \\ 
	 & = \big( \kappa(\vphi) \nabla (\TA(\vphi) \theta ), \nabla (  \TA(\vphi_0)\, \cdot) \big)_{L^2} \in V'. 
\end{align*}
Thus, we define the (unbounded) operator $\A(\vphi)$ on the Hilbert space $V'$ as 
\begin{equation*}
	\mathcal{A}(\vphi) : D(\A(\vphi)) \subset V' \rightarrow V', \quad 
	\theta \mapsto  \big( \kappa(\vphi) \nabla (\TA(\vphi) \theta ), \nabla (  \TA(\vphi_0)\, \cdot) \big)_{L^2} ,
\end{equation*}
where
\begin{equation*}
	D(\mathcal{A}(\vphi_0) ) =  \TA^{-1}(\vphi_0) (H^1(\Omega)) = V \subset V'. 
\end{equation*}
\addtocontents{toc}{\SkipTocEntry}
\subsection*{The spaces} The proper spaces in which we want to find solutions are an essential part of a fixed point problem. Therefore, we set 

\begin{align*}
	Z^1_T &\coloneqq W^{1, r} (0, T;  (W^{1, 3}(\Omega))' ) \cap L^r (0, T; W^{3, \frac{3}{2}}_N (\Omega)),
	\\
	Z^2_T &\coloneqq W^{1, s} (0, T; V') \cap L^s (0, T; V), 
\end{align*}
where 
\begin{equation*}
	W^{3, \frac{3}{2}}_N (\Omega) \coloneqq \{ \vphi \in W^{3, \frac{3}{2}}(\Omega) \colon \partial_n \vphi = 0  \textnormal{ on } \partial \Omega \}, 
\end{equation*}
and equip them with the norms $\norm{\cdot}_{Z^1_T}$, $\norm{\cdot}_{Z^2_T}$, respectively, given by 
\begin{align*}
	\norm{\vphi}_{Z^1_T} &\coloneqq \norm{\pt \vphi}_{L^r(0, T; (W^{1, 3}(\Omega))')} + \norm{\vphi}_{L^r(0, T; W^{3, \frac{3}{2}}_N(\Omega))} + \norm{\vphi(0)}_{  (W^{1,3}(\Omega))' , W^{3, \frac{3}{2}}_N(\Omega))_{1- \frac{1}{r}, r}}, \\ 
	\norm{\theta}_{Z^2_T} &\coloneqq \norm{\pt \theta}_{L^{s}(0, T; V')} + \norm{\theta}_{L^{s}(0, T; V)} + \norm{\theta(0)}_{  (V', V )_{1- \frac{1}{s}, s}}. 
\end{align*}
The following result allows us to find uniform estimates in $T  < \frac{T_0}{2}$ and any fixed $T_0 > 0$.

\begin{lemma}\cite[Lem. 2]{abels2021local} \label{lem:uniform_spaces}
	Suppose $0 < T_0 < \infty$ and assume $X_0, X_1$ are some Banach spaces such that $X_0 \hookrightarrow X_1$ densely. We define for every $0 < T < \tfrac{T_0}{2}$ 
	\begin{equation*}
		X_T \coloneqq L^p(0, T; X_1) \cap W^{1, p} (0, T; X_0), 
	\end{equation*}
	where $1 \leq p < \infty$, which shall be equipped with the norm 
	\begin{equation*}
		\norm{u}_{X_T} \coloneqq \norm{u}_{L^p(0, T; X_1)} + \norm{u}_{W^{1, p}(0, T; X_0)} + \norm{u(0)}_{(X_0, X_1)_{1- \frac{1}{p}, p}}. 
	\end{equation*}
	Then, there exists an extension operator $E \colon X_T \rightarrow X_{T_0}$ and some constant $C > 0$ independent of $T$ such that $Eu_{| (0, T)} = u$ in $X_T$ and 
	\begin{equation*}
		\norm{Eu}_{X_{T_0}} \leq C  \norm{u}_{X_T}
	\end{equation*}
	for every $u \in X_T$ and every $0 < T < \tfrac{T_0}{2}$. Moreover, there exists a constant $\tilde{C}(T_0) > 0$ independent of $T$ such that 
	\begin{equation*}
		\norm{u}_{BUC([0, T]; (X_0, X_1)_{1- \frac{1}{p}, p})}  \leq \tilde{C}(T_0) \norm{u}_{X_T}
	\end{equation*}
	for every $u \in X_T$ and every $0 < T < \tfrac{T_0}{2}$. 
\end{lemma}
Finally, we introduce the affine subspaces
\begin{align*}
	X^1_T &\coloneqq \{  \vphi \in Z^1_T \colon \vphi(0) = \vphi_0 \} , \\  
	X^2_T &\coloneqq \{  \theta \in Z^2_T \colon \theta(0) = \theta_0 \}, 
\end{align*}
and $X_T \coloneqq X^1_T \times X^2_T$, \new{which inherits the distance function from $Z_T$}. Similarly, let 
\begin{equation*}
	Y_T \coloneqq Y^1_T  \times Y^2_T \coloneqq L^r(0, T; (H^1(\Omega))')  \times L^{s}(0, T; V'). 
\end{equation*}

\new{
\begin{remark}
	Note that $X_T$ is an affine space that incorporates the initial conditions, and therefore leads to an affine solution operator $\cL^{-1}_T$. Another possible approach that only includes linear spaces and operators is to reduce the problem to the case of homogeneous initial conditions with the help of a so-called reference solution. This would not make much a difference in the following proof. But without subtracting a reference solution the formulae are shorter.   
\end{remark}
}
\addtocontents{toc}{\SkipTocEntry}
\subsection*{The fixed-point problem} 
Linearizing the highest order operators in \eqref{eq:strong_formulation_phi} and \eqref{eq:strong_w}, we define the mapping $\cL_T$ as 
\begin{equation}\label{def:mathcal_A}
		\cL_T : X_T \rightarrow Y_T,  \quad \quad 
		\cL_T (\vphi, \theta) = 
		\left( 
		\begin{array}{c}
			\pt \vphi + \eps \Delta  (m(\vphi_0) \Delta \vphi)  \\ 
			\pt \theta + \A(\vphi_0) \theta 
		\end{array}
		\right). 
\end{equation}
Consequently, the right hand side is given by 
\begin{align*}
	&\F_T^1 (\vphi, \theta) = 
	\begin{aligned}[t]
		&\eps \Delta  (m(\vphi_0)  \Delta \vphi)  - \eps \nabla  \cdot (m(\vphi)  \nabla \Delta \vphi)  \\ 
		& + \nabla \cdot \Big ( m(\vphi) \nabla \Big( \frac{1}{\eps} \psi'(\vphi) + W_{,\vphi}(\vphi, \E(\bm{u}))  \Big) \Big) \\ 
		& + \nabla \cdot \Big ( m(\vphi) \nabla \Big(
		\frac{M'(\vphi)}{2 M(\vphi)} p^2 - p\, \alpha'(\vphi) \nabla \cdot \bm{u}  \Big) \Big) \\ 
		&+ \new{S_s}(\vphi, \E(\bm{u}), \theta), 
	\end{aligned}  \\ 
	&\F_T^2 (\vphi, \theta) 
	\begin{aligned}[t]
		&= \begin{aligned}[t]
			\big(\A(\vphi_0) - \A(\vphi)\big) \theta  
			&+   \nabla \cdot \Big( \kappa(\vphi) \nabla \Big( \alpha(\vphi) M(\vphi) \nabla \cdot  \big[ \tilde{\mathcal{C}}^{-1}(\vphi) \big( -\nabla(\C(\vphi) \Tau(\vphi)) + \bm{f }, \bm{g}  \big) \big]   \Big)\Big)\\ 
			&+ S_f(\vphi, \E(\bm{u}), \theta), 
		\end{aligned}
	\end{aligned}
\end{align*}
where 
\begin{align*}
	\bm{u} & =  \tilde{\mathcal{C}}^{-1}(\vphi) \big( - \nabla \cdot (   \C (\vphi) \Tau(\vphi) + \alpha(\vphi) M(\vphi) \theta  \bm{I} ) + \bm{f}, \bm{g} \big), \\ 
	p &= M(\vphi)(\theta - \alpha(\vphi) \nabla \cdot \bm{u})
\end{align*}
which we combine to the nonlinear mapping $\F : X_T \rightarrow Y_T$ 
\begin{equation*}
	\F_T(\vphi, \theta) = \left(
	\begin{array}{c}
		\F_T^1(\vphi, \theta) \\ 
		\F_T^2 (\vphi, \theta) \\ 
	\end{array}
	\right). 
\end{equation*}

With these definitions in place we can now formulate our last assumption. 

\begin{enumerate}[label = (\roman*)]
	\myitem{(I.A9)} \label{II.A:source_terms}
	The source terms $S_s, S_f$ are functions $\R \times \R^{n \times n} \times \R  \rightarrow \R$. Moreover, we require that for all $(\vphi_i, \theta_i) \in X_T$ satisfying $\norm{(\vphi_i, \theta_i)}_{X_T} \leq R$, $i = 1, 2$ the following Lipschitz-estimates hold  
	\begin{align*}
		\hspace{1.5cm}
		\norm{S_s(\vphi_1, \bm{u}_1,  \theta_1) - S_s(\vphi_2, \bm{u}_2,  \theta_2)}_{L^r(L^2)}
		&\leq C(T, R) ( \norm{\vphi_1 - \vphi_2}_{\new{Z}_T^1} + \norm{\theta_1 - \theta_2}_{\new{Z}_T^2}),
		\\ 
		\norm{S_f(\vphi_1, \bm{u}_1,  \theta_1) - S_f(\vphi_2, \bm{u}_2,  \theta_2)}_{L^s(L^2)}
		&\leq C(T, R) ( \norm{\vphi_1 - \vphi_2}_{\new{Z}_T^1} + \norm{\theta_1 - \theta_2}_{\new{Z}_T^2}),  
	\end{align*}
	where the constant $C (R, T)$ tends to zero as $T$ goes to zero. 
\end{enumerate}
Let us now state the central results used to prove our main theorem. 

\begin{proposition}\label{prop:contraction_non_ve}
	There exists a constant $C(T, R) > 0$ such that for all $(\vphi_i, \theta_i) \in X_T$, $i = 1, 2$, satisfying $\norm{(\vphi_i, \theta_i)}_{X_T} \leq R$ it holds that 
	\begin{equation}\label{eq:contraction_non_ve}
		\norm{\F_T(\vphi_1, \theta_1) - \F_T(\vphi_2, \theta_2)}_{Y_T} \leq C(T, R) \norm{(\vphi_1 - \vphi_2, \theta_1 - \theta_2)}_{\new{Z}_T} 
	\end{equation}
	Moreover, for all $R > 0$, we have $C(R, T) \rightarrow 0$ as $T \rightarrow 0$. 
\end{proposition}

The proof of this proposition is the subject of \new{Section}~\ref{sec:contraction_non_ve}. 

\begin{proposition}\label{prop:inverse_boundedness_non_ve}
	Let $\cL_{T}, X_T$ and $Y_T$ be defined as above. Then, the operator $\cL_T : X_T \rightarrow Y_T$ is invertible for every $T > 0$ and for any fixed $T_0 > 0$ there exists some constant $C_{\mathcal{L}^{-1}_T}(T_0) > 0$ such that 	\new{
	\begin{equation*}
		\norm{ \cL^{-1}_T (f_1, f_2)}_{\new{Z}_T} \leq C_{\mathcal{L}^{-1}_T}(T_0) \Big(   \norm{ (f_1, f_2)}_{Y_T} + \norm{\vphi_0}_{  (W^{1,3}(\Omega))' , W^{3, \frac{3}{2}}_N(\Omega))_{1- \frac{1}{r}, r}} + \norm{\theta_0}_{  (V', V )_{1- \frac{1}{s}, s}} \Big) 
	\end{equation*}
	and 
	\begin{equation*}
		\big\|  \cL^{-1}_T (f_1, f_2) -  \cL^{-1}_T (\tilde{f_1}, \tilde{f_2}) \big\|_{Z_T} \leq 
		C_{\mathcal{L}^{-1}_T}(T_0) \big\| (f_1, f_2) - (\tilde{f_1}, \tilde{f_2}) \big\|_{Y_T}
	\end{equation*}
for all $T \in (0, T_0]$ and  $(f_1, f_2),(\tilde{f_1}, \tilde{f_2})  \in Y_T$. 
}
\end{proposition}

The proof of this proposition can be found in Section~\ref{section:inverse_non_ve}. 
\par 
\medskip 
\noindent 
Given these results, we can prove our main theorem. The arguments below closely follow Abels and Weber \cite{abels2021local}. 

\begin{proof}[Proof of Theorem \ref{thm:existence_strong_nve}] \label{proof:main_visco_elastic_nve}
	Exploiting the invertibility of $\mathcal{L}_T$ from Proposition~\ref{prop:inverse_boundedness_non_ve}, we observe that the operator $\mathcal{F}_T$ and the spaces $X_T, Y_T$ were defined in such a way, that the system given by \eqref{eq:strong_formulation} is equivalent to 
	\begin{equation}\label{eq:fixed_point_eq_nve}
		(\vphi, \theta) = \mathcal{L}^{-1}_T (\mathcal{F}_T ( \vphi, \theta)) \eqqcolon \mathfrak{T}_T(\vphi, \theta)
		\quad \textnormal{in } X_T. 
	\end{equation}
	Hence, the question whether the Cahn--Hilliard--Biot system is locally well-posed in $X_T$ simplifies to a fixed-point problem, which can be solved utilizing a contraction principle. \\ 
	We start by fixing some arbitrary $(\tilde{\vphi}, \tilde{\theta}) \in X_{\tilde{T}}$, for some $\tilde{T} > 0$, and set
	\begin{equation*}
		\new{
		M \coloneqq C_{\mathcal{L}^{-1}_T}(\tilde{T})  
		\Big( \big\|  (\F_{\tilde{T}} (\tilde{\vphi},  \tilde{\theta})) \big\|_{Y_{\tilde{T}}} + \norm{\vphi_0}_{  (W^{1,3}(\Omega))' , W^{3, \frac{3}{2}}_N(\Omega))_{1- \frac{1}{r}, r}} + \norm{\theta_0}_{  (V', V )_{1- \frac{1}{s}, s}} \Big) < \infty. 
	}
	\end{equation*}
	Moreover, we can choose $R > 0$ such that $(\tilde{\vphi},  \tilde{\theta}) \in \overline{B_R^{Z_{\tilde{T}}}(0)} \cap X_{\tilde{T}}$ and $R > 2M$. An application of Proposition~\ref{prop:contraction_non_ve} implies that there exists a family of constants $C = C(T, R) >0$ such that 
	\begin{equation}\label{eq:main_strong_1_nve}
		\norm{\F_T(\vphi_1,  \theta_1) - \F_T(\vphi_2,  \theta_2)}_{Y_T} 
		\leq C(T, R) \norm{(\vphi_1 - \vphi_2, \theta_1 - \theta_2) }_{\new{Z}_T} 
	\end{equation} 
	for all $({\vphi}_i, {\theta}_i) \in X_T$ with $\| (\vphi, \theta)\|_{X_T} < R$, $i = 1, 2$, satisfying $C(T, R) \rightarrow 0$ as $T \rightarrow 0$. \\ 
	Recalling Proposition~\ref{prop:inverse_boundedness_non_ve}, we know that the constant $\new{C_{\mathcal{L}^{-1}_T}(\tilde{T})}$ \new{only depends on $\tilde{T}$ and not on $T \in (0, \tilde{T}]$.} Thus, we can find some $T >0$ sufficiently small such that  
	\begin{equation}\label{eq:main_strong_2_nve}
		\new{C_{\mathcal{L}^{-1}_T}(\tilde{T})} \ C(T, R) <  \frac{1}{4}.		 
	\end{equation}
	Using these facts, we proceed by verifying the assumptions of Banach's fixed point theorem. Since the relations \eqref{eq:main_strong_1_nve} and \eqref{eq:main_strong_2_nve} only hold for functions $({\vphi}_i, {\theta}_i) \in \overline{B_R^{\new{Z}_T}(0)} \cap X_T$, we need to verify that $\mathfrak{T}_T$ defines a self-map of $\overline{B_R^{X_T}(0)} $, i.e., $\mathfrak{T}_T : \overline{B_R^{X_T}(0)}  \rightarrow \overline{B_R^{X_T}(0)}$, \new{where ${B_R^{X_T}(0)} \coloneqq  {B_R^{\new{Z}_T}(0)} \cap X_T$.}\\ 
	By definition of $R > 0$ and the monotonicity of the involved norms with respect to the time variable, the restriction $(\tilde{\vphi},  \tilde{\theta})_{|[0, T]} \in \overline{B_R^{X_T}(0)}$ satisfies 
	\new{
	\begin{align*}
		\norm{\mathfrak{T}_T (\tilde{\vphi},  \tilde{\theta}) }_{Z_T}  &= \big\| \cL^{-1}_T (\F_T(\tilde{\vphi},  \tilde{\theta})) \big\|_{Z_T} \\ 
		&\leq  C_{\mathcal{L}^{-1}_T}(\tilde{T})  
		\Big( \big\|  (\F_{\tilde{T}} (\tilde{\vphi},  \tilde{\theta})) \big\|_{Y_{\tilde{T}}} + \norm{\vphi_0}_{  (W^{1,3}(\Omega))' , W^{3, \frac{3}{2}}_N(\Omega))_{1- \frac{1}{r}, r}} + \norm{\theta_0}_{  (V', V )_{1- \frac{1}{s}, s}} \Big) 
		= M < \frac{R}{2}. 
	\end{align*}
	}
	For any $(\vphi, \theta) \in  \overline{B_R^{X_T}(0)}$, we compute
	\begin{align*}
		\norm{\mathfrak{T}_T (\vphi, \theta)}_{\new{Z}_T}
		&\leq \norm{\mathfrak{T}_T (\vphi, \theta) - \mathfrak{T}_T (\tilde{\vphi},  \tilde{\theta})}_{\new{Z}_T} + \norm{\mathfrak{T}_T (\tilde{\vphi}, \tilde{\theta})}_{\new{Z}_T}  \\ 
		&\new{=  \norm{ \cL^{-1}_T(\F_T (\vphi, \theta)) - \cL^{-1}_T(\F_T (\tilde{\vphi},  \tilde{\theta}))}_{Z_T} + 
		\norm{\mathfrak{T}_T (\tilde{\vphi}, \tilde{\theta})}_{Z_T} } \\ 
		& \leq \new{C_{\mathcal{L}^{-1}_T}(\tilde{T}) }  \norm{\F_T(\vphi, \theta) - \F_T(\tilde{\vphi},  \tilde{\theta})}_{Y_T}   
		+\frac{R}{2} \\ 
		& \leq   \new{C_{\mathcal{L}^{-1}_T}(\tilde{T}) }\  C(T, R) \norm{(\vphi, \theta) - (\tilde{\vphi}, \tilde{\theta})}_{\new{Z_T}}
		+\frac{R}{2} < R, 
	\end{align*}
	where we made use of the Lipschitz property \eqref{eq:main_strong_1_nve}. \new{Since the initial condition is satisfied by the definition of $\cL^{-1}_T$}, this shows the first claim. \par 
	Secondly, we claim that $\mathfrak{T}_T$ is a contraction. Indeed, for any $({\vphi}_i,  {\theta}_i) \in \overline{B_R^{X_T}(0)}$, $i = 1, 2,$ it holds 
	\begin{align*}
		\norm{ \mathfrak{T}_T ({\vphi}_1 , {\theta}_1) - \mathfrak{T}_T  ({\vphi}_2, {\theta}_2) }_{\new{Z_T}}
		&\leq   \new{C_{\mathcal{L}^{-1}_T}(\tilde{T}) }  C(T, R)  \norm{({\vphi}_1 , {\theta}_1)  - ({\vphi}_2,  {\theta}_2)}_{\new{Z}_T} \\ 
		& \leq  \frac{1}{4}  \norm{({\vphi}_1,  {\theta}_1)  - ({\vphi}_2 , {\theta}_2)}_{\new{Z}_T}. 
	\end{align*}
	Thus, the assumptions of Banach's fixed-point theorem are satisfied and we obtain a unique solution $ (\vphi^*, \theta^*) \in  \overline{B_R^{X_T}(0)}$ to the fixed-point problem \eqref{eq:fixed_point_eq_nve}, which is, by construction, a solution to the Cahn--Hilliard--Biot system given in \eqref{eq:strong_formulation}. \par 
	However, so far the solution $(\vphi^*, \theta^*)$ is merely unique in $\overline{B_R^{X_T}(0)}$, and it remains to show uniqueness in $X_T$. To this end, suppose $(\hat{\vphi}, \hat{\theta}) \in X_T$ is another solution. 
	Then, the arguments above give rise to some $\hat{T}  \in (0, T]$ and $\hat{R} \geq R$, such that \eqref{eq:fixed_point_eq_nve} admits a unique solution in $\overline{B_{\hat{R}}^{X_{\tilde{T}}}(0)}$. Since both solutions satisfy these conditions, they must coincide on the interval $[0, \hat{T}]$, i.e., $(\vphi^*, \theta^*)_{| [0, \hat{T}]} = (\hat{\vphi}, \hat{\theta})_{| [0, \hat{T}]}$ and a standard continuation argument shows that then, this already holds on the whole interval $[0, T]$. 
\end{proof}

\medspace
\section{Existence and continuity of $\mathcal{L}^{-1}$: elastic case}\label{section:inverse_non_ve}

As we have seen in the proof of Theorem \ref{thm:existence_strong_nve} above, reducing the question of well-posedness to a fixed-point problem crucially relies on the operator $\mathcal{L}$ being bijective and possessing an uniformly bounded inverse. Taking advantage of the fact that $\cL$ can be split into two independent components, we will establish theses properties for the two respective evolution operators separately. \\
We would like to point out that while the assertion for $\cL_1$ can be derived from existing theory, the problem for $\cL_2$ is much more interesting since this evolution is governed by the sum of two operators. Exploiting that the aggregate corresponds to a dissipative operator on an appropriate Hilbert space $H$, we will deduce that it generates a strongly continuous analytic semigroup over $H$, and therefore exhibits maximal regularity. 
\par 
In this regard, we start with the first component and the following lemma. 

\begin{lemma} \label{lemma:L_1_non_ve}
	Suppose $\Omega \subset \R^n$ is a bounded domain with $C^4$-boundary and assume that \ref{II.A:m} holds. Let $1 < r < \infty$, $\vphi_0  \in \big( (W^{1, 3}(\Omega))', W^{3, \frac{3}{2}}_N(\Omega)\big)_{1- \frac{1}{r}, r}$ and $f \in L^r(0, T; (W^{1, 3}(\Omega))')$, along with coefficients $\tilde{\vphi} \in  W^{2, p^*}(\Omega)$ for some $p^* > \frac{3}{2}$ if $n = 3$ and $p^* \geq \frac{3}{2}$ if $n \leq 2$. Then, for every $0 < T  < \infty$ there exists a unique
	\begin{equation*}
		\vphi \in   L^r (0, T; W^{3, \frac{3}{2}}_N (\Omega)) \cap W^{1, r} (0, T;(W^{1, 3}(\Omega))')
	\end{equation*}
	such that 
	\begin{alignat*}{2}
		\pt \vphi + \eps \Delta  \big( m(\tilde{\vphi}))  \Delta \vphi \big)  &= f  \quad &&\textnormal{in } (0, T) \times \Omega, \numberthis \label{eq:ch_strong_elastic}\\ 
		\vphi_{|t = 0} &= \vphi_0  \quad &&\textnormal{in } \Omega, 
	\end{alignat*}
	where \eqref{eq:ch_strong_elastic} holds in a $(W^{1, 3}(\Omega))'$-sense, i.e., for almost all $t \in (0, T)$ it holds that 
	\begin{equation*}
		{}_{(W^{1, 3})'} \langle \pt \vphi, \zeta \rangle_{W^{1, 3}} 
		-   \int_\Omega \nabla \big( \eps m(\tilde{\varphi }) \Delta \vphi\big) \cdot  \nabla \zeta \dx 
		= {}_{(W^{1, 3})'} \langle f , \zeta \rangle_{W^{1, 3}}  \quad 
		\textnormal{for all } \zeta \in W^{1, 3}(\Omega). 
	\end{equation*}
	 Moreover, for any $T_0 > 0$ there exists some constant $C_{\mathcal{L}^{-1}_{1, T}}(T_0) > 0$ such that \new{
	 	\begin{equation*}
	 		\norm{ \cL^{-1}_{1, T} (f)}_{Z_T^1} \leq C_{\mathcal{L}^{-1}_{1, T}}(T_0) \Big(   \norm{ f}_{Y_T^1} + \norm{\vphi_0}_{  (W^{1,3}(\Omega))' , W^{3, \frac{3}{2}}_N(\Omega))_{1- \frac{1}{r}, r}}  \Big) 
	 	\end{equation*}
	 	and 
	 	\begin{equation*}
	 		\big\|  \cL^{-1}_{1, T} (f) -  \cL^{-1}_{1, T} (\tilde{f}) \big\|_{Z_T^1} \leq 
	 		C_{\mathcal{L}^{-1}_{1, T}}(T_0) \big\| f - \tilde{f} \big\|_{Y_T^1}
	 	\end{equation*}
	 	for all $T \in (0, T_0]$ and  $f, \tilde{f}  \in Y_T^1$. 
	 }
\end{lemma}

\begin{proof}
	We define the operator 
	\begin{equation*}
		\tilde{\mathcal{L}}_1 :  
		D (	\tilde{\mathcal{L}}_1 ) = \{  u \in W^{4, \frac{3}{2}}(\Omega) : \partial_{\bm{n}} u_{|\partial \Omega} =0 = \partial_{\bm{n}} \Delta u_{|\partial \Omega} \}
		\rightarrow 
		L^{\frac{3}{2}}(\Omega), \quad \vphi \mapsto  \eps \Delta \big( m( \tilde{\vphi})  \Delta \vphi \big), 
	\end{equation*}
	and observe that the right-hand side can be expanded as
	\begin{equation*}
		 \tilde{\mathcal{L}}_1 \vphi = \eps m(\tilde{\varphi }) \Delta^2 \vphi 
		+ \eps  m'(\tilde{\varphi }) \Delta \tilde{\vphi} \,  \Delta \vphi 
		+ \eps m''(\tilde{\vphi}) \nabla \tilde{\vphi} \cdot \nabla \tilde{\vphi}\,  \Delta \vphi
		+2 \eps m'(\tilde{\vphi}) \nabla \tilde{\vphi} \cdot \nabla \Delta \vphi .  
	\end{equation*}
	Due to \ref{A:m} and the regularity assumptions on $\tilde{\vphi}$ it is easy to verify that the coefficients satisfy the assumptions of Theorem \ref{thm:H_infty_calculus}. Moreover, $m$ is uniformly positive and real valued, which entails that the principal symbol
	\begin{equation*}
		\mathcal{A}_{\#} (x, \xi) = \sum_{i, j = 1}^n m(\tilde{\varphi }(\bm{x}) ) \, \xi_i^{2} \xi_j^2
	\end{equation*}
	is uniformly parameter elliptic in $\bm{x} \in \overline{\Omega}$ with arbitrarily small angle of ellipticity $0 < \phi_{\mathcal{A}_\#} < \tfrac{\pi}{2}$. 
	By invoking Theorem \ref{thm:H_infty_calculus}, we infer that for any $\phi < \tfrac{\pi}{2}$ there exists some $\mu_1 \geq 0$ such that this operator admits $\tilde{\mathcal{L}}_1 + \mu_1 \in \mathcal{H}^\infty(L^{\frac{3}{2}}(\Omega))$ with $\phi^\infty_{\tilde{\mathcal{L}}_1 + \mu_1 } \leq \phi$. In particular, this operator has bounded imaginary powers with power angle $\theta_{\tilde{\mathcal{L}}_1 + \mu_1} \new{\leq} \phi^\infty_{\tilde{\mathcal{L}}_1 + \mu_1 } < \tfrac{\pi}{2}$, i.e., $\tilde{\mathcal{L}}_1 + \mu_1 \in \mathcal{BIP}(L^\frac{3}{2}(\Omega))$, implying that the extrapolated fractional power scale of order $m$ generated by $(L^\frac{3}{2}(\Omega), \tilde{\mathcal{L}}_1 + \mu_1 ) \eqqcolon (E, \prescript{}{\frac{3}{2}}{A})$ and the interpolation-extrapolation scale $[(E_\alpha,  \prescript{}{p}{A}_\alpha)\ ;\ \alpha \in [-m, \infty )]$ generated by $(E, \prescript{}{\frac{3}{2}}{A} )$ and the complex interpolation operator $[\cdot, \cdot]_{\vartheta}$ are equivalent, cf.\ Section~\ref{subsec:scales}, \cite[Theorem V 1.5.4]{amann1995linear}. We are interested in the case $\alpha = -\tfrac{1}{4}$, for which we compute 
	\begin{equation*}
		\prescript{}{\frac{3}{2}}{A}_{- \frac{1}{4}} : D( \prescript{}{\frac{3}{2}}{A}_{- \frac{1}{4}}) = E_{\frac{3}{4}} \subset  E_{- \frac{1}{4}}\rightarrow  E_{- \frac{1}{4}}, 
	\end{equation*}
	where $ \prescript{}{p}{A}_{- \frac{1}{4}}$ is the $E_{- \frac{1}{4}}$-realization of $ \prescript{}{p}{A} = \tilde{\mathcal{L}}_1 + \mu_1$. Exploiting the fact that the fractional power scale and the interpolation-extrapolation scale are equivalent, we find 
	\begin{align*}
		E_{\frac{3}{4}} = [E_0, E_1]_{\frac{3}{4}} = [L^\frac{3}{2}(\Omega), D(A)]_{\frac{3}{4}} 
		&= [L^\frac{3}{2}(\Omega),\{  u \in W^{4, \frac{3}{2}}(\Omega) : \partial_n \vphi_{|\partial \Omega} =0 = \partial_n \Delta \vphi_{|\partial \Omega} \}]_{\frac{3}{4}} \\ 
		&= \{  u \in W^{3, \frac{3}{2}}(\Omega) : \partial_n \vphi_{|\partial \Omega} = 0 \} = W^{3, \frac{3}{2}}_N(\Omega),  
	\end{align*}
	where the last identity holds due to results by Seeley \cite{seeley1972interpolation}, see also \cite[VII 2.4]{amann2013analysis}.  \\ 
	To compute $E_{- \frac{3}{4}}$, note that by \cite[Theorem V 1.4.12]{amann1995linear}, see also Section~\ref{subsec:scales}, it holds for $\alpha \in \R$ that $(E_\alpha)' = E^\#_{-\alpha}$ with respect to the duality paring induced by $\langle \cdot , \cdot \rangle_{E, E^\#}$, where $E^\#_{-\alpha}$ is part of the extrapolated fractional power scale $[(E_\alpha^\#,  \prescript{}{\frac{3}{2}}{A}_\alpha^\#) \ ; \ \alpha \in [-m, \infty )]$ generated by the dual pair $(E^\#,  \prescript{}{\frac{3}{2}}{A}^\#)$. An easy computations further shows the dual operator to $ \prescript{}{\frac{3}{2}}{A}$ is given by 
	\begin{align*}
		\prescript{}{\frac{3}{2}}{A}^\# &= \prescript{}{3}{A} = \tilde{\mathcal{L}}_1 + \mu_1 :  
		D (	\prescript{}{3}{A} ) = \{  u \in W^{4, 3}(\Omega) : \partial_n \vphi_{|\partial \Omega} =0 = \partial_n \Delta \vphi_{|\partial \Omega} \} \subset L^{3}(\Omega)
		\rightarrow 
		L^{3}(\Omega),  \\
		\vphi &\mapsto  \eps \Delta  \big( m(\tilde{\varphi })  \Delta \vphi \big) + \mu_1 \vphi. 
	\end{align*}
	Since $	\prescript{}{3}{A}_\alpha \in \mathcal{BIP}(L^{3}(\Omega))$ with $\theta_{\prescript{}{p'}{A}} < \tfrac{\pi}{2}$, cf. \cite[Thm. V 1.4.11]{amann1995linear}, we can argue as above and see $E^\#_{\frac{1}{4}} = [E^\#_0, E^\#_1 ]_{\frac{1}{4}}$. Thus, a similar computation shows 
	\begin{align*}
		E_{-\frac{1}{4}}  = (E^\#_{\frac{1}{4}})' 
		&= [(E^\#_0),( E^\#_1)]'_{\frac{1}{4}} 
		=  [L^{3}(\Omega),\{  u \in W^{4, 3}(\Omega) : \partial_n \vphi_{|\partial \Omega} =0 = \partial_n \Delta \vphi_{|\partial \Omega} \}]'_{\frac{1}{4}}
		= (W^{1, 3}(\Omega))' 
	\end{align*}
	and we arrive at 
	\begin{equation*}
		\prescript{}{\frac{3}{2}}{A}_{- \frac{1}{4}} :  W^{3, \frac{3}{2}}_N(\Omega) \subset  (W^{1, 3}(\Omega))' \rightarrow  (W^{1, 3}(\Omega))',  
	\end{equation*}
	noting that this operator also exhibits bounded imaginary powers with power angle $\theta_{	\prescript{}{\frac{3}{2}}{A}_{- \frac{1}{4}}} < \tfrac{\pi}{2}$ (see \cite[Thm. V 1.5.5]{amann1995linear}), i.e., $\prescript{}{\frac{3}{2}}{A}_{- \frac{1}{4}} \in \mathcal{BIP}(  (W^{1, 3}(\Omega))')$. \\ 
	Finally, consider the operator 
		\begin{equation*}
			\hat{\mathcal{L}}_1 :  
			W^{3, \frac{3}{2}}_N(\Omega) \subset  (W^{1, 3}(\Omega))' \rightarrow  (W^{1, 3}(\Omega))', \quad \vphi \mapsto  \eps \Delta  \big( m(\tilde{\vphi})  \Delta \vphi \big), 
		\end{equation*} 
		and observe that we need to distinguish between two separate cases: if $\mu_1 = 0$ it already holds that $\hat{\mathcal{L}_1}  = { _{\frac{3}{2}} }A_{-\frac{1}{4}} \in  \mathcal{BIP}( (W^{1, 3}(\Omega))')$, and since the power angle $\theta_{\hat{\mathcal{L}_1}}$ is less than $\frac{\pi}{2}$, we even obtain $L^r$-maximal regularity on $\R_+$ by Theorems~\ref{thm:bip->r_sectorial}-\ref{thm:r_sec_max_reg}.\\ 
		Otherwise, we invoke Theorem~\ref{thm:max_reg_semigroups}, which entails maximal regularity for the problem  
		\begin{equation*}
			\pt \vphi + 	\hat{\mathcal{L}_1} \vphi = f(t), \quad t \in (0, T), \quad \vphi(0) = 0, 
		\end{equation*}
		where $f \in L^r(0, T;  (W^{1, 3}(\Omega))')$. Using a standard argument, cf. \cite[Prop.\ 1.3]{ARENDT20071}, \cite[Thm.\ 10]{haselboeck2024existence}, we find that the corresponding Cauchy problem with non-trivial initial data $\vphi_0 \in \big((W^{1, 3}(\Omega))', W^{3, \frac{3}{2}}_N(\Omega)\big)_{1- \frac{1}{r}, r}$ is also well-posed. \\ 
		The second assertion can be shown exactly as in Lemma \ref{lem:L^-1_uniform_T}.
\end{proof}

We would also like to refer to \cite[Thm.\ 6.3.6]{PrussSimonett2016} where the authors prove a similar result for an elliptic operator in non-divergence form and in a more general setting. Let us proceed with maximal regularity for the second component which requires a more subtle approach.

\begin{proposition}\label{prop:analytic_generator}
	Let $V \xhookrightarrow{d} H \cong H' \xhookrightarrow{} V'$ be a Gelfand triplet over $\R$-Hilbert spaces and suppose that the unbounded operator $A : D(A)= V \subset V' \rightarrow V'$ satisfies  
	\begin{enumerate}[label*=(\roman*), noitemsep]
		\item ${}_{V'}\langle Au, v \rangle_{V}  = {}_{V'}\langle Av, u \rangle_{V}$ for all $u, v \in V$;
		\item ${}_{V'}\langle Au, u \rangle_{V}  \geq \alpha \norm{u}_{V}^2 - \beta \norm{u}_{H}^2$ for some $\alpha> 0,\ \beta \geq 0$ and all $u \in V$. 
	\end{enumerate}
	Then there exists a Hilbert space $\hat{V} \cong V$ such that the realization of the operator $-A$ on this space $- \hat{A} : D(\hat{A}) = \hat{V} \subset \hat{V}' \rightarrow \hat{V}'$ generates an analytic semigroup. 
\end{proposition}

\begin{proof}
For now, consider the operator $A_\beta \coloneqq A + \beta\mathcal{I}$. 
	We claim that for all $\lambda \geq 0$, the operator 
	\begin{equation*}
		\lambda \mathcal{I} - (- A_\beta)   :  D(A) = V  \subset V' \rightarrow V'
	\end{equation*}
	is bijective. Indeed, we can define the mappings 
	\begin{align*}
		(u, v) \mapsto {}_{V'}\langle (\lambda + A_\beta) u, v \rangle_{V} = {}_{V'}\langle (\lambda + A + \beta) u, v \rangle_{V}, \quad \textnormal{where}\quad u, v \in V, 
	\end{align*}
	and observe that we obtain bounded, bilinear forms on $V$. Stressing that the canonical embedding $V \hookrightarrow V'$ in a Gelfand triplet exploits the Riesz isomorphism of the Hilbert space $H$, we apply $(ii)$ to find
	\begin{align} \label{eq:a_b_coercive}
		{}_{V'}\langle (\lambda + A + \beta) u, u \rangle_{V} 
		\geq  (\lambda + \beta) (u, u)_{H} + \alpha \norm{u}^2_V - \beta \norm{u}^2_H 
		\geq \alpha \norm{u}^2_V
	\end{align}
	for all $u \in V$ and $\lambda \geq 0$, i.e., the bilinear forms are coercive. 
	\par 
	\medskip 
	In particular, the mapping
	\begin{equation} \label{eq:V_hat_inner}
		V \times V \ni (u, v) \mapsto (u, v)_{\hat{V}} \coloneqq  {}_{V'}\langle A_\beta u, v \rangle_{V}  \in \R
 	\end{equation}
	is a symmetric, positive definite, bilinear form and therefore gives rise to an inner product on $V$. To establish the isomorphism $(\hat{V}, (\cdot, \cdot)_{\hat{V}}^{\sfrac{1}{2}}) \cong V$, it suffices to show that there exist constants $C \geq c > 0$ such that 
	\begin{equation}\label{eq:V_hat_norm_equiv}
		c \norm{u}_{V} \leq \norm{u}_{\hat{V}} \leq C \norm{u}_{V} \quad \textnormal{for all} \quad u \in V. 
	\end{equation}
	Choosing $c = \alpha$, we obtain the first inequality immediately from \eqref{eq:a_b_coercive}, and observing that $A_\beta \in \mathcal{L}(V, V')$ is a bounded operator yields the assertion. Thus, we can consider the Gelfand triple $\hat{V} \xhookrightarrow{d} H \cong H' \xhookrightarrow{} \hat{V}'$ and on it the unbounded operator $\hat{A}_\beta : D(\hat{A}) = \hat{V} \subset \hat{V}' \rightarrow \hat{V}'$, which is the realization of $A_\beta$ on $\hat{V}$. The reason to consider $\hat{V}$ instead of $V$ will become evident when we show that $\hat{A}_\beta$ is self-adjoint, where symmetry is crucial. 
	\par 
	\medskip
	\textit{Claim: $[0, \infty) \subset \rho(\hat{A}_\beta)$:} For any $\lambda \geq 0$, we deduce from \eqref{eq:a_b_coercive} and \eqref{eq:V_hat_norm_equiv} 
	\begin{equation}\label{eq:A_hat_coercive}
			{}_{\hat{V}'}\langle (\lambda + \hat{A}_\beta) u, u \rangle_{\hat{V}} 
		\geq \alpha \norm{u}^2_V
		\geq c \norm{u}^2_{\hat{V}}, 
	\end{equation}
	such that an application of the Lax-Milgram theorem entails that for every $\Phi \in \hat{V}'$ there exists a unique $u = u_{\Phi} \in \hat{V}$ such that 
	\begin{equation*}
		{}_{\hat{V}'}\langle (\lambda + A + \beta) u, v \rangle_{\hat{V}} = \Phi(v) \quad \textnormal{for all } \quad v \in \hat{V}, 
	\end{equation*}
	or equivalently that the operators $(\lambda -(- A_\beta)) : \hat{V} \rightarrow \hat{V}'$ are bijective. Hence, $[0, \infty) \subset \rho(-\hat{A}_\beta)$, where $\rho (- \hat{A}_\beta)$ denotes the resolvent set. 
	\par 
	\medskip 
	 \textit{Closed:} Since $-(\hat{A}_\beta)$ is a densely defined operator whose resolvent set is nonempty, it follows that this operator $-(\hat{A}_\beta)$ is closed, cf. \cite[VII 5.32]{werner2018funktionalanalysis}. 
	 \par
	 \medskip 
	 \textit{Self-adjoint:} To establish that $-\hat{A}_\beta$ is self adjoint, we need to show that $\hat{A}_\beta$ is symmetric on $\hat{V}'$, i.e., 
	 \begin{equation*}
	 	(\hat{A}_\beta u, v)_{\hat{V}'} = 	(u, \hat{A}_\beta  v)_{\hat{V}'}, \quad \textnormal{for all} \quad u, v \in \hat{V} = D(\hat{A}_\beta).
	 \end{equation*}
	 We recall that the inner product on $\hat{V}'$ is defined through the Riesz isomorphism by 
	 \begin{equation*}
	 		(\hat{A}_\beta u, v)_{\hat{V}'} =  (\hat{A}_\beta u, \mathcal{R}_H v)_{\hat{V}'} \coloneqq 
	 		 ( \mathcal{R}^{-1}_{\hat{V}} \hat{A}_\beta u,  \mathcal{R}^{-1}_{\hat{V}} \mathcal{R}_H v)_{\hat{V}},
	 \end{equation*}
	 where due to \eqref{eq:V_hat_inner} we have $ \mathcal{R}^{-1}_{\hat{V}} \hat{A}_\beta u = u$, leading to 
	 \begin{align*}
	 	(\hat{A}_\beta u, v)_{\hat{V}'} 
	 	= ( u,  \mathcal{R}^{-1}_{\hat{V}} \mathcal{R}_H v)_{\hat{V}}
	 	= (u, v)_{H}
	 	= (  \mathcal{R}^{-1}_{\hat{V}} \mathcal{R}_H u, v)_{\hat{V}}
	 	= (u, \hat{A}_\beta v)_{\hat{V}'}. 
	 \end{align*}
 	Hence, the negative operator $-\hat{A}_\beta$ is also symmetric. To conclude we exploit the well-known fact that densely defined, closed and symmetric operators on Hilbert spaces whose resolvent set has a non-empty intersection with the real line (recall that $[0, \infty) \subset  \rho (- \hat{A}_\beta)$) are self-adjoint.\\ 
 	Since the spectrum of a self adjoint operator is always a subset of $\R$, it follows that $\sigma(-\hat{A}_\beta) \subseteq \R_{< 0}$. 
	 \par 
	 \medskip 	 
	 \textit{Analytic semigroup:} 
	 To infer that $\hat{A}_\beta$ generates an analytic semigoup on $\hat{V}'$, 
	 we need to additionally verify that certain resolvent estimates similar to \ref{eq:resolvent_estimates} hold, see \cite[Def.\ 4.1, Thm.\ 4.6 ]{engel2006one}.  \par
	 Relying on profound insights by Weis or Taylor and Lay relating the norm of the resolvent to their spectral range, it was shown in \cite[Cor. 4.7]{engel2006one} that for any normal operator $\mathfrak{A}$ on a Hilbert space satisfying $\sigma(\mathfrak{A}) \subseteq \{ z \in \new{\mathbb{C}} \colon \arg(-z ) < \delta \}$, \new{$\delta \in [0, \frac{\pi}{2})$}, the respective estimates hold, entailing that $\mathfrak{A}$ generates an analytic semigroup. Observing that $-\hat{A}_\beta$ meets these conditions, we find that $- \hat{A}_\beta$ is the generator of an analytic semigroup on $\hat{V}'$.
	 \par 
	 \medskip 	 
	 Finally, we deduce the same for $-\hat{A}$ since shifting the generator of an analytic semigroup by scalar multiples of the identity retains this property, see Theorem \ref{thm:perturbation_semigroup}. 
\end{proof}

As before, and in accordance with \eqref{eq:V_hat_inner}, we can define the bilinear mapping  
\begin{equation*}
	(u, v) \mapsto ( u, v )_{\tilde{V}} \coloneqq
	\int_\Omega \nabla (\TA(\vphi_0) u)  \cdot \nabla (\TA(\vphi_0) v) +  u\, \TA(\vphi_0)v  \dx \quad \textnormal{for }u, v  \in H^1(\Omega), 
\end{equation*}
and note that it is symmetric and positive definite due to $\TA(\vphi_0)$ being positive and self adjoint on $L^2(\Omega)$ and bijective on $H^1(\Omega)$, cf. Corollary~\ref{cor:properties_TA}. In particular, $(u, v)_{\tilde{V}} $ defines an inner product on $H^1(\Omega)$. 

\begin{lemma}\label{lem:norm_equivalence_TH}
	Let $(\tilde{V}, \norm{\cdot}_{\tilde{V}}) \coloneqq \big(H^1(\Omega), ( \cdot , \cdot )_{\tilde{V}} ^{\frac{1}{2}} \big)$, then $\tilde{V}$ is a Hilbert space and $\tilde{V} \cong H^1(\Omega)$. 
\end{lemma}

While this result can be shown solely with estimates from the previous chapter we will shortly discover that the assertion also follows from Proposition \ref{prop:analytic_generator}. \par \medskip 

Finally we will apply Proposition \ref{prop:analytic_generator} to show that $\mathcal{A}(\vphi_0)$ indeed generates an analytic semigroup over some Hilbert space $\tilde{V}' \cong V'$, thus exhibiting maximal regularity. We stress that by suitably choosing the space $H$ in our Gelfand triple, and thereby the canonical Riesz ismorphism, our original problem was  identified with an equivalent problem over a Hilbert space that is isomorphic to $H^1(\Omega)$, such that the realization of the operator $\mathcal{A}(\vphi_0)$ fulfills the assumptions of Proposition \ref{prop:analytic_generator}.

\begin{lemma} \label{lemma:L_2_non_ve}
	Let $1 < s < \infty$ and assume  $\theta_0  \in \big(V', D(\A(\vphi_0))\big)_{1- \frac{1}{s}, s}$ and $f \in L^s(0, T; V')$. Then, for every $0 < T  < \infty$ there exists a unique
	\begin{equation*}
		\theta \in   L^s (0, T; V) \cap W^{1, s} (0, T; V')
	\end{equation*}
	such that 
	\begin{equation}\label{eq:cauchy_theta_non_ve}
			\begin{alignedat}{2}
			\pt \theta +  \A(\vphi_0)\theta  &= f  \quad &&\textnormal{in } (0, T) \times \Omega, \\ 
			\theta(0) &= \theta_0  \quad &&\textnormal{in }\{ 0\} \times \Omega. 
		\end{alignedat}
	\end{equation}
	Moreover, for any $T_0 > 0$ there exists some constant $C_{\mathcal{L}^{-1}_{2, T}}(T_0) > 0$ such that \new{
		\begin{equation*}
			\norm{ \cL^{-1}_{2, T} (f)}_{Z_T^2} \leq C_{\mathcal{L}^{-1}_{2, T}}(T_0) \Big(   \norm{ f}_{Y_T^2} + \norm{\theta_0}_{  (V', V )_{1- \frac{1}{s}, s}} \Big) 
		\end{equation*}
		and 
		\begin{equation*}
			\big\|  \cL^{-1}_{2, T} (f) -  \cL^{-1}_{2, T} (\tilde{f}) \big\|_{Z^2_T} \leq 
			C_{\mathcal{L}^{-1}_{2, T}}(T_0) \big\| f - \tilde{f} \big\|_{Y_T^2}
		\end{equation*}
		for all $T \in (0, T_0]$ and  $f, \tilde{f}  \in Y_T^2$. 
	}
\end{lemma}

\begin{proof}
	We consider the Gelfand triplet $ V \xhookrightarrow{d} H \cong H' \hookrightarrow V'$ and observe that the unbounded operator
	\begin{equation*}
		\mathcal{A}(\vphi_0) : D(\A(\vphi_0)) \subset V' \rightarrow V', \quad 
		\theta \mapsto  \big( \kappa(\vphi_0) \nabla (\TA(\vphi_0) \theta ), \nabla (  \TA(\vphi_0)\, \cdot) \big)_{L^2} 
	\end{equation*}
	satisfies two properties:
	Firstly, it holds for any $u, v \in V$ that 
	\begin{align*}
		{}_{V'}\langle A(\vphi)u, v \rangle_{V} 
		&= \big( \kappa(\vphi_0) \nabla (\TA(\vphi_0) u ), \nabla (  \TA(\vphi_0)\, v) \big)_{L^2}\\
		&=  \big( \kappa(\vphi_0) \nabla (\TA(\vphi_0) v), \nabla (  \TA(\vphi_0)\, u) \big)_{L^2}
		= {}_{V'}\langle A(\vphi)v, u \rangle_{V}. 
	\end{align*}
	Secondly, using \ref{A:kappa} and Lemma~\ref{lem:norm_equivalence_TH}, we estimate for every $u \in V= H^1(\Omega)$
	\begin{align*}
		{}_{V'}\langle A(\vphi)u, u \rangle_{V} 
		&= \big( \kappa(\vphi_0) \nabla (\TA(\vphi_0) u ), \nabla (  \TA(\vphi_0)\, u) \big)_{L^2}
		\geq \underline{\kappa} \norm{\nabla \TA(\vphi_0) u  }_{L^2}^2\\ 
		&\geq   \underline{\kappa} \big( \norm{u}_{\tilde{V}}^2 - \norm{u}_H^2  \big)
		\geq \alpha \norm{u}_{V}^2 - \beta \norm{u}^2_H. 
	\end{align*}
	Hence, the assumptions of Proposition \ref{prop:analytic_generator} are satisfied and we conclude that $-\mathcal{A}(\vphi_0)$ generates a strongly continuous analytic semigroup on some Hilbert space $\hat{V}'  \cong V'$. However, looking more closely, it is evident from \eqref{eq:V_hat_inner} that if the choice $\beta = 1$ is possible, we even have  $\hat{V}' = \tilde{V}'$.  \\ 
	By Theorem~\ref{thm:max_reg_analytic_hilbert} we can therefore infer the existence of a unique $\theta \in L^s(0, T; \tilde{V}) \cap W^{1, s}(0, T; \tilde{V}')$ solving \eqref{eq:cauchy_theta_non_ve}, and the isomorphism between the Hilbert spaces yields the desired assertion.  \\
	Lastly, the second part can be shown exactly as in Lemma \ref{lem:L^-1_uniform_T}. 
\end{proof}

\medspace 

\section{Lipschitz continuity of $\mathcal{F}$: elastic case} \label{sec:contraction_non_ve}

The goal of this section is to establish the Lipschitz estimates postulated in Proposition \ref{prop:contraction_non_ve}. We start by applying the interpolation results from Section \ref{sec:preliminaries} to deduce some important embedding properties that will become very useful in the following computations. 
\medskip 
\par 
From the definition of $X^1_T$ and \eqref{int:BUC}, it follows that 
\begin{align}\label{embedding_X^1_non_ve_1}
	X^1_T &\hookrightarrow BUC ([0, T];  ((W^{1, 3}(\Omega))',W^{3, \sfrac{3}{2}}_N(\Omega))_{1- \frac{1}{r}, r})  \hookrightarrow  BUC ([0, T]; B^{3- \frac{4}{r}}_{\sfrac{3}{2}, r} (\Omega)). 
\end{align}
Note that for $r > 4$, if $n = 3$, it holds that 
\begin{equation*}
	B^{3- \frac{4}{r}}_{\sfrac{3}{2}, r} (\Omega) \hookrightarrow W^{1, q}(\Omega),  
\end{equation*}
for some $q > 3$. 
If $n \leq 2$, we can relax this assumption, merely requiring $r > \frac{12}{5}$ to get 
\begin{equation*}
		B^{3- \frac{4}{r}}_{\sfrac{3}{2}, r} (\Omega) \hookrightarrow W^{1, q} (\Omega), 
\end{equation*}
where $q > 2$. 
Hence there exists some $q > n$, depending on the dimension, such that 
\begin{equation}\label{embedding_X^1_non_ve}
	\begin{aligned}
		X^1_T \hookrightarrow BUC ([0, T]; W^{1, q} (\Omega)). 
	\end{aligned}
\end{equation}
Additionally, it should be noted that for some $\epsilon > 0$, $\vartheta \in (0, 1)$ standard interpolation results, cf. \cite{runst2011sobolev}, lead to the inequality
\begin{equation*}
\norm{\vphi}_{W^{\frac{2n + \epsilon}{3}, \sfrac{3}{2}} } 
	\leq C \norm{ \vphi}_{B^{3 - \frac{4}{r}}_{\sfrac{3}{2}, r} }^{1 - \vartheta} \norm{\vphi}_{(W^{1,3})'}^\vartheta . 
\end{equation*}
Along with the embedding $W^{1, r}(0, T;(W^{1, 3}(\Omega))') \hookrightarrow C^{1 - \frac{1}{r}}([0; T];(W^{1, 3}(\Omega))')$ and \eqref{embedding_X^1_non_ve_1}, we find that Lemma~\ref{lemma:interpolation_hoelder} is applicable, concluding that 
\begin{equation}\label{embedding:hoelder_non_ve}
	X^1_T \hookrightarrow C^{(1 - \frac{1}{r}) \vartheta}([0, T];W^{\frac{2n + \epsilon}{3}, \frac{3}{2}}(\Omega)) 
		\coloneqq C^{ \tilde{\vartheta}}([0, T]; W^{\frac{2n + \epsilon}{3}, \frac{3}{2}}(\Omega))
			 \hookrightarrow C^{ \tilde{\vartheta}}([0, T]; W^{1, q}(\Omega)) 
\end{equation}
for $\tilde{\vartheta} \coloneqq (1 - \frac{1}{r}) \vartheta > 0$ and some $q > n$. Note that $\vartheta$ is indirectly proportional to $r$, i.e., choosing a bigger exponent $r$ leads to better Hölder regularity $\bar{\vartheta}$. \\ 

We proceed by deriving an estimate, which will be quite useful later on.  

\begin{lemma}\label{lemma:general_estimate_non_ve}
	Let $f_i \in W^{1, \sfrac{3}{2}}(\Omega)$ and $ \vphi_i \in W^{1, q}(\Omega)$, where $q > n$, with $\norm{\vphi_i}_{W^{1, q}}  \leq R $ for some $R > 0$ and $i = 1, 2$. Then, it holds for any $\tilde{m} \in C^1(\R)$ that 
	\begin{align*}
		\norm{\nabla \cdot \big(\tilde{m}(\vphi_1) \nabla f_1 \big) -  \nabla \cdot \big(\tilde{m}(\vphi_2) \nabla f_2\big) }_{(W^{1, 3})'}  
		 \leq C (R)  \Big(  \norm{f_1- f_2}_{W^{1, \sfrac{3}{2}}} + \norm{f_2}_{W^{1, \sfrac{3}{2}}} \norm{\vphi_1 - \vphi_2 }_{W^{1, q}}  \Big). 
	\end{align*}
\end{lemma}

\begin{proof}
	Using the fact that $\norm{\nabla \cdot (\cdot)}_{ (W^{1, 3})'}  \leq C \norm{\cdot}_{L^{\sfrac{3}{2}} } $, we compute 
	\begin{align*}
		&\norm{\nabla \cdot \big(\tilde{m}(\vphi_1) \nabla f_1 \big) -  \nabla \cdot \big(\tilde{m}(\vphi_2) \nabla  f_2\big) }_{(W^{1, 3})'} 
		\\
		& \quad 
		\leq  C \norm{ \tilde{m}(\vphi_1) \nabla \big( f_1 - f_2\big) }_{L^{\sfrac{3}{2}}}   + C  \norm{ \big(\tilde{m}(\vphi_1 - \tilde{m}(\vphi_2) \big) \nabla f_2) }_{L^{\sfrac{3}{2}}} \\ 
		& \quad 
		\leq C \norm{ \tilde{m}(\vphi_1) }_{L^\infty} \norm{  \nabla \big( f_1 - f_2\big) }_{L^{\sfrac{3}{2}} }
				+ C \norm{\nabla f_2}_{L^{\sfrac{3}{2}}}  \norm{\tilde{m}(\vphi_1) - \tilde{m}(\vphi_2) }_{L^\infty} . 
	\end{align*}
	Utilizing $\tilde{m} \in C^1(\R)$, we deduce that $\tilde{m}$ is bounded on compact sets and locally Lipschitz continuous. The embedding $W^{1, q}(\Omega) \hookrightarrow C^0(\overline{\Omega})$ therefore yields
	\begin{equation*}
		\norm{\tilde{m}(\vphi_1)}_{L^\infty} \leq C(R) \quad \textrm{and} 
		\quad \norm{\tilde{m}(\vphi_1) - \tilde{m}(\vphi_2)}_{L^\infty} \leq L(R) \norm{\vphi_1  -  \vphi_2}_{W^{1, q} }
	\end{equation*} 
	for all $\norm{\vphi_i}_{W^{1, q}}  \leq R$, $i = 1, 2$. 
	Hence, 
	\begin{equation*}
		\norm{\nabla \cdot \big(\tilde{m}(\vphi_1) \nabla f_1 \big) -  \nabla \cdot \big(\tilde{m}(\vphi_2) f_2\big) }_{(W^{1, 3})'}  
		\leq C (R)  \Big(  \norm{  f_1 - f_2 }_{W^{1, \sfrac{3}{2}}}
		+ \norm{ f_2}_{W^{1, \sfrac{3}{2}}} \norm{\vphi_1 - \vphi_2 }_{W^{1, q}}  \Big), 
	\end{equation*}
	which concludes the proof. 
\end{proof}

Before we start with the proof of Proposition \ref{prop:contraction_non_ve}, let us first examine the displacement $\bm{u}$.

\begin{lemma}\label{lem:contraction_u_non_ve}
	If $\vphi\in  L^\infty( [0, T]; W^{1, q} (\Omega)) $, $q > n$, with $\norm{\vphi}_{L^\infty (W^{1, q})} \leq R$, and $\theta\in L^s(0, T; H^1(\Omega))$, then the function 
	\begin{equation*}
		\bm{u} =  \tilde{\mathcal{C}}^{-1}(\vphi) \big( - \nabla \cdot (\C (\vphi) \Tau(\vphi) + \alpha(\vphi) M(\vphi) \theta \bm{I} ) + \bm{f}, \bm{g} \big)
	\end{equation*}
	belongs to  $L^s(0, T; \bm{H}^2(\Omega))$ and is bounded by a constant $C = C(R)$. Moreover, it holds for almost every $t \in (0, T)$ and all $\vphi_i \in  W^{1, q} (\Omega) $ with $\norm{\vphi_i}_{W^{1, q}} \leq R$, $i = 1, 2$, that 
	\begin{equation*}
		\norm{\bm{u}_1 - \bm{u}_2}_{\bm{H}^2} 
		\leq C(R) \norm{\vphi_1 - \vphi_2}_{W^{1, q}} \big( 1+ \norm{\theta_2}_{H^1}+ \norm{\bm{f}}_{\bm{L}^2} + \norm{ \bm{g}}_{\bm{H}^{\sfrac{1}{2}}(\Gamma_N)}\big) + C(R) \norm{\theta_1 - \theta_2 }_{H^1}. 
	\end{equation*} 
\end{lemma}

\begin{proof}
	Quick computations show that $\norm{f\,g}_{H^1} \leq C \norm{f}_{H^1} \norm{g}_{W^{1, q}}$ and $\norm{f\,g}_{W^{1, q}} \leq C \norm{f}_{W^{1, q}} \norm{g}_{W^{1, q}}$, such that along with Theorem~\ref{thm:compositon} we can estimate
	\begin{align*}
		&\norm{- \nabla \cdot (\C (\vphi) \Tau(\vphi) + \alpha(\vphi) M(\vphi) \theta  \bm{I} ) + \bm{f}}_{\bm{L}^2}
		\\
		& \quad  \leq C \big( \norm{\C (\vphi)}_{\bm{W}^{1, q}} \norm{\Tau(\vphi)}_{\bm{H}^1} +
		 \norm{ \alpha(\vphi) }_{W^{1, q}}  \norm{ M(\vphi)}_{W^{1, q}}\norm{ \theta}_{H^1} \big) + \norm{\bm{f}}_{\bm{L}^2}\\
		&\quad \leq C(R) \big( \norm{\vphi}_{W^{1, q}} + \norm{\vphi}^2_{W^{1, q}}  \norm{ \theta}_{H^1} \big) +\norm{\bm{f}}_{\bm{L}^2}.  
	\end{align*}
	Recalling that $\tilde{\mathcal{C}}^{-1}(\vphi_i) : (\bm{L}^2(\Omega), \bm{H}^{\sfrac{1}{2}}(\Gamma_N) ) \rightarrow \bm{H}^2(\Omega)$ is uniformly bounded for all $\vphi_i$ satisfying the estimate $\norm{\vphi_i}_{ W^{1, q}} \leq R$, see Theorem~\ref{thm:higher_elliptic}, we deduce 
	that
	\begin{align*}
		\norm{\bm{u}}_{L^s(\bm{H}^2)}
		& \leq C(R)  \big( \norm{ - \nabla \cdot (   \C (\vphi) \Tau(\vphi) + \alpha(\vphi) M(\vphi) \theta \bm{I} ) 
			+ \bm{f} }_{L^s(\bm{L}^2)} 
			+ \norm{\bm{g}}_{L^s(\bm{H}^{\sfrac{1}{2}}(\Gamma_N) )} \big) \\ 
		& \leq C(R) \big(\norm{\vphi}_{L^\infty(W^{1, q})} + \norm{\vphi}_{L^\infty(W^{1, q})}^2  \norm{ \theta}_{L^s(H^1)} + \norm{\bm{f} }_{L^s(\bm{L}^2)} + \norm{\bm{g}}_{L^s(\bm{H}^{\sfrac{1}{2}}(\Gamma_N) )}  \big)\\ 
		&\leq C(R) . 
		\numberthis \label{eq:est_u_non_ve}
	\end{align*}
	As for the second part, we derive 
	\begin{align*}
		&\norm{\bm{u}_1 - \bm{u}_2}_{\bm{H}^2}\\ 
		& \quad = \begin{aligned}[t]
			\big\|  \tilde{\mathcal{C}}^{-1}(\vphi_1) &\big( - \nabla \cdot (   \C (\vphi_1) \Tau(\vphi_1) + \alpha(\vphi_1) M(\vphi_1) \theta_1  \bm{I} ) + \bm{f}, \bm{g} \big) \\ 
			&\quad - \tilde{\mathcal{C}}^{-1}(\vphi_2) \big( - \nabla \cdot (   \C (\vphi_2) \Tau(\vphi_2) + \alpha(\vphi_2) M(\vphi_2) \theta_2 \bm{I} ) + \bm{f}, \bm{g} \big) \big\|_{\bm{H}^2}
		\end{aligned} \\ 
		& \quad  \leq \begin{aligned}[t]
			&\big\|\tilde{\mathcal{C}}^{-1}(\vphi_1) \big\|_{\mathcal{L}( (\bm{L}^2, {\bm{H}^{\sfrac{1}{2}} (\Gamma_N)}), \bm{H^2})}  \\ 
			& \qquad \Big( \big\| \C (\vphi_1) \Tau(\vphi_1) - \C (\vphi_2) \Tau(\vphi_2)  \big\|_{\bm{H}^1} 
			+ \big\|\alpha(\vphi_1) M(\vphi_1) \theta_1 -\alpha(\vphi_2) M(\vphi_2) \theta_2   \big\|_{H^1}
			\Big) \\ 
			& + \big\| \tilde{\mathcal{C}}^{-1}(\vphi_1) - \tilde{\mathcal{C}}^{-1}(\vphi_2) \big\|_{ (\bm{L}^2, {\bm{H}^{\sfrac{1}{2}} (\Gamma_N)}), \bm{H^2})}\\
			& \quad \quad \begin{aligned}[t]
				  \Big( \big\|   \C (\vphi_2) \Tau(\vphi_2) \big\|_{\bm{H}^1} &+ \big\| \alpha(\vphi_2) M(\vphi_2) \theta_2 \big\|_{H^1}  
				  +  \norm{\bm{f} }_{\bm{L}^2} + \norm{ \bm{g}}_{\bm{H}^{\sfrac{1}{2}}(\Gamma_N)}  \Big).
			\end{aligned}
		\end{aligned}  \numberthis \label{eq:est_u_non_ve_2} 
	\end{align*}
	In view of this, we need an estimate for the difference of the two operators $\tilde{\mathcal{C}}^{-1}(\vphi_1) - \tilde{\mathcal{C}}^{-1}(\vphi_2)$ to proceed. This is subject of the following claim. 
	\par 
	\medskip 
	\textit{Claim:} There exists a constant $C > 0$ such that 
	\begin{equation}\label{eq:C_^-1_lipschitz}
		\|  \tilde{\mathcal{C}}^{-1}(\phi_1) - \tilde{\mathcal{C}}^{-1}(\phi_2)\|_{\mathcal{L}( (\bm{L}^2, {\bm{H}^{\sfrac{1}{2}} (\Gamma_N)}), \bm{H}^2 )} 
		\leq C \norm{ \phi_1 - \phi_2}_{W^{1, q}}
	\end{equation}
	\indent
	for all 
	$\phi \in W^{1, q}(\Omega)$ with $\norm{\phi_i}_{ W^{1, q}} \leq  R$. \par
	\medskip 
	\textit{Proof of claim:} For any $\bm{w} \in \bm{H}^2(\Omega)$, it holds that 
	\begin{align*}
		\big\| &\nabla \cdot \big( \big[\tilde{\C}(\phi_1) - \tilde{\C}(\phi_2)\big] \E(\bm{w}) \big) \big\|_{\bm{L}^2}
		+ \big\|  \big[\tilde{\C}(\phi_1) -  \tilde{\C}(\phi_2)\big] \E(\bm{w}) \bm{n} \big\|_{\bm{H}^{\sfrac{1}{2}} (\Gamma_N)}\\ 
		&\leq C  \big\|  \big[\tilde{\C}(\phi_1) -  \tilde{\C}(\phi_2)\big] \E(\bm{w}) \big\|_{\bm{H}^1} 
		\leq C  \big\|  \tilde{\C}(\phi_1) -  \tilde{\C}(\phi_2) \big\|_{\bm{W}^{1, q}} \big\| \bm{w}\big\|_{\bm{H}^2}\\ 
		&\leq C \norm{\phi_1 - \phi_2}_{W^{1, q}}  \big\| \bm{w}\big\|_{\bm{H}^2} , 
	\end{align*}
	where the last inequality is a consequence of Theorem \ref{thm:compositon}. In particular, we have
	\begin{equation*}
		\|  \tilde{\mathcal{C}}(\phi_1) - \tilde{\mathcal{C}}(\phi_2)\|_{\mathcal{L}(\bm{H}^2, (\bm{L}^2, {\bm{H}^{\sfrac{1}{2}} (\Gamma_N)})} \leq C \norm{\phi_1 - \phi_2}_{W^{1, q}}
	\end{equation*}
	Invoking Theorem \ref{thm:higher_elliptic}, if follows that  
	\begin{align*}
		& \big\|  \big[ \tilde{\mathcal{C}}^{-1}(\phi_1) - \tilde{\mathcal{C}}^\new{{-1}}(\phi_2) \big] (\tilde{\bm{f}}, \tilde{\bm{g}} )  \big\|_{\bm{H}^2}
		=  \bnorm{  \big[ \tilde{\mathcal{C}}^{-1}(\phi_1) 
			\big[ \tilde{\mathcal{C}}(\phi_2)  -  \tilde{\mathcal{C}}(\phi_1) \big]   \tilde{\mathcal{C}}^{-1}(\phi_2)  \big](\tilde{\bm{f}}, \tilde{\bm{g}} )   }_{\bm{H}^2}\\
		& \quad \leq  \big\| \mathcal{C}^{-1}(\phi_1) \big\|_{\mathcal{L} ( (\bm{L}^2, {\bm{H}^{\sfrac{1}{2}} (\Gamma_N)}), \bm{H}^2)} 
		\big\| 
			\big[ \tilde{\mathcal{C}}(\phi_2)  -  \tilde{\mathcal{C}}(\phi_1) \big]   \tilde{\mathcal{C}}^{-1}(\phi_2)  (\tilde{\bm{f}}, \tilde{\bm{g}} )  
		\big\|_ {(\bm{L}^2, {\bm{H}^{\sfrac{1}{2}} (\Gamma_N)})}
		\\ 
		& \quad \leq   C(R)
		\big\|  \big[ \tilde{\mathcal{C}}(\phi_2)  -  \tilde{\mathcal{C}}(\phi_1)\big\| _{\mathcal{L} (\bm{H}^2,(\bm{L}^2, {\bm{H}^{\sfrac{1}{2}} (\Gamma_N)})}
		\big\| \mathcal{C}^{-1}(\phi_2) \big\|_{\mathcal{L} ((\bm{L}^2, {\bm{H}^{\sfrac{1}{2}} (\Gamma_N)}), \bm{H}^2)}
		\big\| (\tilde{\bm{f}}, \tilde{\bm{g}} ) \big\|_{(\bm{L}^2, {\bm{H}^{\sfrac{1}{2}} (\Gamma_N)})}\\ 
		& \quad \leq C(R) \norm{\phi_1 - \phi_2}_{W^{1, q}} \big\| (\tilde{\bm{f}}, \tilde{\bm{g}} ) \big\|_{(\bm{L}^2, {\bm{H}^{\sfrac{1}{2}} (\Gamma_N)})}, 
	\end{align*}
	for all $\tilde{\bm{f}} \in \bm{L}^2(\Omega), \tilde{\bm{g}} \in \bm{H}^{\sfrac{1}{2}} (\Gamma_N)$. This shows the claim. \hfill $\diamondsuit$\par 
	\medskip 
	Upon inserting this estimate into \eqref{eq:est_u_non_ve_2} and  recalling that the $\norm{\vphi_i}_{W^{1, q}} \leq R$, we find
	\begin{align*}
		\norm{\bm{u}_1 - \bm{u}_2}_{\bm{H}^2} 
		\leq \begin{aligned}[t]
			&C(R) \big( \norm{\vphi_1 - \vphi_2}_{W^{1, q}} + \norm{\theta_1 - \theta_2}_{H^1} + \norm{\theta_2} \norm{\vphi_1 - \vphi_2} _{W^{1, q}} \big) \\ 
			&  +C(R) \norm{ \vphi_1 - \vphi_2}_{W^{1, q}}  \big(1  + \norm{\theta_2}_{H^1} + \norm{\bm{f}}_{\bm{L}^2}  + \norm{ \bm{g}}_{\bm{H}^{\sfrac{1}{2}}(\Gamma_N)} \big), 
		\end{aligned}
	\end{align*}
	and in conjunction with the assumptions \ref{A:f}, \ref{A:g}, we obtain the assertion. 
\end{proof}

In light of the claim we showed as part of Lemma \ref{lem:contraction_u_non_ve}, it is now easy to deduce that the operator $\A$ is also Lipschitz continuous. 

\begin{corollary}\label{cor:TA_lipschitz}
	There exists a constant $C = C(R) > 0$ such that 
	\begin{equation}
		\|  \TA(\phi_1) - \TA(\phi_2)\|_{\mathcal{L}(H^1)} 
		\leq C(R) \norm{ \phi_1 - \phi_2}_{W^{1, q}}
	\end{equation}
	for all 
	$\phi_i \in W^{1, q}(\Omega)$, $q > n$, with $\norm{\phi_i}_{W^{1, q}} \leq  R$. 
\end{corollary}
\begin{proof}
	Let $\vartheta \in H^1(\Omega)$ and define $\hat{\alpha}(\phi_i) \coloneqq \alpha(\phi_i) M(\phi_i)$. We compute with the help of Theorem \ref{thm:compositon} and the estimate \eqref{eq:C_^-1_lipschitz}
	\begin{align*}
		&\norm{ \big(\TA(\phi_1) -  \TA(\phi_2) \big) \vartheta}_{H^1}\\ 
		&\quad = \begin{aligned}[t]
			\big\| M(\phi_1)\big( \vartheta - \alpha(\phi_1) \nabla &\cdot \big[ \tilde{\mathcal{C}}^{-1}(\phi_1) \big( - \nabla  (   \alpha(\phi_1) M(\phi_1) \vartheta ), \bm{0} \big) \big]  \big) \\ 
			& -M(\phi_2)\big( \vartheta - \alpha(\phi_2) \nabla \cdot \big[ \tilde{\mathcal{C}}^{-1}(\phi_2) \big( - \nabla  (   \alpha(\phi_2) M(\phi_2) \vartheta ), \bm{0}  \big) \big]  \big) \big\|_{H^1}
		\end{aligned}\\ 
		&\quad \leq \begin{aligned}[t]
			\norm{ \big(M(\phi_1) - M(\phi_2) \big) \vartheta}_{H^1}
			&+ \bnorm{\hat{\alpha}(\phi_1)\nabla \cdot \big[ \big( \tilde{\mathcal{C}}^{-1}(\phi_1) -  \tilde{\mathcal{C}}^{-1}(\phi_2)\big) \big( - \nabla  (   \hat{\alpha}(\phi_1)  \vartheta ), \bm{0}  \big) \big]  }_{H^1}\\ 
			& + \bnorm{ \big( \hat{\alpha}(\phi_1) - \hat{\alpha}(\phi_2)\big) \nabla \cdot \big[\tilde{\mathcal{C}}^{-1}(\phi_2) \big( - \nabla  (   \hat{\alpha}(\phi_1)  \vartheta ), \bm{0}  \big) \big]  }_{H^1}\\ 
			& + \bnorm{  \hat{\alpha}(\phi_2) \nabla \cdot \big[\tilde{\mathcal{C}}^{-1}(\phi_2) \big( - \nabla  \big(     ( \hat{\alpha}(\phi_1) - \hat{\alpha}(\phi_2))  \vartheta \big), \bm{0}  \big) \big]  }_{H^1}
		\end{aligned}\\ 
		&\quad \leq \begin{aligned}[t]
			C \norm{\phi_1 -\phi_2}_{W^{1, q}} \norm{ \vartheta}_{H^1}
			&+ C(R) \bnorm{  \tilde{\mathcal{C}}^{-1}(\phi_1) -  \tilde{\mathcal{C}}^{-1}(\phi_2) }_{\mathcal{L}((\bm{L}^2, {\bm{H}^{\sfrac{1}{2}} (\Gamma_N)}), \bm{H}^2)} \norm{ \vartheta }_{H^1}\\ 
			& + C(R)\norm{\phi_1 - \phi_2}_{W^{1, q}} \bnorm{ \tilde{\mathcal{C}}^{-1}(\phi_2) }_{\mathcal{L}((\bm{L}^2, {\bm{H}^{\sfrac{1}{2}} (\Gamma_N)}), \bm{H}^2)} \norm{\vartheta}_{H^1}\\ 
			&+ C(R)\bnorm{ \tilde{\mathcal{C}}^{-1}(\phi_2) }_{\mathcal{L}((\bm{L}^2, {\bm{H}^{\sfrac{1}{2}} (\Gamma_N)}), \bm{H}^2)} \norm{ \phi_1 - \phi_2}_{W^{1, q}} \norm{\vartheta}_{H^1} 
		\end{aligned}\\ 
		&\quad  \leq C  \norm{\phi_1 -\phi_2}_{W^{1, q}} \norm{ \vartheta}_{H^1} . 
	\end{align*}
	Here, we also applied Theorem \ref{thm:higher_elliptic}. 
\end{proof}

With the help of Lemma \ref{lem:contraction_u_non_ve} we can further find the following estimates for the pressure $p$. 

\begin{corollary} \label{cor:p_lipschitz}
	If $\vphi\in  L^\infty( [0, T]; W^{1, q} (\Omega))$, $q >n$ and $\theta\in L^s(0, T; H^1(\Omega))$, then the function 
	\begin{equation*}
		p = M(\vphi) (\theta - \alpha(\vphi) \nabla \cdot \bm{u})
	\end{equation*}
	belongs to $L^s(0, T; H^1(\Omega))$ and is bounded by a constant $C = C(R)$. Moreover, if $\norm{\vphi_i}_{L^\infty(W^{1, q})} \leq R$, $i = 1, 2$, then it holds for almost every $t \in (0, T)$ that  
	\begin{equation*}
		\norm{p_1 - p_2}_{H^1} \leq C(R)  \norm{\vphi_1 - \vphi_2}_{W^{1, q}}  \big( 1 + \norm{\theta_2}_{H^1} + \norm{\bm{u}_2}_{\bm{H}^2}  + \norm{\bm{f}}_{\bm{L}^2} + \norm{ \bm{g}}_{\bm{H}^{\sfrac{1}{2}}(\Gamma_N)} \big) + C \norm{\theta_1 - \theta_2}_{H^1}. 
	\end{equation*}
\end{corollary}
\begin{proof}
	The first assertion is an immediate consequence of Lemma~\ref{lem:contraction_u_non_ve} and the definition of $p$.  
	As for the second part, we define $\hat{\alpha}_i = \alpha(\vphi_i) M(\vphi_i) $ to compute
	\begin{align*}
		\norm{p_1 - p_2}_{H^1} 
		&\leq \norm{ M(\vphi_1) \theta_1 -  M(\vphi_2) \theta_2}_{H^1} + \norm{ \hat{\alpha}(\vphi_1) \nabla \cdot \bm{u}_1  - \hat{\alpha}(\vphi_2) \nabla \cdot \bm{u}_2}_{H^1}\\ 
		&\leq  \begin{aligned}[t]
			C(R) \norm{\theta_1 - \theta_2}_{H^1} &+  C\norm{\vphi_1 - \vphi_2}_{W^{1, q}} \norm{\theta_2}_{H^1}\\ 
			&+  C(R) \norm{ \bm{u}_1 - \bm{u}_2}_{\bm{H}^2} 
			+  C \norm{\vphi_1 - \vphi_2}_{W^{1, q}} \norm{\bm{u}_2}_{\bm{H}^2}, 
		\end{aligned}
	\end{align*}
	and obtain the assertion upon inserting the estimate from Lemma \ref{lem:contraction_u_non_ve}. 
\end{proof}

With these preparation in order, we can finally begin the  proof of Proposition \ref{prop:contraction_non_ve}. 

\begin{proof}[Proof of Proposition \ref{prop:contraction_non_ve}]
 To establish Lipschitz continuity for $\F_T$, we need to study the components $\F_T^j$, $j = 1, 2$, separately. 
	\par 
	\medskip 
	\begin{subequations}
	\textit{Ad $\F_T^1$:} For the first component\new{, we use the definition of $\mathcal{F}_T^1$ and find}
	\begin{align*}
		&\norm{\F_T^1(\vphi_1,, \theta_1) - \F_T^1 (\vphi_2, \theta_2)}_{Y^1_T} \\ 
		&\quad \leq   \Big\| \eps\, \Delta  \big( m(\vphi_0)  \Delta (\vphi_1 - \vphi_2 ) \big)
		 - \eps \nabla \cdot \big( m(\vphi_1) \nabla \Delta \vphi_1  - m(\vphi_2) \nabla \Delta \vphi_2  \big) \Big\| _{Y^1_T}  \\ 
		& \quad \quad + \Big\| \nabla \cdot \big(
			m(\vphi_1) \nabla \big(\frac{1}{\eps} \psi'(\vphi_1)  \big)  
			- m(\vphi_2) \nabla \big( \frac{1}{\eps} \psi'(\vphi_2)\big) \big)\Big\|_{Y^1_T}\\ 
		& \quad \quad + \Big\| \nabla \cdot \big(
		m(\vphi_1) \nabla \Wp(\vphi_1, \E(\bm{u}_1))
		- m(\vphi_2) \nabla \Wp(\vphi_2, \E(\bm{u}_2)) \big)\Big\|_{Y^1_T}\\ 
		& \quad \quad + \begin{aligned}[t]
			 \Big\| \nabla \cdot \big( m(\vphi_1) \nabla \big( p_1\, \alpha'(\vphi_1) \nabla \cdot \bm{u}_1  \big) 
			 - m(\vphi_2) \nabla \big( p_2\, \alpha'(\vphi_2) \nabla \cdot \bm{u}_2 \big) \big) \Big\|_{Y^1_T}
		\end{aligned}
		\\ 
		& \quad \quad + \Big\| \nabla \cdot \big(
		m(\vphi_1) \nabla \big(\frac{M'(\vphi_1)}{2 M(\vphi_1)} p_1^2 \big)  
		- m(\vphi_2) \nabla \big(\frac{M'(\vphi_2)}{2 M(\vphi_2)} p_2^2
		\big) \big) \Big\|_{Y^1_T} \\ 
		& \quad \quad + \norm{S_s(\vphi_1, \bm{u}_1,  \theta_1) - S_s(\vphi_2, \bm{u}_2,  \theta_2)}_{Y^1_T}\\ 
		& \quad \eqqcolon I + II + III+ IV + V +VI . 
	\end{align*}
	\textit{Ad $\F_T^1.I$:} Let us observe that $W^{3, \frac{3}{2}}(\Omega) \hookrightarrow W^{2, \frac{3n}{2n-3}}$ for $n \in \{2,3\}$, and set
	\begin{equation*}
		\frac{1}{q^*} = \frac{1}{\sfrac{3}{2}} - \frac{1}{q} > 0. 
	\end{equation*}
	Due to $q > n$, we find that $q^* < \frac{3n}{2n-3}$, giving rise to the interpolation estimate 
	\begin{equation*}
		\norm{f}_{W^{2, q^*}} \leq M \norm{f}_{(W^{1, 3})'}^{1-\hat{\vartheta}}  \norm{f}_{W^{3, \frac{3}{2}}}^{\hat{\vartheta}}
	\end{equation*}
	for some $M > 0$ and $\hat{\vartheta} \in (0, 1)$. Recalling that $\vphi_i \in L^r(0, T; W^{3, \frac{3}{2}}) \cap L^\infty(0, T; (W^{1,3})' )$, \new{an application of the generalized Hölder inequality yields}
	\begin{equation*}
		\norm{\vphi_i }_{L^{r^*}(W^{2, q^*})} \leq M \norm{\vphi_i}_{L^\infty((W^{1, 3})')}^{1-\hat{\vartheta}}  \norm{\vphi_i}_{L^r(W^{3, \frac{3}{2}})}^{\hat{\vartheta}}, 
	\end{equation*}
	where $\frac{1}{r^*} = \frac{\hat{\vartheta}}{r}$. In particular, $\vphi_1 - \vphi_2 \in {L^{r^*}(0, T; W^{2, q^*})}$ for some $r^* > r$ and the norm is bounded by the respective norms in $X_T^1$. 
	Applying Theorem \ref{thm:compositon} along with similar arguments as in Lemma \ref{lemma:general_estimate_non_ve}, we infer
	\begin{align*}
		 &\big\|   \Delta  \big( m(\vphi_0)  \Delta (\vphi_1 - \vphi_2 ) \big)
		 -  \nabla \cdot \big( m(\vphi_1) \nabla \Delta \vphi_1  
		 - m(\vphi_2) \nabla \Delta \vphi_2  \big) \big\| _{(W^{1, 3})'}\\ 
		 & \quad \leq \begin{aligned}[t]
		 	&\big\|   \nabla \cdot  \big( m(\vphi_0)  \nabla\Delta (\vphi_1 - \vphi_2 ) \big)
		 	-  \nabla \cdot \big( m(\vphi_1) \nabla \Delta \vphi_1  
		 	- m(\vphi_2) \nabla \Delta \vphi_2  \big) \big\| _{(W^{1, 3})'}\\ 
		 	&+\big\|  \nabla \cdot \big(  m'(\vphi_0) \nabla \vphi_0\,  \Delta (\vphi_1 - \vphi_2  \big)
		 	\big\|_{(W^{1, 3})'}
		 \end{aligned}
		 \\ 		 
		 & \quad \leq\begin{aligned}[t]
		 	 &\big\|  \nabla \cdot \Big( \big( m(\vphi_0) - m(\vphi_1) \big) \big( \nabla \Delta \vphi_1 -  \new{\nabla}\Delta \vphi_2 \big) \Big)\big\|_{(W^{1, 3})'}
		 	+ \big\| \nabla \cdot  \Big( \big(m(\vphi_1) - m(\vphi_2) \big) \nabla \Delta \vphi_2  \Big) \big\|_{(W^{1, 3})'}\\ 
		 	&+ \norm{m'(\vphi_0) \nabla \vphi_0 \, \Delta (\vphi_1 - \vphi_2) }_{L^{\frac{3}{2}}}
		 \end{aligned}\\ 
		 & \quad \leq \begin{aligned}[t]
		 	 &C(R) \norm{\vphi_0 - \vphi_1}_{W^{1, q}} \norm{\nabla \Delta \vphi_1 - \nabla \Delta \vphi_2}_{L^{\sfrac{3}{2}}}
		 	+C(R)  \norm{\vphi_1 - \vphi_2}_{W^{1, q}} \norm{\nabla \Delta \vphi_2 }_{L^{\sfrac{3}{2}}}\\ 
		 	&+C(R) \norm{m'(\vphi_0)}_{W^{1, q}} \norm{ \vphi_0}_{W^{1, q}} \norm{\vphi_1 - \vphi_2}_{W^{2, q^*}}
		 \end{aligned}\\ 
		 & \quad  \leq C(R) \norm{\vphi_0 - \vphi_1}_{W^{1, q}} \norm{\vphi_1 -  \vphi_2}_{W^{3, \sfrac{3}{2}}}
		 + C(R) \norm{\vphi_1 - \vphi_2}_{W^{1, q}} \norm{ \vphi_2 }_{W^{3, \sfrac{3}{2}}}
		 +C(R)  \norm{\vphi_1 - \vphi_2}_{W^{2, q^*}},
	\end{align*}
	and \eqref{embedding:hoelder_non_ve} leads to  
	\begin{align*}
		 I &\leq C(R) \begin{aligned}[t]
		 	\big( \norm{ \vphi_1(0) - \vphi_1}_{C^0(W^{1, q})} \norm{\vphi_1 -  \vphi_2}_{L^r(W^{3, \sfrac{3}{2}})}
		 	&+ \norm{\vphi_1 - \vphi_2}_{C^0(W^{1, q})}  \norm{ \vphi_2}_{L^r(W^{3, \sfrac{3}{2}})}\\ 
		 	&+\norm{\vphi_1 - \vphi_2}_{L^r(W^{2, q^*})} \big)
		 \end{aligned}\\
		&  \leq C (R) \begin{aligned}[t]
			T^{\tilde{\vartheta}}   \big(   \norm{ \vphi_1}_{C^{\tilde{\vartheta}}(W^{1, q})} \norm{\vphi_1 -  \vphi_2}_{L^r(W^{3, \sfrac{3}{2}})}
			&+  \norm{\vphi_1 - \vphi_2}_{C^{\tilde{\vartheta}}(W^{1, q})}    \big)\\
			&+   C (R) \norm{\vphi_1 - \vphi_2}_{L^{r^*}(W^{2, q^*})}. 
		\end{aligned}
		 \numberthis \label{est:F_1_I_non_ve}
	\end{align*} 
	\\
	\textit{Ad $\F_T^1.II$:} We recall that the functions $\vphi_i$, $i = 1, 2$, are bounded in $X^1_T$ and that \eqref{embedding_X^1_non_ve} provides the embedding $X^1_T \hookrightarrow BUC ([0, T]; W^{1, q}(\Omega))$. Hence, $\norm{\vphi_i}_{W^{1, q}} \leq  C(R)$ and utilizing that $\psi \in C^3(\R)$, see \ref{II.A:psi}, we infer by Theorem \ref{thm:compositon}
	\begin{equation*}
		\norm{	\psi'(\vphi_i)}_{W^{1, q}} < C(R)  \quad \textnormal{and} \quad \norm{\psi'(\vphi_1) - \psi'(\vphi_2)}_{W^{1, q}} \leq  C(R) \norm{\vphi_1 - \vphi_2}_{W^{1, q}} . 
	\end{equation*}
	Lemma~\ref{lemma:general_estimate_non_ve}, along with the embedding $W^{1, q}(\Omega) \hookrightarrow W^{1, \sfrac{3}{2}}$ therefore yields
	\begin{equation*}
		 \Big\| \nabla \cdot \big(
		m(\vphi_1) \nabla \big(\frac{1}{\eps} \psi'(\vphi_1)  \big)  
		- m(\vphi_2) \nabla \big( \frac{1}{\eps} \psi'(\vphi_2)\big) \big)\Big\|_{(W^{1, 3})'}
		\leq C (R) \norm{\vphi_1 - \vphi_2}_{W^{1, q}} , 
	\end{equation*}
	and we conclude 
	\begin{equation}\label{est:F_1_III_non_ve}
		II \leq \norm{ C (R) \norm{\vphi_1 - \vphi_2}_{W^{1, q}}}_{L^r} \leq C(R) T^{\frac{1}{r}} \norm{\vphi_1 - \vphi_2}_{C^0(W^{1, q})}. 
	\end{equation}
	\\ 
	\textit{Ad $\F_T^1.III$:} Recall that in \eqref{eq:est_u_non_ve} we observed $\bm{u}_i \in L^s(0, T; \bm{H}^2(\Omega))$ along with an estimate in the corresponding norm. Unfortunately, that result does not yield an estimate in $L^\infty(0, T; \bm{H}^2(\Omega))$, thus precluding the application of Theorem \ref{thm:compositon}, and requiring us to examine the  growth conditions of the derivative $\nabla \Wp$. By chain and product rule, we infer for any $\vphi \in W^{1, q}(\Omega)$, $\bm{u} \in \bm{H}^2(\Omega)$, 
	\begin{align*}
		\norm{\Wp(\vphi, \E(\bm{u}))}_{W^{1, \sfrac{3}{2}}} 
		\leq \begin{aligned}[t]
			&\norm{\Wp(\vphi, \E(\bm{u}))}_{L^{\sfrac{3}{2}}} 
			+  \big\|  W_{, \vphi \vphi }(\vphi, \E(\bm{u})) \nabla \vphi \big\|_{\bm{L}^{\sfrac{3}{2}}}
			+   \big\|  W_{,\vphi\E}(\vphi, \E(\bm{u}))  \nabla \E(\bm{u})  \big\|_{\bm{L}^{\sfrac{3}{2}}}
		\end{aligned}
	\end{align*}
	where the higher-dimensional multiplication operations should be interpreted appropriately. 
	Recalling the growth conditions \ref{II.A:W_growth_conditions}, it follows that 
	\begin{align} \label{eq:estimate_Wp_nve}
		\norm{\Wp(\vphi, \E(\bm{u}))}_{W^{1, \sfrac{3}{2}}} 
		& \leq C \Big( \begin{aligned}[t]
			\norm{\abs{\vphi}^2 + \abs{\E(\bm{u})}^2 + 1}_{L^{\sfrac{3}{2}}}
			&+ \norm{ (\abs{\vphi}^2 + \abs{\E(\bm{u})}^2 + 1) \abs{\nabla \vphi} }_{L^{\sfrac{3}{2}}}\\ 
			&+ \norm{ (\abs{\vphi} + \abs{\E(\bm{u})} + 1) \abs{ \abs{\nabla \E(\bm{u}) } }}_{L^{\sfrac{3}{2}}}
			\Big) . 
		\end{aligned}
	\end{align}
	If $n = 3$, we obtain 
	\begin{equation}\label{eq:W_p_est_1}
		\norm{  \abs{\E(\bm{u})}^2  \abs{\nabla \vphi}  }_{L^{\sfrac{3}{2}}} \leq \norm{ \E(\bm{u})^2 }_{L^3} \norm{\nabla \vphi}_{L^3} \leq C\norm{\E(\bm{u})}_{\bm{L}^6}^2 \norm{\vphi}_{W^{1, q}} \leq C(R)\norm{\bm
			u}_{\bm{H}^2}^2,
	\end{equation} and if $n = 2$ we use 
	\begin{equation}\label{eq:W_p_est_2}
		\norm{  \abs{\E(\bm{u})}^2  \abs{\nabla \vphi}  }_{L^{\sfrac{3}{2}}} \leq \norm{ \E(\bm{u})^2 }_{L^6} \norm{\nabla \vphi}_{L^2} \leq 
		C(R) \norm{\E(\bm{u})}_{\bm{L}^{12}}^2
		\leq C(R)\norm{\bm{u}}_{\bm{H}^2}^2.
	\end{equation}
	Independently of the dimension, it holds that 
	\begin{equation}\label{eq:W_p_est_3}
		\norm{  \abs{\E(\bm{u})}  \abs{\nabla \E(\bm{u})}  }_{\bm{L}^{\sfrac{3}{2}}} \leq \norm{ \E(\bm{u}) }_{L^6} \norm{\nabla \E(\bm{u})}_{\bm{L}^2} \leq C\norm{\bm{u}}_{\bm{H}^2}^2. 
	\end{equation} 
	Observe that it is not possible to estimate this term in $\bm{L}^2(\Omega)$, which is the reason that we can \textit{not} prove the existence of solutions with $\vphi \in W^{3, 2}(\Omega)$, but only in $W^{3, \sfrac{3}{2}}(\Omega)$. 
	The embedding $W^{1, q}(\Omega) \hookrightarrow L^\infty(\Omega)$, combined with Hölder's inequality and standard arguments leads to 
	\begin{equation*}
			\norm{\Wp(\vphi, \E(\bm{u}))}_{W^{1, \sfrac{3}{2}}} 
			\leq C(R) \big(
				1  + \norm{\bm{u}}_{\bm{H}^2}^2
			\big). 
	\end{equation*}
	Tedious computations relying on the growth conditions \ref{A:W_growth_conditions} combined with the mean value theorem in the multivariable setting further show 
		\begin{align*}
			&\norm{ \Wp(\vphi_1, \E(\bm{u}_1))  - \Wp(\vphi_2, \E(\bm{u}_2)) }_{W^{1, \sfrac{3}{2}}} \\ 
			& \quad \leq C(R)\Big(  \big(1 + \norm{ \bm{u}_1}_{\bm{H}^2}^2 +  \norm{\bm{u}_2}_{\bm{H}^2}^2 \big)  \norm{\vphi_1 - \vphi_2}_{W^{1, q}}		
			+\big( 1 + \norm{\bm{u}_1}_{\bm{H}^2} +  \norm{\bm{u}_2}_{\bm{H}^2} \big) \norm{\bm{u}_1 - \bm{u}_2}_{\bm{H}^2} \Big). \numberthis \label{eq:wp_4_nve}
		\end{align*}
	Recalling that $s = 4r$, Lemma \ref{lemma:general_estimate_non_ve} and Lemma \ref{lem:contraction_u_non_ve} therefore imply 
	\begin{align*}
		III &= \Big\| \nabla \cdot \big(
		m(\vphi_1) \nabla \Wp(\vphi_1, \E(\bm{u}_1))
		- m(\vphi_2) \nabla \Wp(\vphi_2, \E(\bm{u}_2)) \big)\Big\|_{L^r((W^{1, 3})' )} \\ 
		& \leq C(R) 
		\begin{aligned}[t]
			 \big( &\norm{ \Wp(\vphi_1, \E(\bm{u}_1)) -  \Wp(\vphi_2, \E(\bm{u}_2))}_{L^r(W^{1, \sfrac{3}{2}})} \\
			 &+ \norm{\Wp(\vphi_2, \E(\bm{u}_2))}_{L^r(W^{1, \sfrac{3}{2}})} \norm{\vphi_1 - \vphi_2}_{C^0(W^{1, q})} \big)
		\end{aligned}\\ 
		 &\leq  C(R) T^{\frac{1}{2r}} \big(  \norm{\vphi_1 - \vphi_2}_{C^0(W^{1, q})} 
		+  \norm{\bm{u}_1 - \bm{u}_2}_{L^s(\bm{H}^2)}
		 \big)\\ 
		 &\leq  C(R, T)  \big(  \norm{\vphi_1 - \vphi_2}_{C^0(W^{1, q})} 
		 + \norm{\theta_1 - \theta_2}_{L^s(H^1)}
		 \big)
		 \numberthis \label{F_1_IV_non_ve}
	\end{align*}
	with $C(R, T) \rightarrow 0$ as $T \rightarrow 0$. 
	\par 
	\medskip
	\noindent
	\textit{Ad $\F_T^1.IV$:} An application of Lemma \ref{lemma:general_estimate_non_ve} to $(IV)$ along with the multiplication results $W^{1, q} \times H^1(\Omega) \rightarrow  H^1(\Omega)$ and $H^1(\Omega) \times H^1(\Omega) \rightarrow  W^{1, \sfrac{3}{2}}(\Omega)$ gives us  
	\begin{align*}
		&\begin{aligned}
			 \Big\| \nabla \cdot \big( m(\vphi_1) \nabla \big( p_1\, \alpha'(\vphi_1) \nabla \cdot \bm{u}_1  \big) 
			- m(\vphi_2) \nabla \big( p_2\, \alpha'(\vphi_2) \nabla \cdot \bm{u}_2 \big) \big) \Big\|_{(W^{1, 3})'}
		\end{aligned}
		\\ & \quad \leq 
		 \begin{aligned}[t]
			C(R) \Big( 
			 &\norm{\alpha'(\vphi_1) p_1 ( \nabla \cdot \bm{u}_1 - \nabla \cdot \bm{u}_2)}_{W^{1, \sfrac{3}{2}}}
			 + \norm{ \alpha'(\vphi_1) (\nabla \cdot \bm{u}_2) (p_1 - p_2) }_{W^{1, \sfrac{3}{2}}}\\ 
			&+ \norm{ p_2 (\nabla \cdot \bm{u}_2) ( \alpha'(\vphi_1) - \alpha'(\vphi_2)) }_{W^{1, \sfrac{3}{2}}}
			+  \norm {p_2\, \alpha'(\vphi_2) \nabla \cdot \bm{u}_2}_{W^{1, \sfrac{3}{2}}} \norm{\vphi_1 - \vphi_2}_{W^{1,q}} 
				\Big). 
		\end{aligned}\\ 
		& \quad \leq  \begin{aligned}[t]
			C(R) \Big( 
			& \norm{\alpha'(\vphi_1)}_{W^{1, q}} \norm{p_1}_{H^1}\norm{\bm{u}_1 -  \bm{u}_2}_{\bm{H}^2}
			+ \norm{\alpha'(\vphi_1)}_{W^{1, q}} \norm{\bm{u}_2}_{\bm{H}^2} \norm{ p_1 - p_2}_{H^1}\\ 
			&+   \norm{p_2}_{H^1}  \norm{\bm{u}_2}_{\bm{H}^2}  \norm{\vphi_1 - \vphi_2}_{W^{1, q}} 
			+  \norm{p_2}_{H^1}  \norm{ \bm{u}_2}_{\bm{H}^2}   \norm{\vphi_1 - \vphi_2}_{W^{1, q}} 
			\Big), 
		\end{aligned}
	\end{align*}
	and we deduce with Corollary \ref{cor:p_lipschitz}, \ref{II.A:g}, \ref{II.A:f}, and Lemma \ref{lem:contraction_u_non_ve} 
	\begin{align*}
		IV \leq 
			C(R) T^{\frac{1}{2r}}   \big(  \norm{\vphi_1 - \vphi_2}_{C^0(W^{1, q})} +   \norm{\theta_1 - \theta_2}_{L^s(H^1)} 
			 \big) . 
			\numberthis \label{est:F_1_V_non_ve}  
	\end{align*}
	\\ 
	\textit{Ad $\F_T^1.V$:} We define $\hat{M}(\vphi_i) \coloneqq \frac{M'(\vphi_i)}{2 M(\vphi_i)}$ and apply Lemma \ref{lemma:general_estimate_non_ve} along with Theorem \ref{thm:compositon} to find 
		\begin{align*}
		&\begin{aligned}
			\Big\| \nabla \cdot \big( m(\vphi_1) \nabla \big( \hat{M}(\vphi_1)  p_1^2 \big) 
			- m(\vphi_2) \nabla \big(\hat{M}(\vphi_2) p_2^2 \big) \big) \Big\|_{L^2}
		\end{aligned}
		\\ & \quad \leq 
		\begin{aligned}[t]
			C(R) \Big( 
			&\norm{ \hat{M}(\vphi_1) p_1 (p_1 - p_2)}_{H^2}
			+ \norm{  \hat{M}(\vphi_1) p_2 (p_1 - p_2 )}_{H^2}\\ 
			&+ \norm{ p_2^2 (  \hat{M}(\vphi_1) -  \hat{M} (\vphi_2)) }_{H^2}
			+  \norm {\hat{M}(\vphi_2) \, p_2^2}_{H^2} \norm{\vphi_1 - \vphi_2}_{H^2} 
			\Big). 
		\end{aligned}\\ 
		& \quad \leq  \begin{aligned}[t]
			C(R) \Big( 
			& \norm{\vphi_1}_{H^2} \norm{p_1}_{H^2} \norm{ p_1 - p_2}_{H^2}
			+ \norm{\vphi_1}_{H^2} \norm{p_2}_{H^2} \norm{ p_1 - p_2}_{H^2}\\ 
			&+   \norm{p_2}_{H^2}^2  \norm{\vphi_1 - \vphi_2}_{H^2} 
			+  \norm{p_2}_{H^2}^2  \norm{\vphi_1 - \vphi_2}_{H^2} 
			\Big), 
		\end{aligned}
		\end{align*}
	and conclude with help from Corollary \ref{cor:p_lipschitz} and \ref{II.A:g}, \ref{II.A:f}, that 
	\begin{align*}
		V&\leq C(R) T^{\frac{1}{2r}}   \big(  \norm{\vphi_1 - \vphi_2}_{C^0(H^2)} +   \norm{\theta_1 - \theta_2}_{L^s(H^2)} 
		\big) . 
		\numberthis \label{est:F_1_VI_non_ve}
	\end{align*}
	\\ 
	\textit{Ad $\F_T^1.VI$:} The Lipschitz assumption \ref{II.A:source_terms} exactly states
	\begin{align*}
		VI &= 
		 \norm{S_s(\vphi_1, \bm{u}_1,  \theta_1) - S_s(\vphi_2, \bm{u}_2,  \theta_2)}_{L^r(L^2)}
		 \leq C(T, R) ( \norm{\vphi_1 - \vphi_2}_{\new{Z}_T^1} + \norm{\theta_1 - \theta_2}_{\new{Z}_T^2}), \numberthis \label{est:F_1_VII_non_ve}
	\end{align*}
	where the constant $C (R, T)$ tends to zero as $T$ goes to zero. 
	\par 
	\medskip 
	Combining the estimates \eqref{est:F_1_I_non_ve}-\eqref{est:F_1_VII_non_ve} with the embeddings \eqref{embedding_X^1_non_ve} and \eqref{embedding:hoelder_non_ve} finally leads to 
	\begin{align*}
		&\norm{\F_T^1 (\vphi_1, \theta_1) - \F_T^1 (\vphi_2, \theta_2)}_{Y^1_T}
		  \leq C(R, T) \Big( \norm{\vphi_1 - \vphi_2}_{\new{Z}^1_T} 
		+ \norm{\theta_1 - \theta_2}_{\new{Z}^2_T}
		  \Big)
		  	\numberthis \label{est:F_1_non_ve}
		\end{align*}
	where $C(R, T) > 0$ is a constant depending on $R$ and $T$ such that $C(R, T) \rightarrow 0$ as $T \rightarrow 0$. 
	\par 
	\end{subequations}
		\par 
		\medskip 
		\noindent
	\textit{Ad $\F_T^2$:} We need to repeat the same process for the second component $\F_2$
	\begin{align*}
		&\norm{\F_T^2 (\vphi_1, \theta_1) - \F_T^2 (\vphi_2,  \theta_2)}_{Y^2_T}\\ 
		& \ \  \leq \norm{ \A(\vphi_0)  \big(\theta_1 - \theta_2 \big) -  \big(\A(\vphi_1)\theta_1 - \A(\vphi_2) \theta_2  \big) }_{Y^2_T}\\ 
		& \quad  + \big\|
				\begin{aligned}[t]
					&\nabla \cdot \Big( \kappa(\vphi_1) \nabla \Big( \alpha(\vphi_1) M(\vphi_1) \nabla \cdot  \big[ \tilde{\mathcal{C}}^{-1}(\vphi_1) \big( -\nabla(\C(\vphi_1) \Tau(\vphi_1)) + \bm{f }, \bm{g}  \big) \big]   \Big)\Big)\\ 
					&- \nabla \cdot \Big( \kappa(\vphi_2) \nabla \Big( \alpha(\vphi_2) M(\vphi_2) \nabla \cdot  \big[ \tilde{\mathcal{C}}^{-1}(\vphi_2) \big( -\nabla(\C(\vphi_2) \Tau(\vphi_2)) + \bm{f }, \bm{g}  \big) \big]   \Big)\Big)\big\| _{Y^2_T}\\  
				\end{aligned}\\ 
		& \quad + \norm{S_f(\vphi_1, \E(\bm{u}_1), \theta_1)  - S_f(\vphi_2, \E(\bm{u}_2), \theta_2)}_{Y_T^2}\\
		& \ \  \eqqcolon I + II + III  . 
	\end{align*}
	\begin{subequations}
	\textit{Ad $\F_T^2.I$:}  Due to the additional dependency of $\TA$, \new{see \eqref{def:mathcal_A}}, on $\vphi_i$, the first term has to be expanded as
	\begin{align*}
		& \norm{ \A(\vphi_0)  \big(\theta_1 - \theta_2 \big) -  \big(\A(\vphi_1)\theta_1 - \A(\vphi_2) \theta_2  \big) }_{(H^1)'}\\ 
		& \quad =  \norm{ \big( \A(\vphi_0) -  \A(\vphi_1) \big)  \big(\theta_1 - \theta_2 \big) -  \big(\A(\vphi_1) - \A(\vphi_2)  \big) \theta_2}_{(H^1)'} \\ 
		& \quad \leq \begin{aligned}[t]
			\big\| - \nabla \cdot \big( \big( \kappa(\vphi_0) \nabla \TA(\vphi_0) &-   \kappa(\vphi_1) \nabla \TA(\vphi_1) \big) (\theta_1 - \theta_2) \big)\big\|_{(H^1)'} \\
			&+  \big\| - \nabla \cdot \big( \big( \kappa(\vphi_1) \nabla \TA(\vphi_1) -   \kappa(\vphi_2) \nabla \TA(\vphi_2) \big) \theta_2 \big)\big\|_{(H^1)'}
		\end{aligned}\\ 
		& \quad \leq \begin{aligned}[t]
			&\Big( \big\|  (\kappa(\vphi_0) - \kappa(\vphi_1)) \nabla \TA(\vphi_0) \big\|_{\mathcal{L}(H^1, L^2)} + \big\| \kappa(\vphi_1) \nabla [\TA(\vphi_0) - \TA(\vphi_1) ] \big\|_{\mathcal{L}(H^1, L^2)} \Big) \norm{\theta_1 - \theta_2}_{H^1}\\ 
			& \quad + \Big( \big\|  (\kappa(\vphi_1) - \kappa(\vphi_2)) \nabla \TA(\vphi_1) \big\|_{\mathcal{L}(H^1, L^2)} + \big\| \kappa(\vphi_2) \nabla [\TA(\vphi_1) - \TA(\vphi_2) ] \big\|_{\mathcal{L}(H^1, L^2)} \Big) \norm{ \theta_2}_{H^1}
		\end{aligned}\\ 
		& \quad \leq \begin{aligned}[t]
			&C(R) \big( \norm{\vphi_0 - \vphi_1}_{W^{1, q}} \big\|\TA(\vphi_0)\big\|_{\mathcal{L}(H^1)} + \norm{\kappa(\vphi_1)}_{W^{1, q}} \big\|\TA(\vphi_0) - \TA(\vphi_1)\big\|_{\mathcal{L}(H^1)} \big) \norm{\theta_1 - \theta_2}_{H^1}\\ 
			& \quad +C(R) \big( \norm{\vphi_1 - \vphi_2}_{W^{1, q}} \big\|\TA(\vphi_1)\big\|_{\mathcal{L}(H^1)} + \norm{\kappa(\vphi_2)}_{W^{1, q}} \big\|\TA(\vphi_1) - \TA(\vphi_2)\big\|_{\mathcal{L}(H^1)} \big) \norm{ \theta_2}_{H^1}
		\end{aligned}\\ 
	\end{align*}
	Thus, we deduce with Corollary \ref{cor:properties_TA} and Corollary \ref{cor:TA_lipschitz} that 
		\begin{align*}
		& \norm{ \A(\vphi_0)  \big(\theta_1 - \theta_2 \big) -  \big(\A(\vphi_1)\theta_1 - \A(\vphi_2) \theta_2  \big) }_{(H^1)'}\\ 
		& \quad \leq \begin{aligned}[t]
			&C(R) \big(  \norm{ \vphi_0 -\vphi_1}_{W^{1, q}} \norm{\theta_1 - \theta_2}_{H^1} 
			+ \norm{ \vphi_1 -\vphi_2}_{W^{1, q}}  \norm{ \theta_2}_{H^1} \big), 
		\end{aligned}
	\end{align*}
	and exploiting \eqref{embedding:hoelder_non_ve}, we conclude as in \eqref{est:F_1_I_non_ve} 
	\begin{align*}
		I &\leq C(R) T^{\bar{\vartheta}}  
		 \big(  \norm{\vphi_1}_{C^{{\bar{\vartheta}}}(W^{1, q})} \norm{\theta_1 - \theta_2}_{L^s(H^1)}
		 + \norm{\vphi_1 - \vphi_2}_{C^{{\bar{\vartheta}}}(W^{1, q})} \norm{ \theta_2}_{L^s(H^1)} \big) \\ 
		 & \leq C(R) T^{\bar{\vartheta}}  
		 \big( \norm{\theta_1 - \theta_2}_{L^s(H^1)}
		 + \norm{\vphi_1 - \vphi_2}_{C^{{\bar{\vartheta}}}(W^{1, q})}  \big). 
		 \numberthis \label{est:F_2_I_non_ve}
	\end{align*}
	\\ 
	\textit{Ad $\F_T^2.II$:} For the second term we note that $\kappa$ and $m$ have similar properties, such that a result related to Lemma \ref{lemma:general_estimate_non_ve} also holds for $\kappa$. In conjunction with \eqref{eq:C_^-1_lipschitz}, a tedious computation shows
	\begin{align*}
		&\begin{aligned}[t]
			&\Big\| \nabla \cdot \Big( \kappa(\vphi_1) \nabla \Big( \alpha(\vphi_1) M(\vphi_1) \nabla \cdot  \big[ \tilde{\mathcal{C}}^{-1}(\vphi_1) \big( -\nabla(\C(\vphi_1) \Tau(\vphi_1)) + \bm{f }, \bm{g}  \big) \big]   \Big)\Big)\\ 
			&\quad - \nabla \cdot \Big( \kappa(\vphi_2) \nabla \Big( \alpha(\vphi_2) M(\vphi_2) \nabla \cdot  \big[ \tilde{\mathcal{C}}^{-1}(\vphi_2) \big( -\nabla(\C(\vphi_2) \Tau(\vphi_2)) + \bm{f }, \bm{g}  \big) \big]   \Big)\Big)\Big\| _{(H^1)'}\\  
		\end{aligned}\\ 
		& \quad  \leq C(R) \bnorm{\vphi_1 -\vphi_2}_{W^{1, q}}  ( 1+ \norm{\bm{f}}_{\bm{L}^2}+ \norm{\bm{g}}_{\bm{H}^{\sfrac{1}{2}}(\Gamma_N) } ). 
	\end{align*}
	Hence, we infer with help from the embedding result \eqref{embedding:hoelder_non_ve} that 
	\begin{align*}
		II \leq C(R) T^{\frac{1}{s}}  \norm{\vphi_1 - \vphi_2 }_{C^0(W^{1, q})} 
		+ C(R) T^{\bar{\vartheta}} \norm{\vphi_1- \vphi_2}_{C^{ T^{\bar{\vartheta}}}(H^2)}  ( \norm{\bm{f}}_{L^s(\bm{H}^1)}  +  \norm{\bm{g}}_{L^s(\bm{H}^{\sfrac{3}{2}}(\Gamma_N) )}). 
		 \numberthis \label{est:F_2_II_non:ve}
	\end{align*}
	\\ 
	\textit{Ad $\F_T^2.III$:} The Lipschitz assumption \ref{II.A:source_terms} exactly states
\begin{align*}
	III &= 
	\norm{S_f(\vphi_1, \bm{u}_1,  \theta_1) - S_f(\vphi_2, \bm{u}_2,  \theta_2)}_{L^s(L^2)} 
	\leq C(T, R) ( \norm{\vphi_1 - \vphi_2}_{\new{Z}_T^1} + \norm{\theta_1 - \theta_2}_{\new{Z}_T^2}), 
	\numberthis \label{est:F_2_III_non_ve}
\end{align*}
where the constant $C (R, T)$ tends to zero as $T$ goes to zero. 
\par 
\medskip 
	In summary, \eqref{est:F_2_I_non_ve}-\eqref{est:F_2_III_non_ve} lead to 
		\begin{align*}
		&\norm{\F_T^2 (\vphi_1,  \theta_1)  - \F_T^2 (\vphi_2, \theta_2)}_{Y^2_T} 
		\leq C(R, T) \Big( \norm{\vphi_1 - \vphi_2}_{\new{Z}^1_T} 
		+ \norm{\theta_1 - \theta_2}_{\new{Z}^2_T}
		\Big)
		\numberthis \label{est:F_3_non_ve}
	\end{align*}
	for some constant $C(R, T) > 0$ such that $C(R, T) \rightarrow 0$ as $T \rightarrow 0$, which concludes the proof. 
	\end{subequations}
\end{proof}

\medspace
\part{The visco-elastic case}\label{part:ve}

In the second part of this exposition, we turn our attention to the Cahn--Hilliard--Biot system with visco-elastic regularization of Kelvin--Voigt type. While our strategy for establishing short-time existence of a unique solution stays essentially the same, our choice to work under minimal assumptions with respect to spatial regularity leads to many difficulties requiring more advanced theory. \\ 
We begin in Section \ref{sec:assumptions} with a discussion of $\bm{W}^{1 ,p}$-regularity for elliptic systems, thereby explaining our choices on regularity, and give (geometric) examples under which the assumptions we introduce next are satisfied. Rewriting the problem as a fixed-point equation in Section \ref{sec:proof} then allows for the application of Banach's fixed point theorem, yielding the desired result. To this end, Lipschitz estimates for the corresponding right-hand side in non-standard spaces are needed, which can be found in Section \ref{sec:contraction}. Moreover, just like in Part~\ref{part:nve}, it is crucial to establish maximal regularity for the linearized system, with the main difference that we are no longer in a Hilbert space setting, and the necessity to make use of Banach-scales.

\section{Assumptions: visco-elastic case}\label{sec:assumptions} 

As in the previous case, proving Theorem \ref{thm:existence_strong} is only possible under precise assumptions on the domain $\Omega$ as well as the given functions. Since it is our goal to establish (short time) existence of a unique solution under minimal assumptions with respect to spatial regularity, we start with a discussion of elliptic $\bm{W}^{1, q}$-regularity for systems, before presenting the comprehensive list of assumptions.

\addtocontents{toc}{\SkipTocEntry}
\subsection*{$\textbf{W}^{1,q}$-regularity for elliptic systems}

Considering that elliptic equations and systems have been under intensive investigation for several decades, it is only natural that the theory is incredibly rich with extensive results concerning the regularity of solutions. Nevertheless, we want to give a brief overview, highlighting findings most closely related the problems of this exposition, particularly referencing \cite{MR4770785, bociu2016analysis}, which served as a foundation of the following discussion. \\ 
Here, and in the rest of this paper, we assume that the subset $\Gamma_D \subseteq \partial \Omega$ is relatively closed and set $\Gamma_N\coloneqq \partial \Omega \setminus \Gamma_D$. 
\par 
\medskip 
By the Lax-Milgram theorem, using estimates relying on Korn's inequality, operators of the form 
\begin{equation*}
	- \nabla \cdot \big( \C(\vphi) \E(\bm{v}) \big) +1 : 
	\bm{W}^{1, 2}_{\Gamma_D}(\Omega) \rightarrow \bm{W}^{-1, 2}_{\Gamma_D}(\Omega), 
	\quad \bm{v} \mapsto \Big( \bm{w} \mapsto \int_\Omega  \C(\vphi) \E(\bm{v}) : \E(\bm{w}) + \bm{v} \bm{w} \dx \Big)
\end{equation*}
are topological isomorphisms between the spaces $\bm{W}^{1, 2}_{\Gamma_D}(\Omega)$ and $\bm{W}^{-1, 2}_{\Gamma_D}(\Omega)$, and in case $\mathcal{H}^{n-1}(\Gamma_D)  \neq \emptyset$, we can even omit the "$+ 1$". For $q > 2$, these operators are obviously still bijective on $\bm{W}^{-1, q}_{\Gamma_D}(\Omega)$, however, in general the domain is not $\bm{W}^{1, q}_{\Gamma_D}(\Omega)$, c.f. \cite[Thm. 5.6]{Haller15}. Note that there exist perturbation results, e.g. \cite{groger1989aw}, that establish invertibility for some suitably small $q > 2$. \\ 
Due to the counterexample of Shamir \cite{shamir1968regularization}, we know that, in general, $W^{1, 4}(\Omega)$-regularity can \textit{not} be expected even for smooth data if segments with different types of boundary conditions are allowed to meet arbitrarily. However, \textit{some} regularity can be retained by imposing geometric conditions on the domain.\par 
\medskip 
In the scalar case, classical results on $W^{1, q}$-regularity for pure Dirichlet or Neumann boundary conditions can be found in \cite[Ch. 15]{MR125307} and \cite[pp. 156-157]{MR202511}, see also \cite{MR2314060} and the references therein for more general assumptions. \\ 
For mixed boundary conditions where the different parts meet, $W^{1, q}$-regularity for some $q > 3$ and a wide range of geometric constellations relevant to real world problems is established in \cite{MR3342697}. Moreover, Savaré (\cite{MR1452171}) showed $H^{\frac{3}{2}-\eps}(\Omega)$-regularity for solutions, assuming only that the boundary is of class $C^{1, 1}$. \par 
Since classical De Giorgi-Moser-Nash iteration is no longer applicable, the situation is even more difficult for systems. Here, we would like to reference results on the Lamé system in $\bm{H}^{\frac{3}{2} + \eps }(\Omega)$ over polyhedral domains where the different boundary parts do not meet tangentially by Nicaise \cite{MR1205404}. \\ 
Considering only pure Dirichlet boundary conditions or when $\overline{\Gamma}_D \cap \overline{\Gamma}_N = \emptyset$, higher regularity can be recovered, see \cite{mclean2000, Valent_Boundary, MR3038752, MR1790411}. More recent results on $\bm{W}^{1, q}$-regularity include \cite{MR2399163} (clamped boundary) and \cite{MR2880228} (Neumann boundary, restrictions on $q$), see also the references cited therein. 
\par 
\medskip 

\begin{enumerate}[label = {(II.A\arabic*)} ]
	\item  \label{A:domain}
	Let $\Omega \subset \mathbb{R}^n$ be a bounded and connected $C^{4}$-domain in dimension $n \leq 3$. 
	 Additionally, we assume that the subset $\Gamma_D \subset \partial \Omega$ is relatively closed and satisfies $\mathcal{H}^{n-1} (\Gamma_D) > 0$. Finally, we set $\Gamma_N\coloneqq \partial \Omega \setminus \Gamma_D$. 
	\item  \label{A:psi} 
	The potential $\psi : \R \rightarrow \R$ is of class $\psi \in C^2(\R)$.  
	\item \label{A:W}
	The elastic free energy density $W \in C^1(\mathbb{R} \times \mathbb{R}^{n \times n}_{sym})$ is of the form 
	\begin{equation*}
		W(\vphi', \E') = \C(\vphi') (\E' - \Tau(\vphi')):  (\E' - \Tau(\vphi')), \quad \textnormal{for all } \vphi' \in \R,\, \E' \in \R^{n \times n}_{sym}, 
	\end{equation*} 
	where $\C: \R \rightarrow \mathcal{L}(\R^{n \times n}_{sym})$ is a differentiable tensor-valued function whose derivative $\C'$ is locally Lipschitz continuous and bounded. We require it to fulfill the standard assumptions of linear elasticity, i.e., $\C(\vphi')$ is  symmetric and uniformly positive definite on $\R^{n \times n}_{sym}$, mapping symmetric matrices to symmetric matrices, i.e., 
	\begin{align*}
		\E: \C(\vphi') \E &\geq c |\E|^2,\\ 
		\mathcal{D} : \C(\vphi') \E &= \C(\vphi')\mathcal{D} : \E
	\end{align*}
	for all symmetric matrices $\E, \mathcal{D} \in \R^{n\times n}_{sym}$ and all $\vphi' \in \R$.  \\ 
	The eigenstrain $\Tau : \R \rightarrow \R^{n \times n}_{sym}$ is a differentiable, matrix-valued function with locally Lipschitz continuous and bounded derivative $\Tau'$.
	\item \label{A:C_nu} The modulus of visco-elasticity $\C_{\nu}$ satisfies the same assumptions as the elasticity tensor $\C$. 
	\item \label{A:iso}
	We assume that for all $\vphi \in C^\beta (\overline{\Omega})$, $\beta > 0$, the operators 
	\begin{align*}
		\mathcal{B}(\vphi) : \X (\Omega) \rightarrow \Xd (\Omega), \quad &\bm{v} \mapsto  \nabla \cdot \big( \C_\nu(\vphi) \E(\bm{v}) \big), \\ 
			\mathcal{C}(\vphi) : \X (\Omega) \rightarrow \Xd(\Omega), \quad &\bm{v} \mapsto \nabla \cdot \big( \C(\vphi) \E(\bm{v}) \big)
	\end{align*}
	are topological isomorphisms with inverses that are uniformly bounded in $\norm{\vphi}_{C^\beta}$. Here, $\X(\Omega)$ is a subspace of $\bm{W}^{1, q}(\Omega)$, where $q = 3$ if $n = 2$ and $q > 3$ for some $\delta > 0$ if $n = 3$, cf. Remark \ref{rem:bc}. 
	\item \label{A:g}
	The function $\bm{g} : \Gamma_N \rightarrow \mathbb{R}^n$, modeling applied outer forces, fulfills $\bm{g} \in  L^r(0, T; \bm{W}^{- (1-\frac{1}{q'}), q}(\Gamma_N))$.
	\item \label{A:f}
	The function $\bm{f} \colon \Omega \rightarrow \R^n$, modeling body forces, satisfies $\bm{f} \in \bm{L}^r(0, T; \bm{W}^{-1, q}_{\Gamma_D}(\Omega))$.
	\item \label{A:m}
	There exists a constant $\underline{m} > 0$ such that the mobility $m \in C^{2}(\mathbb{R})$ fulfills 
	\begin{equation*}
		\underline{m} \leq m(z) \quad \textrm{for all} \quad z \in \mathbb{R}. 
	\end{equation*}
	\item \label{A:kappa}
	There exist a constant $\underline{\kappa} > 0$ such that the function $\kappa \in C^{2}(\mathbb{R})$ fulfills 
	\begin{equation*}
		\underline{\kappa} \leq \kappa(z) \quad \textrm{for all} \quad z \in \mathbb{R}. 
	\end{equation*}
	\item \label{A:phase_coefficients}
	The maps $\alpha, M : \R \rightarrow \R$ are of class $C^{3}(\R)$. 
	\item \label{A:source_terms}
	The source terms $S_s, S_f$ are real functions $\R \times \R^{n \times n} \times \R  \rightarrow \R$. Moreover, we require the following Lipschitz-estimates for all $(\vphi_i, \bm{u}_i, \theta_i) \in X_T$ satisfying $\norm{(\vphi_i, \bm{u}_i, \theta_i)}_{X_T} \leq R$
	\begin{align*}
		\hspace{1.5cm}
		\norm{S_s(\vphi_1,  \bm{u}_1, \theta_1) - S_s(\vphi_2,  \bm{u}_2, \theta_2)}_{L^r (W^{-2, p}_N) } 
		&\leq C(R, T) \norm{ (\vphi_1 - \vphi_2, \bm{u}_1 - \bm{u}_2, \theta_1 - \theta_2) }_{X_T}, 
		\\ 
			\norm{S_f(\vphi_1,  \bm{u}_1, \theta_1) - S_f(\vphi_2,  \bm{u}_2, \theta_2)}_{L^r(W^{-2, q}_N)} 
		&\leq C(R, T) \norm{ (\vphi_1 - \vphi_2, \bm{u}_1 - \bm{u}_2, \theta_1 - \theta_2) }_{X_T}, 
	\end{align*}
	where $C(R, T)$ tends to zero as $T$ goes to zero.
\end{enumerate}

\new{
	\begin{figure}[h!]
		\centering
		\begin{tikzpicture}[scale=2, \ifdraft{blue}{}]
			
			\draw[thick]
			(0,0)
			.. controls (1.2,-0.4) and (2.2,-0.2) .. (2.6,0.6)
			.. controls (3.0,1.4) and (2.2,2.2) .. (1.4,2.4)
			.. controls (0.6,2.6) and (-0.2,1.8) .. (-0.3,1.0)
			.. controls (-0.4,0.4) and (-0.2,0.1) .. (0,0);
			
			\draw[ultra thick]
			(2.6,0.6)
			.. controls (3.0,1.4) and (2.2,2.2) .. (1.4,2.4);
			
			\node at (1.2,1.2) {$\Omega$};
			\node[above right] at (2.2,2.2) {$\Gamma_D$};
			\node[left] at (-0.4,1.4) {$\Gamma_N$};
			
		\end{tikzpicture}
		\caption{Sketch of an admissible domain in 2D}
	\end{figure}
}

We end this section with some comments on our assumptions. 

\begin{remark}\label{rem:bc} 
	\begin{enumerate}[label*=(\roman*)]
		\item The result in \cite[Thm. 4.1]{Valent_Boundary} shows that \ref{A:iso} is certainly satisfied if we restrict ourselves to clamped boundary conditions, i.e., $\Gamma_D = \partial \Omega$. Combined with the  proof of \cite[Thm. 5.3]{Valent_Boundary} for the Neumann case, the assertion can also be derived for mixed boundary conditions (i.e., $\mathcal{H}^{n-1} (\Gamma_D), \mathcal{H}^{n-1} (\Gamma_N)  > 0$) where $\overline{\Gamma}_D \cap \overline{\Gamma}_N = \emptyset$. 
		For a rigorous proof, one can proceed exactly as in Theorem \ref{thm:elastic_non_symmetric} after deriving the analogous result to Theorem \ref{thm:w^2_p_non_div} in the space $\bm{W}^{1, p}(\Omega)$ instead of $\bm{W}^{2, p}(\Omega)$. For the latter we refer to \cite[Thm. 6.3.6]{PrussSimonett2016} and the proof of Theorem \ref{thm:w^2_p_non_div}. \\
		Unlike the scalar case (cf. \cite{MR3342697}), to our best knowledge, there are no explicit geometric assumptions or compatibility conditions known for which \ref{A:iso} holds if the different boundary parts meet and $\Omega$ is not a polyhedron. In that regard, our exposition should be understood as motivation for further investigation. 
		\item \label{remark:assumptions}
		For any $\bm{v} \in \bm{W}^{1, q'}_{\Gamma_D}(\Omega)$, well-known trace theorems imply $ \bm{v}_{| \Gamma_N}  \in \bm{W}^{1-  \frac{1}{q'}, q'}_{\Gamma_D}(\Gamma_N)$ and we deduce the embedding 
			\begin{equation}\label{eq:trace_dual_space}
				\bm{W}^{- (1-\frac{1}{q'}), q}(\Gamma_N) \hookrightarrow \bm{W}^{-1, q}_{\Gamma_D}(\Omega). 
			\end{equation}
		In particular, the sum $\bm{f} + \bm{g}$ is well-defined, since we assume that $\bm{f} \in \bm{L}^r(0, T; \new{\bm{W}^{-1, q}_{\Gamma_D}(\Omega)})$ and $\bm{g} \in  L^r(0, T; \bm{W}^{- (1-\frac{1}{q'}), q}(\Gamma_N)) \hookrightarrow  L^r(0, T; \bm{W}^{-1, q}_{\Gamma_D}(\Omega)) $. 
		\item The assumptions above further imply the following growth conditions:  \\ 
		\begin{enumerate}[label = (\roman*)]
			\myitem{(II.A3.2)}\label{A:W_growth_conditions}
			There exists a constant $C_2 > 0$ such that for all $\vphi' \in \mathbb{R}$ and all symmetric $\E' \in  \mathbb{R}^{n \times n}$,
			\begin{align*}
				\abs{W(\vphi', \E') } & \leq C_2 \left( \abs{\E'}^2 + \abs{\vphi'}^2 +1 \right), \\
				\abs{W_{,\vphi}(\vphi', \E') } & \leq C_2 \left( \abs{\E'}^2 + \abs{\vphi'}^2 +1 \right), \\
				\abs{W_{,\E}(\vphi', \E')} & \leq C_2 \left( \abs{\E'} + \abs{\vphi'} +1 \right). 
			\end{align*}
		\end{enumerate}
	\end{enumerate}
	
\end{remark}

\medspace

\section{Proof of Theorem \ref{thm:existence_strong}} \label{sec:proof}

The strategy of this proof is essentially identical to Part \ref{part:nve}: after eliminating the functions $p$ and $\mu$ from the system by substituting the corresponding expression for the respective variables, we linearize the highest-order evolution operators, aiming to reduce the question of well-posedness to a fixed-point argument.\\ 
In contrast to previous arguments, we do not eliminate $\bm{u}$ but use elliptic regularity theory, or rather Assumption \ref{A:iso}, to infer the existence of the inverse operator to $\mathcal{B} (\vphi_0)$, and derive an evolution equation governed by a bounded operator for the displacement.\\ 
After showing maximal regularity for the linearized system and deriving Lipschitz estimates for the corresponding right-hand sides in the respective spaces, the existence of a unique solution will be established by means of a contraction principle.  

\addtocontents{toc}{\SkipTocEntry}
\subsection*{The spaces}
In analogy to the notation of the previous part, we start by introducing relevant function spaces, setting
\begin{align*}
	Z^1_T &\coloneqq W^{1, r} (0, T; W^{-2, p}_N (\Omega)) \cap L^r (0, T; W^{2, p}_N (\Omega)),
	 \\
	Z^2_T &\coloneqq W^{1, r} (0, T; \X (\Omega)), \\ 
	Z^3_T & \coloneqq W^{1, 3r} (0, T; W^{-2, q}_N (\Omega)) \cap L^{3r} (0, T; L^q(\Omega)),
\end{align*}
where 
\begin{equation*}
	W^{2, p}_N (\Omega) \coloneqq \{ \vphi \in W^{2, p} \colon \partial_n \vphi = 0   \textnormal{ on } \partial \Omega \}, 
\end{equation*}
and $W^{-2, p}_N (\Omega))  \coloneqq \big(W^{2, p'}_N (\Omega)\big)'$. These spaces are equipped with the norms $\norm{\cdot}_{Z^1_T}$, $\norm{\cdot}_{Z^2_T}$ and $\norm{\cdot}_{Z^3_T}$, respectively, defined by
\begin{align*}
	\norm{\vphi}_{Z^1_T} &\coloneqq \norm{\pt \vphi}_{L^r(0, T; W^{-2, p}_N(\Omega))} + \norm{\vphi}_{L^r(0, T; W^{2, p}_N(\Omega))} + \norm{\vphi(0)}_{  (W^{-2, p}_N(\Omega), W^{2, p}_N(\Omega))_{1- \frac{1}{r}, r}}, \\ 
	\norm{\bm{u}}_{Z^2_T} & \coloneqq \norm{\bm{u}}_{W^{1, r} (0, T; \X (\Omega))} + \norm{\bm{u}(0)}_{\X (\Omega)}, \\ 
	\norm{\theta}_{Z^3_T} &\coloneqq \norm{\pt \theta}_{L^{3r}(0, T; W^{-2, q}_N(\Omega))} + \norm{\theta}_{L^{3r}(0, T; L^q(\Omega))} + \norm{\theta(0)}_{  (W^{-2, q}_N(\Omega), L^q(\Omega))_{1- \frac{1}{3r}, 3r}}. 
\end{align*}
As before, Lemma \ref{lem:uniform_spaces} allows us to find uniform estimates in $T  < \frac{T_0}{2}$ for any given $T_0 > 0$. Given these definition, we introduce the affine subspaces 
\begin{align*}
	X^1_T &\coloneqq \{  \vphi \in Z^1_T \colon \vphi(0) = \vphi_0 \} , \\  
	X^2_T &\coloneqq \{  \bm{u} \in Z^2_T \colon \bm{u}(0) = \bm{u}_0 \},  \\  
	X^3_T &\coloneqq \{  \theta \in Z^3_T \colon \theta(0) = \theta_0 \}
\end{align*}
and set $X_T \coloneqq X^1_T \times X^2_T \times X^3_T$. Similarly, let 
\begin{equation*}
	Y_T \coloneqq Y^1_T \times Y^2_T \times Y^3_T \coloneqq L^r(0, T; W^{-2, p}_N(\Omega)) \times L^r (0, T; \Xd (\Omega)) \times L^{3r}(0, T; W^{-2, q}_N(\Omega)). 
\end{equation*}
\par 

\addtocontents{toc}{\SkipTocEntry}
\subsection*{The fixed-point problem}
After this preparation, we are finally able to start with the proof by defining the mapping 
\begin{equation}\label{eq:linearized_viscous}
		\cL_T : X_T \rightarrow Y_T,  \quad \quad 
		\cL_T (\vphi, \bm{u}, \theta) = 
		\left( 
		\begin{array}{c}
			\pt \vphi + \eps \Delta  (m(\vphi_0)  \Delta \vphi)  \\ 
			\pt  \bm{u} + \mathcal{A}(\vphi_0) \bm{u} \\ 
			\pt \theta - \nabla \cdot (\kappa(\vphi_0) M(\vphi_0)  \nabla \theta) 
		\end{array}
		\right), 
\end{equation}
where the operator $\mathcal{A}$ is given by 
\begin{align*}
	\mathcal{A}(\vphi) &: \X(\Omega) \rightarrow \X(\Omega), 
	\quad \bm{v} \mapsto \mathcal{B}^{-1}(\vphi) \mathcal{C}(\vphi) \bm{v}. 
\end{align*}
Note that the mapping above is well-defined on account of the multiplication results in Section~\ref{sec:preliminaries}.
Moreover, we set 
\begin{align*}
	&\F_T^1 (\vphi, \bm{u}, \theta) = 
	\begin{aligned}[t]
		&\eps \Delta  (m(\vphi_0)  \Delta \vphi)  - \eps \nabla \cdot (m(\vphi) \nabla \Delta \vphi) \\ 
		& + \nabla \cdot \Big ( m(\vphi) \nabla \Big( \frac{1}{\eps} \psi'(\vphi) + W_{,\vphi}(\vphi, \E(\bm{u}))  \Big) \Big) \\ 
		& + \nabla \cdot \Big ( m(\vphi) \nabla \Big(
		\frac{M'(\vphi)}{2} (\theta - \alpha(\vphi) \nabla \cdot \bm{u} )^2 - M(\vphi) (\theta - \alpha(\vphi) \nabla \cdot \bm{u} ) \alpha'(\vphi) \nabla \cdot \bm{u}  \Big) \Big) \\ 
		&+ S_s(\vphi, \E(\bm{u}), \theta), 
	\end{aligned}  \\ 
	&\F_T^2 (\vphi, \bm{u}, \theta) = \big(  \mathcal{A}(\vphi_0) - \mathcal{A}(\vphi) \big) \bm{u} + 
		\mathcal{B}^{-1} (\vphi) \Big( \nabla \cdot \big( \C(\vphi) \Tau (\vphi)  + \alpha(\vphi) M(\vphi) (\theta - \alpha(\vphi) \nabla \cdot \bm{u})  \textbf{I} \big)+ \bm{f}+  \bm{g}\Big),  \\ 
	&\F_T^3 (\vphi, \bm{u}, \theta) = 
		\begin{aligned}[t]
			& - \nabla \cdot \Big( \kappa(\vphi_0) M(\vphi_0)  \nabla \theta \Big)
			+ \nabla \cdot  \Big( \kappa(\vphi) \nabla \big( M(\vphi)  \theta \big) \Big)   \\ 
			&  \quad + \nabla \cdot \Big(\kappa(\vphi) \nabla \big(M(\vphi) \alpha(\vphi) \nabla \cdot \bm{u} \big) \Big) + S_f(\vphi, \E(\bm{u}), \theta), 
		\end{aligned}
\end{align*}
and combine these to obtain the non-linear mapping $\F_T : X_T \rightarrow Y_T$ defined by
\begin{equation*}
	\F_T(\vphi, \bm{u}, \theta) = \left(
		\begin{array}{c}
			\F_T^1 (\vphi, \bm{u}, \theta) \\ 
				\F_T^2 (\vphi, \bm{u}, \theta) \\ 
					\F_T^3 (\vphi, \bm{u}, \theta)
		\end{array}
	\right)
\end{equation*}
for all $(\vphi, \bm{u}, \theta) \in X_T$. 
The application of Banach's fixed point theorem is only possible if we can show that the right-hand side describes a contraction, which is the result of the following proposition.  

\begin{proposition}\label{prop:contraction}
	There exists a constant $C(T, R) > 0$ such that for all $(\vphi_i, \bm{u}_i, \theta_i) \in X_T$, $i = 1, 2$, satisfying $\norm{(\vphi_i,\bm{u}_i, \theta_i)}_{X_T} \leq R$ it holds that 
	\begin{equation}\label{eq:contraction}
		\norm{\F_T(\vphi_1, \bm{u}_1, \theta_1) - \F_T(\vphi_2, \bm{u}_2, \theta_2)}_{Y_T} \leq C(T, R) \norm{(\vphi_1 - \vphi_2, \bm{u}_1 - \bm{u}_2, \theta_1 - \theta_2)}_{\new{Z}_T}
	\end{equation}
	Moreover, for all $R > 0$, we have $C(R, T) \rightarrow 0$ as $T \rightarrow 0$. 
\end{proposition}

The proof of this proposition is the subject of Section~\ref{sec:contraction}. \\ 
Rewriting the problem above as a fixed-point equation further requires invertibility of the evolution operators in \eqref{eq:linearized_viscous}. 

\begin{lemma}\label{lem:inverse_boundedness}
	Let $\cL_T, X_T$ and $Y_T$ be defined as above. Then, the operator $\cL_T : X_T \rightarrow Y_T$ is invertible for every $T > 0$ and for any fixed $T_0 > 0$ there exists some constant $C_{\cL^{-1}_T}(T_0) > 0$ such that \new{
		\begin{align*}
			&\norm{ \cL^{-1}_T (f_1, f_2, f_3)}_{Z_T} \\ 
			 &\quad \leq C_{\mathcal{L}^{-1}_T}(T_0) \Big(  
			 \begin{aligned}[t]
			 	  &\norm{ (f_1, f_2, f_3)}_{Y_T} \\
			 	& \quad +\norm{\vphi_0}_{  (W^{-2, p}_N(\Omega), W^{2, p}_N(\Omega))_{1- \frac{1}{r}, r}}
			 	+\norm{\bm{u}_0}_{\X (\Omega)}
			 	+ \norm{\theta_0}_{  (W^{-2, q}_N(\Omega), L^q(\Omega))_{1- \frac{1}{3r}, 3r}}\Big) 
			 \end{aligned}
		\end{align*}
		and 
		\begin{equation*}
			\big\|  \cL^{-1}_T (f_1, f_2, f_3) -  \cL^{-1}_T (\tilde{f_1}, \tilde{f_2}, \tilde{f_3}) \big\|_{Z_T} \leq 
			C_{\mathcal{L}^{-1}_T}(T_0) \big\| (f_1, f_2, f_3) - (\tilde{f_1}, \tilde{f}_2, \tilde{f}_3) \big\|_{Y_T}
		\end{equation*}
		for all $T \in (0, T_0]$ and  $(f_1, f_2, f_3),(\tilde{f_1}, \tilde{f_2}, \tilde{f}_3)  \in Y_T$. 
	}
\end{lemma}

The proof of this lemma can be found in Section~\ref{section:inverse}. \\ 
Relying on these two results, it is now possible to prove Theorem~\ref{thm:existence_strong}. 

\begin{proof}[Proof of Theorem \ref{thm:existence_strong}] \label{proof:main_visco_elastic}
		Using Proposition \ref{prop:contraction} and Lemma \ref{lem:inverse_boundedness}, this result follows analogously to the proof of Theorem \ref{thm:existence_strong_nve}, see Section \ref{proof:main_visco_elastic_nve}. 
\end{proof}

\medspace

\section{Existence and continuity of $\mathcal{L}^{-1}$: visco-elastic case} \label{section:inverse}

The aim of this section is to prove Lemma \ref{lem:L^-1_uniform_T}, i.e., the existence and continuity of the solution operator $\mathcal{L}^{-1}$. To this end, we will consider the components  $\mathcal{L}_i$, $i = 1, 2, 3$, of $\mathcal{L}$ separately, showing that these are invertible and that their inverses can be bounded independently of the time interval. \\ 
While the evolution in the second component is governed by a bounded operator, thus exhibiting maximal regularity on all bounded intervals, treatment of $\cL_1$ and $\cL_3$ requires more subtle techniques. For these, an application of known results also yields maximal regularity, and even implies that they admit a bounded $\mathcal{H}^\infty$-calculus, albeit not in the right spaces, hence necessitating an application of Banach-scales. \\ 
Note that the following proof heavily relies on definitions and results about the $\mathcal{H}^\infty$-calculus, $\mathcal{R}$-sectorial operators, operators with bounded imaginary powers and Banach-scales. A summary of necessary  prerequisites can be found in Section~\ref{sec:preliminaries}. 

\begin{lemma}\label{lemma:L_1}
	Suppose $\Omega \subset \R^n$ is a bounded domain with $C^4$-boundary, \new{ $1 < r < \infty$}, and assume that \ref{A:m} holds. Let $\vphi_0  \in \big(W^{-2, p}_N(\Omega), W^{2, p}_N(\Omega)\big)_{1- \frac{1}{r}, r}$,  and $f \in L^r(0, T; W^{-2, p}_{N}(\Omega))$, along with coefficients satisfying $\tilde{\varphi }\in W^{1, q}(\Omega) \cap W^{2, p}(\Omega)$. Then, for every ${0 < T  < \infty}$, there exists a unique
	\begin{equation*}
		\vphi \in   L^r (0, T; W^{2, p}_N (\Omega)) \cap W^{1, r} (0, T; W^{-2, p}_N (\Omega))
	\end{equation*}
	such that 
	\begin{alignat*}{2}
		\pt \vphi + \eps \Delta  \big( m(\tilde{\varphi })  \Delta \vphi \big)  &= f  \quad &&\textnormal{in } (0, T) \times \Omega, \numberthis \label{eq:ch_strong}\\ 
		\vphi_{| t = 0 } &= \vphi_0  \quad &&\textnormal{in } \Omega,
	\end{alignat*}
	where \eqref{eq:ch_strong} holds in a $W^{-2, p}_N$-sense, i.e., for almost all $t \in (0, T)$ it holds that 
	\begin{align*}
		{}_{W^{-2, p}_N} \langle \pt \vphi, \zeta \rangle_{W^{2, p'}_N} 
		+   \int_\Omega \eps m(\tilde{\varphi }) \Delta\vphi \,  \Delta \zeta \dx = {}_{W^{-2, p}_N} \langle f , \zeta \rangle_{W^{2, p'}_N}  \quad 
		\textnormal{for all } \zeta \in W^{2, p'}_N(\Omega). 
	\end{align*}
	\new{
	Denoting by $\mathcal{L}_{1, T}^{-1}$ the solution operator, then there exists a constant $C(T) > 0$ such that 
	\begin{align*}
		\norm{ \cL^{-1}_{1, T} (f) }_{Z_T^1}  \leq C(T) \Big(   \norm{ f}_{Y_T^1} + \norm{\vphi_0}_{  (W^{-2, p}_N(\Omega), W^{2, p}_N(\Omega))_{1- \frac{1}{r}, r}}
		 \Big) 
	\end{align*}
	for all $f \in Y_T^1$. 
	}
\end{lemma}

\begin{proof}
	We proceed as in the proof of Lemma \ref{lemma:L_1_non_ve}. 
	First of all, consider the operator 
	\begin{equation*}
		\tilde{\mathcal{L}}_1 :  
		D (	\tilde{\mathcal{L}}_1 ) = \{  u \in W^{4, p}(\Omega) : \partial_{\bm{n}} u_{|\partial \Omega} =0 = \partial_{\bm{n}} \Delta u_{|\partial \Omega} \}
		\rightarrow 
		L^p(\Omega), \quad \vphi \mapsto  \eps \Delta  \big( m(\tilde{\varphi })  \Delta \vphi \big), 
	\end{equation*}
	which we may rewrite as 
	\begin{equation*}
		\tilde{\mathcal{L}}_1 \vphi = \eps m(\tilde{\varphi }) \Delta^2 \vphi 
		+ \eps m'(\tilde{\varphi }) \Delta \tilde{\vphi} \,  \Delta \vphi 
		+ \eps m''(\tilde{\vphi}) \nabla \tilde{\vphi} \cdot \nabla \tilde{\vphi}\,  \Delta \vphi
		+  2 \eps m'(\tilde{\varphi }) \nabla \tilde{\vphi} \cdot \nabla  \Delta \vphi .  
	\end{equation*}
	Since $q > n$, standard embedding theorems imply 
	\begin{equation*}
		\tilde{\varphi }  \in W^{1, q}(\Omega) 
		 \hookrightarrow C^\beta(\bar{\Omega}), \quad 
		 \textnormal{for some } \beta > 0, 
	\end{equation*}
	and we deduce with the help of assumption \ref{A:m} the regularity
	\begin{equation*}
			m(\tilde{\varphi }) \in C^\beta(\bar{\Omega})  \quad \textnormal{as well as}\quad 
			 m'(\tilde{\varphi }) \Delta \tilde{\varphi },  m''(\tilde{\vphi}) \nabla \tilde{\vphi} \cdot \nabla \tilde{\vphi} \in L^{p}(\Omega)
			  \quad \textnormal{and }\quad  m'(\tilde{\vphi}) \nabla \tilde{\vphi} \in L^{q}(\Omega). 
	\end{equation*}
	Since $m(\cdot)$ is a uniformly positive, real valued function, the principal symbol
	\begin{equation*}
		\mathcal{A}_{\#} (x, \xi) = \sum_{i, j = 1}^n m(\tilde{\varphi })(\bm{x}) \, \xi_i^{2} \xi_j^2
	\end{equation*}
	is, independently of $\bm{x} \in \bar{ \Omega}$, parameter elliptic admitting an arbitrarily small angle of ellipticity \linebreak$0 < \phi_{\mathcal{A}_\#} < \tfrac{\pi}{2}$. 
	Now it is easy to verify the smoothness and ellipticity conditions required in Theorem~\ref{thm:H_infty_calculus}, which yields that for any $\phi < \tfrac{\pi}{2}$ there exists some $\mu_1 \geq 0$ such that $\tilde{\mathcal{L}}_1 + \mu_1 \in \mathcal{H}^\infty(L^p(\Omega))$ with $\phi^\infty_{\tilde{\mathcal{L}}_1 + \mu_1 } \leq \phi$. Thus, this operator has bounded imaginary powers with power angle $\theta_{\tilde{\mathcal{L}}_1 + \mu_1} \new{\leq} \phi^\infty_{\tilde{\mathcal{L}}_1 + \mu_1 } < \tfrac{\pi}{2}$, i.e., $\tilde{\mathcal{L}}_1 + \mu_1 \in \mathcal{BIP}(L^p(\Omega))$. As we have observed before, this implies that the extrapolated fractional power scale of order $m$ generated by $(L^p(\Omega), \tilde{\mathcal{L}}_1 + \mu_1 ) \eqqcolon (E, \prescript{}{p}{A})$ and the interpolation-extrapolation scale $[(E_\alpha,  \prescript{}{p}{A}_\alpha)\ ;\ \alpha \in [-m, \infty )]$ generated by $(E, \prescript{}{p}{A} )$ and the complex interpolation operator $[\cdot, \cdot]_{\vartheta}$ are equivalent, cf. Section~\ref{subsec:scales}, \cite[Theorem V 1.5.4]{amann1995linear}. 
	Let us therefore turn our attention to the case $\alpha = -\tfrac{1}{2}$, for which it holds that 
	\begin{equation*}
		 \prescript{}{p}{A}_{- \frac{1}{2}} : D( \prescript{}{p}{A}_{- \frac{1}{2}}) = E_{\frac{1}{2}} \subset  E_{- \frac{1}{2}}\rightarrow  E_{- \frac{1}{2}}, 
	\end{equation*}
	where $ \prescript{}{p}{A}_{- \frac{1}{2}}$ is the $E_{- \frac{1}{2}}$-realization of $ \prescript{}{p}{A} = \tilde{\mathcal{L}}_1 + \mu_1$. As in Lemma \ref{lemma:L_1_non_ve}, we use complex interpolation and results by Seeley \cite{seeley1972interpolation}, see also \cite[VII 2.4]{amann2013analysis}, to compute 
	\begin{align*}
		E_{\frac{1}{2}} = [E_0, E_1]_{\frac{1}{2}} = [L^p(\Omega), D({}_\new{p}A)]_{\frac{1}{2}} 
		&= [L^p(\Omega),\{  u \in W^{4, p}(\Omega) : \partial_n \vphi_{|\partial \Omega} =0 = \partial_n \Delta \vphi_{|\partial \Omega} \}]_{\frac{1}{2}} \\ 
		&= \{  u \in W^{2, p}(\Omega) : \partial_n \vphi_{|\partial \Omega} = 0 \} = W^{2, p}_N(\Omega). 
	\end{align*}
	In remains to determine $E_{- \frac{1}{2}}$. In Section~\ref{subsec:scales}, see also \cite[Theorem V 1.4.12]{amann1995linear} we have already observed that $(E_\alpha)' = E^\#_{- \alpha}$ for all $\alpha \in \R$, where $E^\#_{-\alpha}$ is part of the extrapolated fractional power scale generated by the dual pair $(E^\#,  \prescript{}{p}{A}^\#)$, i.e., $[(E_\alpha^\#,  \prescript{}{p}{A}_\alpha^\#) \ ; \ \alpha \in [-m, \infty )]$. 
	Since the dual operator of $\prescript{}{p}{A}$ is given by 
	\begin{align*}
		\prescript{}{p}{A}^\# &= \prescript{}{p'}{A} = \tilde{\mathcal{L}}_1 + \mu_1 :  
		D (	\prescript{}{p'}{A} ) = \{  u \in W^{4, p'}(\Omega) : \partial_n \vphi_{|\partial \Omega} =0 = \partial_n \Delta \vphi_{|\partial \Omega} \} \subset L^{p'}(\Omega)
		\rightarrow 
		L^{p'}(\Omega),  \\
		\vphi &\mapsto  \eps \Delta  \big( m(\tilde{\varphi })  \Delta \vphi \big) + \mu_1 \vphi, 
	\end{align*}
	and satisfies $\prescript{}{p'}{A}_\alpha \in \mathcal{BIP}(L^{p'}(\Omega))$ with $\theta_{\prescript{}{p'}{A}} < \tfrac{\pi}{2}$, cf. \cite[Thm. V 1.4.11]{amann1995linear}, we can argue as before and see $E^\#_{\frac{1}{2}} = [E^\#_0, E^\#_1 ]_{\frac{1}{2}}$. Interpolation therefore yields 
	\begin{align*}
		E_{-\frac{1}{2}}  = (E^\#_{\frac{1}{2}})' 
		&= [(E^\#_0),( E^\#_1)]_{\frac{1}{2}} 
			=  [L^{p'}(\Omega),\{  u \in W^{4, p'}(\Omega) : \partial_n \vphi_{|\partial \Omega} =0 = \partial_n \Delta \vphi_{|\partial \Omega} \}]'_{\frac{1}{2}}\\ 
		& = (W^{2, p'}_N(\Omega))' 
		= W^{-2, p}_N(\Omega), 
	\end{align*}
	leading to 
	\begin{equation*}
		\prescript{}{p}{A}_{- \frac{1}{2}} :  W^{2, p}_N(\Omega) \subset  W^{-2, p}_N(\Omega) \rightarrow  W^{-2, p}_N(\Omega). 
	\end{equation*}
	By \cite[Prop. V 1.5.5]{amann1995linear} this operator also exhibits bounded imaginary powers with power angle $\theta_{	\prescript{}{p}{A}_{- \frac{1}{2}}} < \tfrac{\pi}{2}$, i.e., $\prescript{}{p}{A}_{- \frac{1}{2}} \in \mathcal{BIP}( W^{-2, p}_N(\Omega))$. \\ 
	Finally, we can return to the operator 
			\begin{equation}\label{eq:maxreg_inverse_homogeneous}
			\hat{\mathcal{L}}_1 :  
			W^{2, p}_N(\Omega) \subset  W^{-2, p}_N(\Omega) \rightarrow  W^{-2, p}_N(\Omega), \quad \vphi \mapsto  \eps \Delta  \big( m(\tilde{\vphi})  \Delta \vphi \big). 
		\end{equation} 
		By treating cases $\mu_1 = 0$ and $\mu_1 > 0$ exactly as in Lemma \ref{lemma:L_1_non_ve}, we obtain maximal regularity for the problem  
	\begin{equation*} 
		\pt \vphi + 	\hat{\mathcal{L}_1} \vphi = f(t), \quad t \in (0, T), \quad \vphi(0) = 0, 
	\end{equation*}
	where $f \in L^r(0, T;  W^{-2, p}_N(\Omega))$ \new{along with the estimate 
	\begin{align*}
		\norm{ \hat{\cL}^{-1}_{1, T} (f) }_{Z_T^1}  \leq C(T)  \norm{ f}_{Y_T^1} 
	\end{align*}
	due to the bounded inverse theorem. Here, $ \hat{\cL}^{-1}_{1, T}$ denotes the solution operator and $C(T) >$ is a constant depending on $T > 0$.} Using a standard argument, cf. \cite[Prop. 1.3]{ARENDT20071}, \cite[Lemma 10]{haselboeck2024existence}, we find that the corresponding Cauchy problem with non-trivial initial data satisfying $\vphi_0 \in \big(W^{-2, p}_N(\Omega), W^{2, p}_N(\Omega)\big)_{1- \frac{1}{r}, r}$ is also well-posed, \new{and that the asserted estimate holds}. 
\end{proof}

Adapting parts of the proof above, we proceed by showing a similar result for $\mathcal{L}_3$. 

\begin{lemma}\label{lemma:L_3}
		Suppose $\Omega \subset \R^n$ is a domain with compact $C^4$-boundary and assume that \ref{A:kappa}, \ref{A:phase_coefficients} hold, along with $1 < r < \infty$. Let the initial value $\theta_0  \in \big(W^{-2, q}_N(\Omega), L^q(\Omega)\big)_{1- \frac{1}{3r}, 3r}$ and the right hand side $f \in L^{3r}(0, T; W^{-2, \new{q}}_{N}(\Omega))$ be given, and fix some $\vphi_0 \in W^{1, q}(\Omega)$ for some $q > n$. Then, for every $0 < T  < \infty$ there exists a unique
	\begin{equation*}
		\theta \in   L^{3r}(0, T; L^q (\Omega)) \cap   W^{1, 3r} (0, T; W^{-2, q}_N (\Omega))
	\end{equation*}
	such that 
	\begin{alignat*}{2}
		\pt \theta -  \nabla \cdot \big( \kappa(\vphi_0) M(\vphi_0) \nabla \theta \big) &= f  \quad &&\textnormal{in } (0, T) \times \Omega, \numberthis \label{eq:fluid_strong}\\ 
		\theta_{|t =0} &= \theta_0  \quad &&\textnormal{in } \Omega,
	\end{alignat*}
	where \eqref{eq:fluid_strong} holds in a $W^{-2, \new{q}}_N$-sense, i.e., for almost all $t \in (0, T)$ it holds that
	\begin{equation*}
		{}_{W^{-2, \new{q}}_{N}} \langle \pt \theta, \xi \rangle_{W^{2, \new{q}}_N} 
		- \int_\Omega  \theta\,  \nabla \cdot \big(  \kappa(\vphi_0) M(\vphi_0)  \nabla \xi \big)  \dx = 
			{}_{W^{-2, \new{q}}_{N}} \langle f , \xi \rangle_{W^{2, \new{q}}_N} 
			\quad \textnormal{for all } \xi \in W^{2, \new{q}}_N(\Omega). 
	\end{equation*}
	\new{
		Denoting by $\mathcal{L}_{3, T}^{-1}$ the solution operator, then there exists a constant $C(T) > 0$ such that 
		\begin{align*}
			\norm{ \cL^{-1}_{3, T} (f) }_{Z_T^3}  \leq C(T) \Big(   \norm{ f}_{Y_T^3} + \norm{\theta_0}_{  (W^{-2, q}_N(\Omega), L^q(\Omega))_{1- \frac{1}{3r}, 3r}}
			\Big) 
		\end{align*}
		for all $f \in Y_T^3$. 
	}
\end{lemma}

\begin{proof}
	Similar to above, we need to study the properties of the operator 
	\begin{align*}
		\tilde{\mathcal{L}}_3 &:  
		D (	\tilde{\mathcal{L}}_3 ) = \{  u \in W^{2, q}(\Omega) : \partial_{\bm{n}} u_{|\partial \Omega} = 0  \}
		\rightarrow 
		L^q(\Omega),\\
		 \theta &\mapsto  - \nabla \cdot \big( \kappa(\vphi_0) M(\vphi_0) \nabla \theta \big) = 
		 - \kappa(\vphi_0) M(\vphi_0) \Delta \theta -  \big( \kappa'(\vphi_0)M(\vphi_0) + \kappa(\vphi_0)M'(\vphi_0)  \big)  \nabla \vphi_0 \cdot \nabla \theta, 
	\end{align*}
	where 
	\begin{equation*}
		 \kappa(\vphi_0) M(\vphi_0) \in C^\beta(\overline{\Omega}), \quad 
		  \Big( \kappa'(\vphi_0)M(\vphi_0) + \kappa(\vphi_0)M'(\vphi_0)  \Big)  \new{\nabla \vphi_0}\in L^q(\Omega)
	\end{equation*}
	due to assumptions \ref{A:kappa} and \ref{A:phase_coefficients}. Recall that in Section~\ref{subsec:BIP}, we use the convention $D = i (\partial_1, \ldots, \partial_n)$. Keeping this in mind, we find the principal symbol of this operator given by
	\begin{equation*}
		\mathcal{A}_{\#}  (x, \xi ) = \sum_{i = 1}^n  \kappa(\vphi_0) M(\vphi_0)  \xi_i^2  \quad \textnormal{ for all } x \in \Omega, \xi \in \R^n,  
	\end{equation*} 
	which is parameter elliptic with arbitrarily small angle of ellipticity $0 < \phi_{\mathcal{A}_\#} < \tfrac{\pi}{2}$ independently of $\bm{x} \in \Omega$. Analogously to the proof of Lemma~\ref{lemma:L_1}, and with the help of Theorem~\ref{thm:H_infty_calculus}, we deduce the existence of some $\mu_2 \geq 0$ and $\phi < \frac{\pi}{2}$ such that the operator $\prescript{}{q}{B} =\tilde{\mathcal{L}}_3 + \mu_2 \in \mathcal{BIP}(L^q(\Omega))$ with power angle $\theta_B < \tfrac{\pi}{2}$. Putting $F \coloneqq L^q(\Omega)$, the considerations in Section~\ref{subsec:scales} show that $[(F_\alpha, \prescript{}{q}{B}_\alpha) \ ; \ \alpha \in [-m, \infty)]$, the extrapolated fractional power scales of order $m \in \N$ generated by $(F, \prescript{}{q}{B})$, are well-defined, such that we may concentrate on  
	\begin{equation*}
		 \prescript{}{\new{q}}{B}_{- 1} : D( \prescript{}{p}{B}_{- 1} ) = F_0 \subset  F_{-1}\rightarrow  F_{-1}.
	\end{equation*}
	By definition, it holds that $F_0 = F = L^q(\Omega)$. Observing that the dual operator $\prescript{}{q}{B}^\# = (\prescript{}{q}{B})'$ is given by 
	\begin{align*}
		\prescript{}{q}{B}^\#  &= \prescript{}{q'}{B} : D(\prescript{}{q'}{B}) = \{ u \in W^{2, q'}(\Omega) \colon \partial_{\bm{n}} u_{| \partial \Omega} = 0 \} \subset  L^{q'}(\Omega)
		\rightarrow L^{q'}(\Omega),\\ 
		\theta &\mapsto - \nabla \cdot \big( \kappa(\vphi_0) M(\vphi_0) \nabla \theta \big) + \mu_2 \theta, 
	\end{align*}
	we can apply the same arguments as above to obtain
	\begin{equation*}
		F_{-1} = (F^\#_1) =  \{ u \in W^{2, q'}(\Omega) \colon \partial_{\bm{n}} u_{| \partial \Omega} = 0 \} ' = W^{-2, q}_N(\Omega)
	\end{equation*}
	and hence, 
	\begin{equation*}
		 \prescript{}{\new{q}}{B}_{- 1} : D( \prescript{}{p}{B}_{- 1} ) = L^q(\Omega) \subset  W^{-2, q}_N(\Omega) \rightarrow  W^{-2, q}_N(\Omega).
	\end{equation*}
	Now the assertion follows exactly as in Lemma \ref{lemma:L_1}. 
\end{proof}

Having established invertibility of these operators, we proceed by showing that their inverses can be bounded independently of the time interval $[0, T]$ for all $T < T_0$, given a fixed $T_0$, cf. \cite[Lem.~7]{abels2021local}. 

\begin{lemma}\label{lem:L^-1_uniform_T}
	Under the assumptions of Lemma \ref{lemma:L_1} and for any $T_0 > 0$ there exists some constant $C_{\mathcal{L}^{-1}_{1, T}}(T_0) > 0$ such that \new{
		\begin{equation*}
			\norm{ \cL^{-1}_{1, T} (f)}_{Z_T^1} \leq C_{\mathcal{L}^{-1}_{1, T}}(T_0) \Big(   \norm{ f}_{Y_T^1} + \norm{\vphi(0)}_{  (W^{-2, p}_N(\Omega), W^{2, p}_N(\Omega))_{1- \frac{1}{r}, r}} \Big) 
		\end{equation*}
		and 
		\begin{equation*}
			\big\|  \cL^{-1}_{1, T} (f) -  \cL^{-1}_{1, T} (\tilde{f}) \big\|_{Z_T^1} \leq 
			C_{\mathcal{L}^{-1}_{1, T}}(T_0) \big\| f - \tilde{f} \big\|_{Y_T^1}
		\end{equation*}
		for all $T \in (0, T_0]$ and  $f, \tilde{f}  \in Y_T^1$. 
	}
\end{lemma}

\begin{proof}
	\new{
	We begin with the second equation and observe that for any $T > 0$ the difference $\cL^{-1}_{1, T} (f) -  \cL^{-1}_{1, T} (\tilde{f})$ solves the associated problem with homogeneous initial conditions  
	\begin{equation*} 
		\pt \vphi + 	\hat{\mathcal{L}_1} \vphi = f - \tilde{f}, \quad t \in (0, T), \quad \vphi(0) = 0, 
	\end{equation*}
	where the linear operator $\hat{\cL}_{1, T}$ is defined as in \eqref{eq:maxreg_inverse_homogeneous}. Since this problem has a unique solution, see Lemma~\ref{lemma:L_1}, we conclude that there exists some $C = C(T)$ --  the operator norm of $\hat{\cL}_{1, T}^{-1}$ -- such that 
	\begin{equation}\label{eq:estimate_L_depending_T}
			\big\|  \cL^{-1}_{1, T} (f) -  \cL^{-1}_{1, T} (\tilde{f}) \big\|_{Z_T^1} 
			= \norm{\hat{\cL}_{1, T}^{-1} (f - \tilde{f})}_{Z_T^1}
			\leq C(T) \big\| f - \tilde{f} \big\|_{Y_T^1}. 
	\end{equation}	
	Now we proceed with the estimates that are independent of $T$ and set for any $f \in L^r(0, T; W^{2, p}_N(\Omega))$
	\begin{equation*}
		\tilde{f}(t) \coloneqq 
		\begin{cases}
			f(t) & \textnormal{if } t \in [0, T], \\ 
			0 & \textnormal{if } t \in (T, T_0). 
		\end{cases}
	\end{equation*}
	Since the operator ${\mathcal{L}}_{1, T} : X_T^1 \to Y_T^1$ is invertible for every $0 < T < \infty$, we can  find the unique solutions $\vphi \in X^1_T$ and $\tilde{\vphi} \in X^1_{T_0}$ to 
	\eqref{eq:ch_strong} on the respective intervals. Observing that the two functions solve the same equation on $(0, T) \times \Omega$, uniqueness tells us that they need to coincide, i.e., $\tilde{\vphi}_{|(0, T) } = \vphi$. Thus, 
	\begin{align*}
		\norm{\mathcal{L}_{1, T}^{-1} (f)}_{Z^1_T}  
		= \norm{\vphi}_{Z^1_T} 
		&\leq \norm{\tilde{\vphi}}_{Z^1_{T_0}} 
		= \norm{\mathcal{L}_{1, T_0}^{-1} (\tilde{f})}_{Z^1_{T_0}}  \\ 
		&\leq C_{\cL_{1, T_0}^{–1}}(T) \Big(   \norm{ f}_{Y_T^1} + \norm{\vphi_0}_{  (W^{-2, p}_N(\Omega), W^{2, p}_N(\Omega))_{1- \frac{1}{r}, r}}
		\Big). 
	\end{align*}
	Finally, we note that the second equation can be derived analogously when combined with \eqref{eq:estimate_L_depending_T}. 
	}
\end{proof}

\begin{lemma}
	Under the assumptions of Lemma \ref{lemma:L_3}, and for any $T_0 > 0$ there exists some constant $C_{\mathcal{L}^{-1}_{3, T}}(T_0) > 0$ such that \new{
		\begin{equation*}
			\norm{ \cL^{-1}_{3, T} (f)}_{Z_T^3} \leq C_{\mathcal{L}^{-1}_{3 T}}(T_0) \Big(   \norm{ f}_{Y_T^3} +  \norm{\theta(0)}_{  (W^{-2, q}_N(\Omega), L^q(\Omega))_{1- \frac{1}{3r}, 3r}}\Big) 
		\end{equation*}
		and 
		\begin{equation*}
			\big\|  \cL^{-1}_{3, T} (f) -  \cL^{-1}_{3, T} (\tilde{f}) \big\|_{Z_T^3} \leq 
			C_{\mathcal{L}^{-1}_{3, T}}(T_0) \big\| f - \tilde{f} \big\|_{Y_T^3}
		\end{equation*}
		for all $T \in (0, T_0]$ and  $f, \tilde{f}  \in Y_T^3$. 
	}
\end{lemma}
\begin{proof}
	The assertion follows analogously to the proof of Lemma~\ref{lem:L^-1_uniform_T}. 
\end{proof}

Finally, it only remains to show these properties for $\mathcal{L}_2$. 

\begin{lemma}\label{lemma:L_2}
	Let $1 < r <\infty$. Suppose that assumptions \ref{A:iso} hold and let the initial value $\bm{u}_0  \in \X(\Omega)$ and the right hand side $\bm{f} \in L^{r}(0, T; \X(\Omega))$ be given along with some $\tilde{\vphi} \in W^{1, q^*}(\Omega)$, $q^* > n$. Then, for every $0 < T  < \infty$ there exists a unique $\bm{u} \in   W^{1, r} (0, T; \X(\Omega))$
	such that 
	\begin{alignat*}{2}
	  \pt \bm{u} + \mathcal{A}(\tilde{\vphi}) \bm{u}&= \bm{f}  \quad &&\textnormal{in } (0, T) \times \Omega, \numberthis \label{eq:displacement_strong}\\ 
		\bm{u}_{|t = 0} &= \bm{u}_0  \quad && \Omega. 
	\end{alignat*}
Moreover,
and for any $T_0 > 0$ there exists some constant $C_{\mathcal{L}^{-1}_{2, T}}(T_0) > 0$ such that \new{
	\begin{equation*}
		\norm{ \cL^{-1}_{2, T} (f)}_{Z_T^2} \leq C_{\mathcal{L}^{-1}_{2 T}}(T_0) \Big(   \norm{ f}_{Y_T^2} +  \norm{\bm{u}(0)}_{\X (\Omega)} \Big) 
	\end{equation*}
	and 
	\begin{equation*}
		\big\|  \cL^{-1}_{2, T} (f) -  \cL^{-1}_{2, T} (\tilde{f}) \big\|_{Z_T^2} \leq 
		C_{\mathcal{L}^{-1}_{2, T}}(T_0) \big\| f - \tilde{f} \big\|_{Y_T^2}
	\end{equation*}
	for all $T \in (0, T_0]$ and  $f, \tilde{f}  \in Y_T^2$. 
}
\end{lemma}

\begin{proof}
	Due to our assumptions \ref{A:iso}, it holds that
		\begin{align*}
			\mathcal{A}(\tilde{\vphi}) &: \X(\Omega) \rightarrow \X(\Omega), 
			\quad \bm{v} \mapsto \mathcal{B}^{-1}(\tilde{\vphi}) \mathcal{C}(\tilde{\vphi}) \bm{v}
		\end{align*}
		is a well-defined automorphism of $\X(\Omega)$ and therefore in particular closed. A well-know result, cf. \cite{Dore93}, shows that under these conditions $\mathcal{A}(\vphi_0)$ is of class $\mathcal{MR}_s$ on every bounded interval and all $s \in (1, \infty)$. In particular, $\mathcal{A}(\vphi_0) \in \mathcal{MR}$, cf. \cite{ARENDT20071}. As before, we can apply the standard argument, cf. \cite[Prop. 1.3]{ARENDT20071}, \cite[ Lemma 10]{haselboeck2024existence}, and deduce that the corresponding Cauchy-problem on bounded intervals with non-trivial initial data $\bm{u}_0 \in \X (\Omega)$ is also well-posed. 
\\ 
Finally, the last assertion follows as in Lemma~\ref{lem:L^-1_uniform_T}. 
\end{proof}

\medspace

\section{Lipschitz continuity of $\mathcal{F}$: visco-elastic case} \label{sec:contraction}

Before we can start to derive the Lipschitz estimates, some preparations are in order. Firstly, we exploit the general interpolation results we recalled in Section~\ref{sec:preliminaries} to find bounds for $\vphi$ and $\theta$ in suitable spaces. 
From the definition of $X^1_T$ and \eqref{int:BUC}, it follows that 
\begin{align}\label{embedding_X^1_1}
	X^1_T &\hookrightarrow BUC ([0, T];  (W^{-2, p}_N(\Omega), W^{2, p}_N(\Omega))_{1- \frac{1}{r}, r})  \hookrightarrow  BUC ([0, T]; B^{2- \frac{4}{r}}_{p, r} (\Omega)). 
\end{align}
Since $q > n $, we can find some $r > \frac{4q}{q-n} = \frac{8p}{2p-n} > 1$. Hence, there exists a $\delta > 0$ such that $2p-n - 4 p \delta > 0 $ and $r > \frac{8p}{2p-n - 4 p \delta}$. For this $\delta$, it follows that $2-\frac{4}{r} - \frac{n}{p} >  1 + 2 \delta- \frac{n}{2p}$ and we obtain the continuous embedding $B^{2- \tfrac{4}{r}}_{p, r} (\Omega)  \hookrightarrow W^{1+ 2 \delta, q}(\Omega)$. Moreover, since $q > n$, we further get $W^{1, q}(\Omega) \hookrightarrow C^\beta(\overline{\Omega})$ for some $\beta > 0$ and hence, 
\begin{equation}\label{embedding_X^1}
	X^1_T \hookrightarrow BUC ([0, T]; W^{1+ 2 \delta, q} (\Omega)) \hookrightarrow BUC ([0, T]; W^{1, q} (\Omega)) \hookrightarrow BUC ([0, T]; C^\beta (\overline{\Omega})). 
\end{equation}
Furthermore, the choice $r > \frac{8p}{2p-n} $ and standard interpolation results give rise to the inequality 
\begin{equation*}
	\norm{\vphi}_{W^{1, q}} \leq C \norm{ \vphi}_{B^{2 - \frac{4}{r}}_{p, r} }^{1 - \vartheta} \norm{\vphi}_{W^{-2, p}}^\vartheta 
\end{equation*}
for some $\vartheta \in (0, 1)$. Along with the embedding $W^{1, r}(0, T; W^{-2, p}(\Omega)) \hookrightarrow C^{1 - \frac{1}{r}}([0; T]; W^{-2, p}(\Omega))$ and \eqref{embedding_X^1_1}, we find Lemma~\ref{lemma:interpolation_hoelder} applicable, deducing
\begin{equation}\label{embedding:hoelder}
	X^1_T \hookrightarrow C^{(1 - \frac{1}{r}) \vartheta}([0, T]; W^{1, q}(\Omega)) \coloneqq C^{ \tilde{\vartheta}}([0, T]; W^{1, q}(\Omega)) 
\end{equation}
for $\tilde{\vartheta} \coloneqq (1 - \frac{1}{r}) \vartheta > 0$. \\ 
On the other hand, we analogously find that 
\begin{align*}
	X^3_T \hookrightarrow BUC ([0, T];  (W^{-2, q}_N(\Omega), L^q(\Omega))_{1- \frac{1}{3r}, 3r})
	&\hookrightarrow BUC ([0, T];  W^{-\delta, q}(\Omega))\\ 
	&\hookrightarrow BUC ([0, T];  W^{-1, q}(\Omega)), 
	\numberthis
	\label{embedding_X^3}
\end{align*}
where 
$(W^{-2, q}_N(\Omega), L^q(\Omega))_{1- \frac{1}{3r}, 3r} = B^{- \frac{2}{3r}}_{q, 3r}  \hookrightarrow W^{-\delta, q}(\Omega) \hookrightarrow W^{-1, q}(\Omega)$, if $\delta < 1$, since $- \frac{2}{3r} - \frac{n}{q} > - \new{\delta} - \frac{n}{q}$ for any $r > \frac{2}{3\delta}$. 
\par 
\medskip 
The following lemma provides a general estimate, which we will employ multiple times in the proof of Proposition \ref{prop:contraction}.  

\begin{lemma}\label{lemma:general_estimate}
	Let $f_i \in L^p(\Omega)$ and $\vphi_i \in W^{1, q}(\Omega)$ with $\norm{\vphi_i}_{W^{1, p}}  \leq R $ for some $R > 0$ and $i = 1, 2$. Then, it holds that 
	\begin{align*}
		\norm{\nabla \cdot \big(m(\vphi_1) \nabla f_1 \big) -  \nabla \cdot \big(m(\vphi_2) \nabla f_2\big) }_{W^{-2, p}_N}  
		 \leq C (R)  \Big(  \norm{  f_1 - f_2 }_{L^p} + \norm{ f_2}_{L^p} \norm{\vphi_1 - \vphi_2 }_{W^{1, q}}  \Big)
	\end{align*}
\end{lemma}

\begin{proof}
	Exploiting the multiplication result $W^{-1, p}(\Omega) \times W^{1, q}(\Omega) \hookrightarrow W^{-1, p}(\Omega)$,  cf. \eqref{mult:-1p}, along with the fact that $\norm{\nabla \cdot (\cdot)}_{W^{-2, p}_ N }  \leq C \norm{\cdot}_{W^{-1, p}}$, we compute 
	\begin{align*}
		&\norm{\nabla \cdot \big(m(\vphi_1) \nabla f_1 \big) -  \nabla \cdot \big(m(\vphi_2) \nabla f_2\big) }_{W^{-2, p}_N} 
		= \norm{\nabla \cdot \Big(m(\vphi_1) \nabla \big( f_1 - f_2\big)  + \big(m(\vphi_1 - m(\vphi_2) \big) \nabla f_2) \Big) }_{W^{-2, p}_N}\\ 
		& \quad 
		\leq  C \norm{ m(\vphi_1) \nabla \big( f_1 - f_2\big) }_{W^{-1, p}}   + C  \norm{ \big(m(\vphi_1 - m(\vphi_2) \big) \nabla f_2) }_{W^{-1, p}} \\ 
		& \quad 
		\leq C \norm{ m(\vphi_1) }_{W^{1, q}} \norm{  \nabla \big( f_1 - f_2\big) }_{W^{-1, p}}
				+ C \norm{\nabla f_2}_{W^{-1, p}} \norm{m(\vphi_1) - m(\vphi_2) }_{W^{1, q}} . 
	\end{align*}
	Taking advantage of $m \in C^2(\R)$ along with the fact that $1 - n/ q  > 0$ due to the assumptions  on $q$, an application of Theorem~\ref{thm:compositon} yields 
	\begin{equation*}
		\norm{m(\vphi_1)}_{W^{1, q}} \leq C(R) \quad \textrm{and} \quad \norm{m(\vphi_1) - m(\vphi_2)}_{W^{1, q}} \leq L(R) \norm{\vphi_1  -  \vphi_2}_{W^{1, q}}
	\end{equation*} 
	for all $\vphi_i \in W^{1, q}(\Omega)$ with $\norm{\vphi_i}_{W^{1, q}} \leq R$. 
	Moreover, we note that $ \norm{  \nabla f_1  }_{W^{-1, p}} \leq C  \norm{ f_1 }_{L^p}$. 
	Applying these inequalities to the estimate above leads to 
	\begin{equation*}
		\norm{\nabla \cdot \big(m(\vphi_1) \nabla f_1 \big) -  \nabla \cdot \big(m(\vphi_2) \nabla f_2\big) }_{W^{-2, p}_N}  
		\leq C (R)  \Big(  \norm{  f_1 - f_2 }_{L^p}
		+ \norm{ f_2}_{L^p} \norm{\vphi_1 - \vphi_2 }_{W^{1, q}}  \Big), 
	\end{equation*}
	which concludes the proof. 
\end{proof}

We proceed with the proof of Proposition \ref{prop:contraction}. \new{Note the the estimates for $\F^1_T$ and $\F^3_T$ are often very similar to those in Section~\ref{sec:contraction_non_ve} and that their numbering is often identical (apart from the fact that $\F^3_T$ in this section corresponds to $\F^2_T$ in preceding part). However, due to the weaker spaces used here, some estimates or a little bit more subtle and require the multiplication results we recalled in the preliminaries.}

\begin{proof}[Proof of Proposition \ref{prop:contraction}]
	Let $(\vphi_i, \bm{u}_i,  \theta_i) \in X_T$ such that $\norm{(\vphi_i, \bm{u}_i,  \theta_i)}_{X_T} < R$ for $i = 1, 2$. In order to establish Lipschitz continuity for $\F_T$, we study the components $\F_T^j$, $j = 1, 2, 3$, separately. 
	\par 
	\medskip 
	\begin{subequations}
	\textit{Ad $\F_T^1$:} First of all, we find
	\begin{align*}
		&\norm{\F_T^1 (\vphi_1, \bm{u}_1, \theta_1) - \F_T^1 (\vphi_2, \bm{u}_2, \theta_2)}_{Y^1_T} \\ 
		&\quad \leq   \Big\| \eps\, \Delta  \big( m(\vphi_0)  \Delta (\vphi_1 - \vphi_2 ) \big)
		 - \eps \nabla \cdot \big( m(\vphi_1) \nabla \Delta \vphi_1  - m(\vphi_2) \nabla \Delta \vphi_2  \big) \Big\| _{Y^1_T}  \\ 
		& \quad \quad + \Big\| \nabla \cdot \big(
			m(\vphi_1) \nabla \big(\frac{1}{\eps} \psi'(\vphi_1)  \big)  
			- m(\vphi_2) \nabla \big( \frac{1}{\eps} \psi'(\vphi_2)\big) \big)\Big\|_{Y^1_T}\\ 
		& \quad \quad + \Big\| \nabla \cdot \big(
		m(\vphi_1) \nabla \Wp(\vphi_1, \E(\bm{u}_1))
		- m(\vphi_2) \nabla \Wp(\vphi_2, \E(\bm{u}_2)) \big)\Big\|_{Y^1_T}\\ 
		& \quad \quad + \begin{aligned}[t]
			 \Big\| \nabla \cdot \big( &m(\vphi_1) \nabla \big( M(\vphi_1) (\theta_1 - \alpha(\vphi_1) \nabla \cdot \bm{u}_1 ) \alpha'(\vphi_1) \nabla \cdot \bm{u}_1 \big) \\
			 - &m(\vphi_2) \nabla \big( M(\vphi_2) (\theta_2 - \alpha(\vphi_2) \nabla \cdot \bm{u}_2 ) \alpha'(\vphi_2) \nabla \cdot \bm{u}_2 \big) \big) \Big\|_{Y^1_T}
		\end{aligned}
		\\ 
		& \quad \quad + \Big\| \nabla \cdot \big(
		m(\vphi_1) \nabla \big( \frac{M'(\vphi_1)}{2} (\theta_1 - \alpha(\vphi_1) \nabla \cdot \bm{u}_1 )^2 \big)  
		- m(\vphi_2) \nabla \big( \frac{M'(\vphi_2)}{2} (\theta_2 - \alpha(\vphi_2) \nabla \cdot \bm{u}_2 )^2 
		\big) \big) \Big\|_{Y^1_T} \\ 
		& \quad \quad + \norm{S_s(\vphi_1, \bm{u}_1,  \theta_1) - S_s(\vphi_2, \bm{u}_2,  \theta_2)}_{Y^1_T}\\ 
		& \quad \eqqcolon I + II + III+ IV + V +VI . 
	\end{align*}
	\textit{Ad $\F_T^1.I$:} 
	The triangle inequality and the multiplication results \eqref{mult:-1p} and \eqref{eq:mutl_delta_p} imply 
	\begin{align*}
		 &\big\|  \Delta  \big( m(\vphi_0)  \Delta (\vphi_1 - \vphi_2 ) \big)
		 - \nabla \cdot \big( m(\vphi_1) \nabla \Delta \vphi_1  - m(\vphi_2) \nabla \Delta \vphi_2 \big) \big\| _{W^{-2, p}_N}\\ 
		 & \quad \leq\begin{aligned}[t]
		 	&\big\| \nabla \cdot   \big( m(\vphi_0) \nabla \Delta (\vphi_1 - \vphi_2 ) \big)
		 	-  \nabla \cdot \big( m(\vphi_1) \nabla \Delta \vphi_1  - m(\vphi_2) \nabla \Delta \vphi_2 \big) \big\| _{W^{-2, p}_N}\\ 
		 	&+ \big\|  \nabla \cdot \big( m'(\vphi_0) \nabla \vphi_0 \, \Delta (\vphi_1 - \vphi_2) \big) \big\|_{W^{-2, p}}
		 \end{aligned}
		 \\ 
		 & \quad \leq \begin{aligned}[t]
		 	& \big\| \big( m(\vphi_0) - m(\vphi_1) \big) \big( \nabla \Delta \vphi_1 - \nabla \Delta \vphi_2 \big)  \big\| _{W^{-1, p}}
		 	+  \big\| \big(m(\vphi_1) - m(\vphi_2) \big) \nabla \Delta \vphi_2   \big\| _{W^{-1, p}}\\ 
		 	&+\norm{ m'(\vphi_0) \nabla \vphi_0 \, \Delta(\vphi_1 - \vphi_2) }_{W^{-1, p}}
		 \end{aligned}\\ 
		 & \quad  \leq \begin{aligned}[t]
		 	&\norm{m(\vphi_0) - m(\vphi_1)}_{W^{1, q}} \norm{\nabla \Delta \vphi_1 - \nabla \Delta \vphi_2}_{W^{-1, p}}
		 	+ \norm{m(\vphi_1) - m(\vphi_2)}_{W^{1, q}} \norm{\nabla \Delta \vphi_2 }_{W^{-1, p}} \\ 
		 	&+ \norm{m'(\vphi_0)}_{W^{1, q}} \norm{\nabla \vphi_0}_{W^{1, p}} \norm{\Delta (\vphi_1 - \vphi_2)}_{W^{\tilde{\delta}, p}}
		 \end{aligned}\\ 
		 & \quad  \leq C(R) \norm{\vphi_0 - \vphi_1}_{W^{1, q}} \norm{\vphi_1 -  \vphi_2}_{W^{2, p}}
		 + C(R) \norm{\vphi_1 - \vphi_2}_{W^{1, q}} \norm{ \vphi_2 }_{W^{2, p}}
		 + C(R) \norm{\vphi_1 - \vphi_2}_{W^{2-\tilde{\delta}, p}}. 
	\end{align*}
	Recall the interpolation estimate 
		\begin{equation*}
			\norm{f}_{W^{2-\tilde{\delta}, p}} \leq M \norm{f}_{W^{-2, p}}^{1-\hat{\vartheta}}  \norm{f}_{W^{2, p}}^{\hat{\vartheta}}
		\end{equation*}
		for some $M > 0$ and $\hat{\vartheta} \in (0, 1)$ and remember that $\vphi_i \in L^r(0, T; W^{2, p} \cap L^\infty(0, T; W^{-2,p} )$. Thus, standard interpolation results for Bochner spaces imply 
		\begin{equation*}
			\norm{\vphi_i }_{L^{r^*}(W^{2- \tilde{\delta}, p})} \leq M \norm{\vphi_i}_{L^\infty((W^{-1, p})}^{1-\hat{\vartheta}}  \norm{\vphi_i}_{L^r(W^{2, p})}^{\hat{\vartheta}}, 
		\end{equation*}
		where $\frac{1}{r^*} = \frac{\hat{\vartheta}}{r}$. In particular, $ \vphi_1 - \vphi_2 \in {L^{r^*}(0, T; W^{2- \tilde{\delta}, q})}$ for some $r^* > r$ and the norm is bounded by the respective norms in $X_T^1$. 
		\new{Using the definition of Hölder spaces, a quick computation shows that
			\begin{equation*}
				\norm{\vphi_1(0) - \vphi}_{C^0{W^{1, p}}} \leq T^{\tilde{\vartheta}} \norm{ \vphi_1}_{C^{\tilde{\vartheta}}(W^{1, p})}. 
			\end{equation*}
	Moreover, the analogous relation holds for $\vphi_1 - \vphi_2$, and we note that since $\vphi_i \in X^1_T$, their initial values coincide. In particular, it holds that $\vphi_1(0) - \vphi_2(0) =0$.}
	Using the embedding $BUC ([0, T], W^{1, q}(\Omega)) \xhookrightarrow{} C^0([0, T]; W^{1, q} (\Omega))$, where the latter space denotes continuous functions on $[0, T]$ with values in $W^{1, q}(\Omega)$, along with \eqref{embedding:hoelder}, this estimate leads to 
	\begin{align*}
		 I &\leq C(R) \begin{aligned}[t]
		 	 \big(  &\norm{ \vphi_1(0) - \vphi_1}_{C^0(W^{1, q})} \norm{\vphi_1 -  \vphi_2}_{L^r(W^{2, p})}
		 	+ \norm{\vphi_1 - \vphi_2}_{C^0(W^{1, q})}  \norm{ \vphi_2}_{L^r(W^{2, p})} \big) \\ 
		 	&+ \norm{\vphi_1 - \vphi_2}_{L^r(W^{2- \tilde{\delta}, p})} 
		 \end{aligned} \\ 
		&  \leq C (R) \begin{aligned}[t]
			\Big( &\norm{ \vphi_1(0) - \vphi_1}_{C^0(W^{1, q})} \norm{\vphi_1 -  \vphi_2}_{L^r(W^{2, p})}
			+  \norm{(\vphi_1  - \vphi_2)  - (\vphi_1 (0) - \vphi_2(0))}_{C^0(W^{1, q})}  \\
			&+ T^{\frac{1}{(r^*)'}}\norm{\vphi_1 - \vphi_2}_{L^{r^*}(W^{2- \tilde{\delta}, p})}
			 \Big)
		\end{aligned} \\ 
		&  \leq \begin{aligned}[t]
			C (R)& T^{\tilde{\vartheta}}   \Big(   \norm{ \vphi_1}_{C^{\tilde{\vartheta}}(W^{1, q})} \norm{\vphi_1 -  \vphi_2}_{L^r(W^{2, p})}
			+  \norm{\vphi_1 - \vphi_2}_{C^{\tilde{\vartheta}}(W^{1, q})}    \Big) \\ 
			&+ C(R)T^{\frac{1}{(r^*)'}}\norm{\vphi_1 - \vphi_2}_{L^{r^*}(W^{2- \tilde{\delta}, p})}.
		\end{aligned} 
		 \numberthis \label{est:F_1_I}
	\end{align*} 
	\\ 
	\textit{Ad $\F_T^1.II$:} Recall that the functions $\vphi_i$, $i = 1, 2$, are bounded in $X^1_T$ and that $X^1_T \hookrightarrow BUC ([0, T]; C^\beta (\overline{\Omega}))$. In particular, $\norm{\vphi_i}_{C(\OT)}  < C(R)$ which, along with $\psi' \in C^1(\R)$, implies  
	\begin{equation*}
		\abs{	\psi'(\vphi_i)} < C  \quad \textnormal{and} \quad \abs{\psi'(\vphi_1) - \psi'(\vphi_2)} \leq  C \abs{\vphi_1 - \vphi_2}. 
	\end{equation*}
	Moreover, we note that $\norm{\vphi_i}_{W^{1, q}} \leq C(R)$ due to the embedding $X^1_T \hookrightarrow BUC ([0, T]; W^{1, q}(\Omega))$. Applying Lemma~\ref{lemma:general_estimate} therefore yields
	\begin{equation*}
		 \Big\| \nabla \cdot \big(
		m(\vphi_1) \nabla \big(\frac{1}{\eps} \psi'(\vphi_1)  \big)  
		- m(\vphi_2) \nabla \big( \frac{1}{\eps} \psi'(\vphi_2)\big) \big)\Big\|_{W^{-2, p}_N}
		\leq C (R) \norm{\vphi_1 - \vphi_2}_{W^{1, q}} , 
	\end{equation*}
	leading to 
	\begin{equation}\label{est:F_1_III}
		II \leq \norm{ C (R) \norm{\vphi_1 - \vphi_2}_{W^{1, q}} }_{L^r} \leq C(R) T^{\frac{1}{r}} \norm{\vphi_1 - \vphi_2}_{BUC (W^{1, q}))}. 
	\end{equation}
	\\ 
	\textit{Ad $\F_T^1.III$:} Due to Remark \ref{A:W_growth_conditions} we know that $\abs{W_{,\vphi}(\vphi', \E') } \leq C_2 \left( \abs{\E'}^2 + \abs{\vphi'}^2 +1 \right)$ and since $\vphi_i, \bm{u}_i$ are bounded in $X^j_T$, $j = 1, 2$, respectively, Lemma~\ref{lemma:general_estimate} is applicable, yielding 
	\begin{align*}
		&\Big\| \nabla \cdot \big(
		m(\vphi_1) \nabla \Wp(\vphi_1, \E(\bm{u}_1))
		- m(\vphi_2) \nabla \Wp(\vphi_2, \E(\bm{u}_2)) \big)\Big\|_{W^{-2, p}_{N}} \\ 
		& \quad \leq C(R) \big( \norm{ \Wp(\vphi_1, \E(\bm{u}_1)) - \Wp(\vphi_2, \E(\bm{u}_2))}_{L^p} 
		+ \norm{\Wp(\vphi_2, \E(\bm{u}_2))}_{L^p}  \norm{\vphi_1 - \vphi_2}_{W^{1, q}}   \big), 
	\end{align*}
	where 
	\begin{equation*}
		\norm{ W_{,\vphi}(\vphi_2, \E(\bm{u}_2)) }_{L^p} \leq  C \big( \norm{\vphi_2}_{L^q}^2 +  \norm{\bm{u}_2}_{\X}^2  +1 \big)  \leq C(R) 
	\end{equation*}
	due to $\vphi_2 \in X^1_T \hookrightarrow BUC([0, T]; W^{1, q}(\Omega))$ and $\bm{u}_2 \in X^2_T \hookrightarrow  C^0([0, T]; \X(\Omega))$. Taking advantage of the explicit form of the elastic energy density $W$ postulated in \ref{A:W}, along with the computations from \cite[Lem. 17, p.17]{haselboeck2024existence}, we find 
		\begin{align*}
		&| \Wp(\vphi_1 , \E(\bm{u}_1)) - \Wp(\vphi_2, \E(\bm{u}_2))|  \\ 
		& \quad \leq   C \Big( \abs{ \E(\bm{u}_1) }  +  \abs{\E(\bm{u}_2)} + \abs{\E(\bm{u}_2)}^2 + 1 \Big) \abs{\vphi_1 - \vphi_2}
		+ C \Big( \abs{ \E(\bm{u}_1) }  +  \abs{\E(\bm{u}_2)}  + 1\Big)  \abs{\E(\bm{u}_1) - \E(\bm{u}_2) }. 
	\end{align*}
	Thus, 
	\begin{align*}
		 &\norm{ \Wp(\vphi_1, \E(\bm{u}_1)) - \Wp(\vphi_2, \E(\bm{u}_2))}_{L^p}  \\ 
		 & \  \leq C \norm{\Big( \abs{ \E(\bm{u}_1) }  +  \abs{\E(\bm{u}_2)} + \abs{\E(\bm{u}_2)}^2 + 1 \Big) \abs{\vphi_1 - \vphi_2} }_{L^p}
		 + C \norm{ \Big( \abs{ \E(\bm{u}_1) }  +  \abs{\E(\bm{u}_2)}  + 1\Big)  \abs{\E(\bm{u}_1) - \E(\bm{u}_2) } }_{L^p} \\ 
		 & \ \leq C(R) \Big( \norm{ \bm{u}_1}_{\X} + \norm{\bm{u}_2}^2_{\X} + 1 \Big) \norm{\vphi_1 - \vphi_2}_{L^\infty }
		 + C(R) \Big( \norm{\bm{u}_1}_{\X} + \norm{\bm{u}_2}_{\X} + 1 \Big) \norm{\bm{u}_1 - \bm{u}_2}_{\X} \\ 
		 & \ \leq C(R) \big(   \norm{\vphi_1 - \vphi_2}_{W^{1, q}} +   \norm{\bm{u}_1 - \bm{u}_2}_{\X}\big). 
	\end{align*}
	In summary, this leads to 
	\begin{align*}
		III = \Big\| \nabla \cdot \big(
		m(\vphi_1) \nabla \Wp(\vphi_1, \E(\bm{u}_1))
		&- m(\vphi_2) \nabla \Wp(\vphi_2, \E(\bm{u}_2)) \big)\Big\|_{L^r(W^{-2, p}_{N})} \\ 
		 &  \leq C(R) T^{\frac{1}{r}} \big(   \norm{\vphi_1 - \vphi_2}_{C^0(W^{1, q})} +   \norm{\bm{u}_1 - \bm{u}_2}_{C^0(\X) }\big). 
		 \numberthis \label{F_1_IV}
	\end{align*}
	\\ 
	\textit{Ad $\F_T^1.IV$:} We apply Lemma \ref{lemma:general_estimate} to $(IV)$, arriving at 
	\begin{align*}
		&\begin{aligned}
			 \Big\|\nabla \cdot \Big( &m(\vphi_1) \nabla \big( M(\vphi_1) \alpha'(\vphi_1)  \theta_1  \nabla \cdot \bm{u}_1-  M(\vphi_1)\alpha(\vphi_1) \alpha'(\vphi_1) (\nabla \cdot \bm{u}_1)^2  \big) \\
			- &m(\vphi_2) \nabla \big( M(\vphi_2) \alpha'(\vphi_2)  \theta_2  \nabla \cdot \bm{u}_2-  M(\vphi_2)\alpha(\vphi_2) \alpha'(\vphi_2) (\nabla \cdot \bm{u}_2)^2  \big) \Big) \Big\|_{W^{-2, p}_{N}}
		\end{aligned}
		\\ & \quad \leq 
		 \begin{aligned}[t]
			C(R) \Big( &\norm{M(\vphi_1) \alpha'(\vphi_1)  \theta_1  \nabla \cdot \bm{u}_1
				- M(\vphi_2) \alpha'(\vphi_2)  \theta_2  \nabla \cdot \bm{u}_2 }_{L^p}\\ 
			+ &\norm{
				M(\vphi_1)\alpha(\vphi_1) \alpha'(\vphi_1) (\nabla \cdot \bm{u}_1)^2 
				-  M(\vphi_2)\alpha(\vphi_2) \alpha'(\vphi_2) (\nabla \cdot \bm{u}_2)^2 }_{L^p} \\ 
			+ & \norm{ M(\vphi_2) \alpha'(\vphi_2)  \theta_2  \nabla \cdot \bm{u}_2-  M(\vphi_2)\alpha(\vphi_2) \alpha'(\vphi_2) (\nabla \cdot \bm{u}_2)^2 }_{L^p} \norm{\vphi_1 - \vphi_2}_{W{1,q}} 
				\Big). 
		\end{aligned}
	\end{align*}
	Let $\tilde{M}(\vphi_1) \coloneqq M(\vphi_i) \alpha'(\vphi_1)$ and $\hat{M}(\vphi_i) \coloneqq M(\vphi_i)\alpha(\vphi_i) \alpha'(\vphi_i)$ for $i = 1, 2$. We recall the embedding $\vphi_i \in  X^1_T \hookrightarrow  BUC([0, T]; C^\beta(\overline{\Omega}))$ and that the functions $M, \alpha, \alpha'$ are continuously differentiable, which leads to
	\begin{gather*}
		\abs{\tilde{M}(\vphi_1) - \tilde{M}(\vphi_2)} \leq C(R) \abs{\vphi_1 - \vphi_2},\\ 
		\abs{\hat{M}(\vphi_1) - \hat{M}(\vphi_2)} \leq C(R) \abs{\vphi_1 - \vphi_2}.
	\end{gather*}
	Using this, we estimate 
	\begin{align*}
		&\norm{ \tilde{M}(\vphi_1)  \theta_1  \nabla \cdot \bm{u}_1 - \tilde{M}(\vphi_2) \theta_2  \nabla \cdot \bm{u}_2  }_{L^p}\\ 
			& \quad \leq \norm{\tilde{M}(\vphi_1)  \theta_1 ( \nabla \cdot \bm{u}_1 -  \nabla \cdot \bm{u}_2)   }_{L^p}
				+ \norm{\tilde{M}(\vphi_1)  \nabla \cdot \bm{u}_2 ( \theta_1 - \theta_2 ) }_{L^p}
				+ \norm{
					 \theta_2  \nabla \cdot \bm{u}_2   (\tilde{M} (\vphi_1) - \tilde{M}(\vphi_2))}_{L^p} \\ 
			& \quad \leq 
			\begin{aligned}[t]
				|| \tilde{M}(\vphi_1) ||_{L^\infty} \norm{\theta_1}_{L^q} \norm{ \bm{u}_1 -  \bm{u}_2}_{\bm{W}^{1, q}} 
				&+ ||\tilde{M}(\vphi_1)||_{L^\infty} \norm{ \bm{u}_2}_{\bm{W}^{1, q}}  \norm{\theta_1 - \theta_2}_{L^q}\\
				&+ \norm{\theta_2}_{L^q} \norm{\bm{u}_2}_{\bm{W}^{1, q}} \norm{\vphi_1- \vphi_2}_{W^{1, q}}
			\end{aligned}\\ 
			& \quad \leq C(R) \Big(  \norm{\theta_1}_{L^q} \norm{ \bm{u}_1 -  \bm{u}_2}_{\X}  + \norm{\theta_1 - \theta_2}_{L^q} + \norm{\theta_2}_{L^q} \norm{\vphi_1- \vphi_2}_{W^{1, q}}\Big). 
			\numberthis \label{est:F_1_help_1}
	\end{align*}
	A similar computation yields 
	\begin{align*}
		\norm{
			\hat{M}(\vphi_1) (\nabla \cdot \bm{u}_1)^2 
			-  \hat{M}(\vphi_1) (\nabla \cdot \bm{u}_2)^2 }_{L^p} 
			\leq C(R) \Big( \norm{ \bm{u}_1 -  \bm{u}_2}_{\X} +  \norm{\vphi_1- \vphi_2}_{W^{1, q}} \Big)
			\numberthis \label{est:F_1_help_2}
	\end{align*} 
	and due to our embedding theorems, we deduce
	\begin{equation*}
		\norm{ M(\vphi_2) \alpha'(\vphi_2)  \theta_2  \nabla \cdot \bm{u}_2-  M(\vphi_2)\alpha(\vphi_2) \alpha'(\vphi_2) (\nabla \cdot \bm{u}_2)^2 }_{L^p}
		\leq C(R) \Big( \norm{\theta_2}_{L^q} + 1\Big). 
	\end{equation*}
	Therefore, 
	\begin{align*}
		&\begin{aligned}
			\Big\|\nabla \cdot \big( &m(\vphi_1) \nabla \big( M(\vphi_1) \alpha'(\vphi_1)  \theta_1  \nabla \cdot \bm{u}_1-  M(\vphi_1)\alpha(\vphi_1) \alpha'(\vphi_1) (\nabla \cdot \bm{u}_1)^2  \big) \\
			- &m(\vphi_2) \nabla \big( M(\vphi_2) \alpha'(\vphi_2)  \theta_2  \nabla \cdot \bm{u}_2-  M(\vphi_2)\alpha(\vphi_2) \alpha'(\vphi_2) (\nabla \cdot \bm{u}_2)^2  \big)  \Big\|_{W^{-2, p}_{N}}
		\end{aligned} \\ 
		&\quad  \leq C (R) \Big(  \big(\norm{\theta_1}_{L^q} + 1\big)  \norm{ \bm{u}_1 -  \bm{u}_2}_{\X} + \norm{\theta_1 - \theta_2}_{L^q} +  \big(\norm{\theta_2}_{L^q} + 1\big)  \norm{\vphi_1- \vphi_2}_{W^{1, q}} \Big). 
	\end{align*}
	and hence, 
	\begin{align*}
		IV &\leq C(R)\begin{aligned}[t] \Big(
				( &T^{\frac{2}{3r}}\norm{\theta_1}_{L^{3r}(L^q)} + T^\frac{1}{r})\norm{ \bm{u}_1 -  \bm{u}_2}_{C^0(\X) } 
				+ T^{\frac{2}{3r}} \norm{\theta_1 - \theta_2}_{L^{3r}(L^q)}   \\ 
				&+ ( T^{\frac{2}{3r}}\norm{\theta_2}_{L^{3r}(L^q)} + T^\frac{1}{r})  \norm{\vphi_1- \vphi_2}_{C^0(W^{1, q})} \Big)
				\end{aligned}\\ 
			&\leq C(R, T) \Big( \norm{ \bm{u}_1 -  \bm{u}_2}_{C^0(\X) } 
			+  \norm{\theta_1 - \theta_2}_{L^{3r}(L^q)}   
			+   \norm{\vphi_1- \vphi_2}_{C^0(W^{1, q})} \Big), 
			\numberthis \label{est:F_1_V}  
	\end{align*}
	with $C(R, T) \rightarrow 0$ as $T \rightarrow 0$. 
	\par 
	\medskip
	\textit{Ad $\F_T^1.V$:} We observe 
	\begin{align*}
		&M'(\vphi_1) (\theta_1 - \alpha(\vphi_1) \nabla \cdot \bm{u}_1 )^2 
		- M'(\vphi_2) (\theta_2 - \alpha(\vphi_2) \nabla \cdot \bm{u}_2 )^2\\ 
		& \quad = \begin{aligned}[t]
			\Big( M'(\vphi_1) \theta_1^2 - M'(\vphi_2) \theta_2^2 \Big) 
			&+ 2 \Big( M'(\vphi_2) \alpha(\vphi_2) \theta_2 \nabla \cdot \bm{u}_2 -  M'(\vphi_1) \alpha(\vphi_1)\theta_1\nabla \cdot \bm{u}_1\Big)\\ 
			&+ \Big( M(\vphi_1) \alpha(\vphi_1)^2 (\nabla \cdot \bm{u}_1)^2 -  M(\vphi_2) \alpha(\vphi_2)^2 (\nabla \cdot \bm{u}_2)^2  \Big)
		\end{aligned}
	\end{align*}
	and calculate for the first term
	\begin{align*}
		\norm{M'(\vphi_1) \theta_1^2 - M'(\vphi_2) \theta_2^2 }_{L^p} 
		\leq C(R) \big( \norm{\theta_1}_{L^q} +\norm{\theta_2}_{L^q} \big) \norm{\theta_1 - \theta_2}_{L^q} + C(R) \norm{\theta_2}_{L^q}^2 \norm{\vphi_1 - \vphi_2}_{W^{1, q}}. 
	\end{align*}
	Note that the two remaining terms can be estimated analogously to \eqref{est:F_1_help_1} and \eqref{est:F_1_help_2}, respectively, 
	\new{yielding 
	\begin{align*}
		&\norm{M'(\vphi_2) \alpha(\vphi_2) \theta_2 \nabla \cdot \bm{u}_2 -  M'(\vphi_1) \alpha(\vphi_1)\theta_1\nabla \cdot \bm{u}_1 }_{L^p}\\ 
		& \quad \leq
		C(R) \Big(  \norm{\theta_1}_{L^q} \norm{ \bm{u}_1 -  \bm{u}_2}_{\X}  + \norm{\theta_1 - \theta_2}_{L^q} + \norm{\theta_2}_{L^q} \norm{\vphi_1- \vphi_2}_{W^{1, q}}\Big)
	\end{align*}
	and 
	\begin{align*}
		\norm{M(\vphi_1) \alpha(\vphi_1)^2 (\nabla \cdot \bm{u}_1)^2 -  M(\vphi_2) \alpha(\vphi_2)^2 (\nabla \cdot \bm{u}_2)^2 }_{L^p} 
		\leq C(R) \Big( \norm{ \bm{u}_1 -  \bm{u}_2}_{\X} +  \norm{\vphi_1- \vphi_2}_{W^{1, q}} \Big). 
	\end{align*}
	}
	Once again, we apply Lemma~\ref{lemma:general_estimate} and \new{find after inserting these estimates that }
	\begin{align*}
		&\Big\| \nabla \cdot \big(
		m(\vphi_1) \nabla \big( \frac{M'(\vphi_1)}{2} (\theta_1 - \alpha(\vphi_1) \nabla \cdot \bm{u}_1 )^2 \big)  
		- m(\vphi_2) \nabla \big( \frac{M'(\vphi_2)}{2} (\theta_2 - \alpha(\vphi_2) \nabla \cdot \bm{u}_2 )^2 
		\big) \big) \Big\|_{W^{-2, p}_{N}}\\
		& \quad \leq \begin{aligned}[t]
			&C(R) \norm{M'(\vphi_1) (\theta_1 - \alpha(\vphi_1) \nabla \cdot \bm{u}_1 )^2 
				- M'(\vphi_2) (\theta_2 - \alpha(\vphi_2) \nabla \cdot \bm{u}_2 )^2 }_{L^p}\\ 
			& \quad + C(R) \norm{M'(\vphi_2) (\theta_2 - \alpha(\vphi_2) \nabla \cdot \bm{u}_2 )^2}_{L^p} \norm{\vphi_1 - \vphi_2}_{W^{1, p}}
		\end{aligned}\\ 
		& \quad \leq \begin{aligned}[t]
			C(R) \Big(  \big( \norm{\theta_1}_{L^q} +\norm{\theta_2}_{L^q} + 1 \big) \norm{\theta_1 - \theta_2}_{L^q} 
			&+ \big(\norm{\theta_1}_{L^q} + 1\big)  \norm{ \bm{u}_1 -  \bm{u}_2}_{\bm{W}^{1, q}} \\ 
			&+  \big(\norm{\theta_2}_{L^q} + 1\big) \norm{\vphi_1 - \vphi_2}_{W^{1, q}}
			\Big).
		\end{aligned}
	\end{align*}
	In particular, we arrive at
	\begin{align*}
		V&\leq C(R) \Big(  T^{\frac{1}{3r}}\big( \norm{\theta_1}_{L^{3r}(L^q)} +\norm{\theta_2}_{L^{3r}(L^q)} + 1 \big) \norm{\theta_1 - \theta_2}_{L^{3r}(L^q)} \\ 
		& \quad  +( T^{\frac{2}{3r}}\norm{\theta_1}_{L^{3r}(L^q)} + T^\frac{1}{r})\norm{ \bm{u}_1 -  \bm{u}_2}_{C^0(\X) }  \
		  + ( T^{\frac{2}{3r}}\norm{\theta_2}_{L^{3r}(L^q)} + T^\frac{1}{r})  \norm{\vphi_1- \vphi_2}_{C^0(W^{1, q})} 
		\Big)\\ 
		& \leq C(R, T) \Big(   \norm{\theta_1 - \theta_2}_{L^{3r}(L^q)} 
		 \norm{ \bm{u}_1 -  \bm{u}_2}_{C^0(\X) }  
		+  \norm{\vphi_1- \vphi_2}_{C^0(W^{1, q})} 
		\Big),
		\numberthis \label{est:F_1_VI}
	\end{align*}
	with $C(R, T) \rightarrow 0$ as $T \rightarrow 0$. 
	\par 
	\medskip 
	\textit{Ad $\F_T^1.VI$:} By the Lipschitz assumption \ref{A:source_terms} we further get 
	\begin{align*}
		VI &\leq 
		 \norm{S_s(\vphi_1,  \bm{u}_1, \theta_1) - S_s(\vphi_2,  \bm{u}_2, \theta_2)}_{L^r (W^{-2, p}_N) }  
		 \leq C(R, T) \norm{ (\vphi_1 - \vphi_2, \bm{u}_1 - \bm{u}_2, \theta_1 - \theta_2) }_{X_T}, 
		\numberthis \label{est:F_1_VII}
	\end{align*}
	where $C(T, R)\rightarrow 0$ as $T$ tends to zero. 
	\\ 
	Combining the estimates \eqref{est:F_1_I}-\eqref{est:F_1_VII} with the embeddings \eqref{embedding_X^1}, \eqref{embedding:hoelder} and \eqref{embedding_X^3} finally leads to 
	\begin{align*}
		&\norm{\F_T^1 (\vphi_1, \bm{u}_1, \theta_1) - \F_T^1 (\vphi_2, \bm{u}_2, \theta_2)}_{Y^1_T}
		  \leq C(R, T) \Big( \norm{\vphi_1 - \vphi_2}_{\new{Z}^1_T} 
		+  \norm{\bm{u}_1 - \bm{u}_2}_{\new{Z}^2_T}   
		+ \norm{\theta_1 - \theta_2}_{\new{Z}^3_T}
		  \Big)
		  	\numberthis \label{est:F_1}
		\end{align*}
	where $C(R, T) > 0$ is a constant depending on $R$ and $T$ such that $C(R, T) \rightarrow 0$ as $T \rightarrow 0$. 
	\par 
	\end{subequations}
	\begin{subequations}
	\medskip 
	\textit{Ad $\F_T^2$:} As before, we find
	\begin{align*}
		&\norm{\F_T^2 (\vphi_1, \bm{u}_1, \theta_1)  - \F_T^2 (\vphi_2, \bm{u}_2, \theta_2)}_{Y^2_T} \\ 
		 & \quad = \big\|  \big(\A(\vphi_0) - \A(\vphi_1)\big) \bm{u}_1 -  \big(\A(\vphi_0) - \A(\vphi_2)\big) \bm{u}_2 \big\|_{Y^2_T} \\ 		 
		& \qquad \quad   +\begin{aligned}[t]
			 \Big\|
			&\mathcal{B}^{-1} (\vphi_1) \Big( \nabla \cdot \big( \C(\vphi_1) \Tau (\vphi_1)  + \alpha(\vphi_1) M(\vphi_1) (\theta_1 - \alpha(\vphi_1) \nabla \cdot \bm{u}_1)  \textbf{I} \big) \new{+ \bm{f} +  \bm{g}} \Big) \\ 
			&- 	\mathcal{B}^{-1} (\vphi_2) \Big( \nabla \cdot \big( \C(\vphi_2) \Tau (\vphi_2)  + \alpha(\vphi_2) M(\vphi_2) (\theta_2 - \alpha(\vphi_2) \nabla \cdot \bm{u}_2) \textbf{I} \big)+ \bm{f} +  \bm{g} \Big) 
			\Big\|_{Y^2_T} 
		\end{aligned}\\ 
		& \quad \leq  \begin{aligned}[t]
			&   \Big\| \norm{\A(\vphi_0) - \A(\vphi_1) }_{\cL(\X)} \norm{\bm{u}_1 - \bm{u}_2}_{\X}  \Big\|_{L^r}
			 +  \Big\|  \norm{ \big(\A(\vphi_1)- \A(\vphi_2)}_{\cL(\X)} \norm{ \bm{u}_2 }_{\X}   \Big\|_{L^r}  \\ 
 			&+ C\Big\|  \begin{aligned}[t]
				\nabla \cdot \Big( &\C(\vphi_1) \Tau (\vphi_1)  + \alpha(\vphi_1) M(\vphi_1) (\theta_1 - \alpha(\vphi_1) \nabla \cdot \bm{u}_1) \textbf{I} \\ 
				&-   \C(\vphi_2) \Tau (\vphi_2)  + \alpha(\vphi_2) M(\vphi_2) (\theta_2 - \alpha(\vphi_2) \nabla \cdot \bm{u}_2)  \textbf{I} \Big)  \Big\| _{L^r(\Xd)}
			\end{aligned}\\ 
			& +\begin{aligned}[t]
				 & \Big\| \norm{\mathcal{B}^{-1} (\vphi_1) - \mathcal{B}^{-1} (\vphi_2)}_{\mathcal{L}(\Xd, \X)} \\ 
				& \quad \norm{\nabla \cdot \big( \C(\vphi_2) \Tau (\vphi_2)  + \alpha(\vphi_2) M(\vphi_2) (\theta_2 - \alpha(\vphi_2) \nabla \cdot \bm{u}_2) \textbf{I} \big)  + \bm{f} +  \bm{g}}_{\Xd} \Big\|_{L^r}
			\end{aligned}
		\end{aligned}\\ 
		& \quad \eqqcolon I + II + III + IV. 
	\end{align*}
	\\ 
	\textit{Ad $\F_T^2.I$:} Recall that by assumption $\A(\vphi) = \new{ \mathcal{B}^{-1}(\vphi) \mathcal{C}(\vphi)}$ is bounded independently of $\vphi \in \new{C^{\beta}(\overline\Omega)}$, $\beta > 0$, cf. \ref{A:iso}. Hence, 
	\begin{equation*}
		I = \Big\| \norm{\A(\vphi_0) - \A(\vphi_1) }_{\cL(\Xd) } \norm{ \bm{u}_1 - \bm{u}_2}_{\X}  \Big\|_{L^r}  
		 \leq  C T^{\frac{1}{r} } \norm{ \bm{u}_1 - \bm{u}_2 }_{C^0(\X)} . 
	\end{equation*}
	\\ 
	\textit{Ad $\F_T^2.III$:} For $(III)$ we can argue similarly as before, \new{cf. ($\F_T^1.IV$)}, 
	\begin{align*}
		&\begin{aligned}[t]
			\Big\| \nabla \cdot \Big( &\C(\vphi_1) \Tau (\vphi_1)  + \alpha(\vphi_1) M(\vphi_1) (\theta_1 - \alpha(\vphi_1) \nabla \cdot \bm{u}_1) \textbf{I} \\ 
			&-   \C(\vphi_2) \Tau (\vphi_2)  + \alpha(\vphi_2) M(\vphi_2) (\theta_2 - \alpha(\vphi_2) \nabla \cdot \bm{u}_2)  \textbf{I} \Big)  \Big\| _{\Xd}
	\end{aligned}\\ 
	& \quad \leq \begin{aligned}[t]
		\norm{ \C(\vphi_1) \Tau (\vphi_1) 	-   \C(\vphi_2) \Tau (\vphi_2)}_{L^q}
		&+ \norm{ \alpha(\vphi_1) M(\vphi_1) \theta_1 - \alpha(\vphi_2) M(\vphi_2) \theta_2 }_{L^q}\\ 
		&+ \norm{ \alpha(\vphi_1)^2 M(\vphi_1) \nabla \cdot \bm{u}_1  - \alpha(\vphi_2)^2 M(\vphi_2) \nabla \cdot \bm{u}_2}_{L^q}
	\end{aligned}\\ 
	& \quad \leq C(R) \Big( (1 + \norm{\theta_2}_{L^q}) \norm{ \vphi_1 - \vphi_2 }_{W^{1, q}} + \norm{\bm{u}_1 - \bm{u}_2}_{\X} + \norm{\theta_1 - \theta_2}_{L^q}  \Big), 
	\end{align*}
	and conclude 
	\begin{align*}
		III \leq C(R) \Big( (T^{\frac{1}{r}} + T^{\frac{2}{3r}})  \norm{ \vphi_1 - \vphi_2 }_{C^0(W^{1, q})} 
		+ T^{\frac{1}{r}} \norm{\bm{u}_1 - \bm{u}_2}_{C^0(\bm{W}^{1, q})} 
		+ T^{\frac{2}{3r}}\norm{\theta_1 - \theta_2}_{L^{3r}(L^q)}  \Big) . 
		\numberthis \label{est:F_2_III}
	\end{align*}
	\\ 
	\textit{Ad $\F_T^2.IV$:} To estimate $(IV)$, we first note that 
	\begin{align*}
		 &\norm{\nabla \cdot \big( \C(\vphi_2) \Tau (\vphi_2)  + \alpha(\vphi_2) M(\vphi_2) (\theta_2 - \alpha(\vphi_2) \nabla \cdot \bm{u}_2)  \textbf{I}\big) +  \bm{f} + \bm{g}} _{\Xd}  \\ 
		 &\quad\leq C(R) ( \norm{\theta_2}_{L^q}+  \norm{\bm{f}}_{\Xd} + \norm{\bm{g}}_{\bm{W}^{-(1-\frac{1}{q'})}(\Gamma_N) }   + 1 ). 
	\end{align*}
	To estimate the operator norm, we start by examining the following difference
	\begin{align*}
		& \norm{\mathcal{B} (\vphi_1) - \mathcal{B} (\vphi_2)}_{\mathcal{L}(\X, \Xd )} 
		\leq  \sup_{\norm{\bm{f}}_{ \X} = 1 } \sup_{ \norm{\bm{\eta}}_{\Xp} = 1} \abs{ \int_\Omega [\C_\nu (\vphi_1) - \C_\nu(\vphi_2)] \E(\bm{f} )\colon \E(\eta) \dx  }\\ 
		& \quad \leq  \sup_{\norm{\bm{f}}_{ \X} = 1 }  \norm{  [\C_\nu (\vphi_1) - \C_\nu(\vphi_2)] \E(\bm{f} )}_{\bm{L}^q}
		\leq  \sup_{\norm{\bm{f}}_{\X} = 1 }   \norm{\C_\nu (\vphi_1) - \C_\nu(\vphi_2)}_{\bm{L}^\infty } \norm{\E(\bm{f} )}_{\bm{L}^q}\\ 
		& \quad \leq C \norm{\vphi_1 - \vphi_2}_{W^{1, q} }. 
		\numberthis \label{eq:B_Lipschitz}
	\end{align*}
	Here, we used the Lipschitz continuity of $\C_\nu$, \new{which follows from $\C'_\nu$ being bounded, cf.\ \ref{A:C_nu}}, and the embedding $W^{1, q}(\Omega) \hookrightarrow L^\infty(\Omega)$. 
	Exploiting this estimate along with the uniform bound on $\mathcal{B}^{-1}(\cdot)$ imposed in \ref{A:iso}, we deduce 
	\begin{align*}
		&\norm{\mathcal{B}^{-1} (\vphi_1) - \mathcal{B}^{-1} (\vphi_2)}_{\mathcal{L}(\Xd, \X)} 
		 = \norm{ {B}^{-1} (\vphi_1)  [ \mathcal{B}(\vphi_2) - \mathcal{B} (\vphi_1)]   \mathcal{B}^{-1} (\vphi_2) }_{\cL(\Xd, \X)} \\ 
		& \quad \leq 
		\norm{ {B}^{-1} (\vphi_1)}_{\cL(\Xd, \X)}
		 \norm{ \mathcal{B}(\vphi_2) - \mathcal{B} (\vphi_1) }_{\cL(\X, \Xd)}   \norm{\mathcal{B}^{-1} (\vphi_2) }_{\cL(\Xd, \X)}\\ 
		 & \quad \leq C \norm{\vphi_1 - \vphi_2}_{W^{1, q}} 
		 	\numberthis \label{eq:B^-1_Lipschitz}
	\end{align*}
	and conclude 
	\begin{align*}
		IV \leq &C(R) \big(  T^{\frac{2}{3r}}   \norm{\theta_2}_{L^{3r}(L^q)} +    T^{\frac{1}{r}}  \big) \norm{\vphi_1 - \vphi_2}_{C^0(W^{1, q}) }\\ 
		 &+ C(R) T^{\tilde{\vartheta}}  \big( \norm{\bm{f}}_{L^r(\Xd) } + \norm{\bm{g}}_{L^r( \bm{W}^{-(1-\frac{1}{q'})}(\Gamma_N) )} \big)   \norm{\vphi_1 - \vphi_2}_{C^{\tilde{\vartheta}}(W^{1, q})}.  
		\numberthis \label{est:F_2_IV}
	\end{align*}
	\\ 
	\textit{Ad $\F_T^2.II$:} Finally, we use the definition of $\A$,  Lipschitz-continuity of $\mathcal{B}^{-1}$ and $\mathcal{C}$, see \eqref{eq:B_Lipschitz}, \eqref{eq:B^-1_Lipschitz}, and the uniform boundedness of these operators and their inverses to compute 
	\begin{align*}
		\big\| &\A(\vphi_1)- \A(\vphi_2)\big\|_{\cL(\X)}
		= \norm{ \mathcal{B}^{-1}(\vphi_1) \mathcal{C}(\vphi_1) - \mathcal{B}^{-1}(\vphi_2) \mathcal{C}(\vphi_2) }_{\cL(\X)} \\ 
		& \leq \norm{\mathcal{B}^{-1}(\vphi_1)}_{\cL(\Xd, \X)} \norm{\mathcal{C}(\vphi_1) - \mathcal{C}(\vphi_2) }_{\cL(\X, \Xd)} \\ 
		& \quad + \norm{ \mathcal{B}^{-1}(\vphi_1) - \mathcal{B}^{-1}(\vphi_2)}_{\cL(\Xd, \X)} \norm{ \mathcal{C}(\vphi_2)}_{\cL(\X, \Xd)}\\ 
		& \leq C \norm{\vphi_1 - \vphi_2}_{W^{1, q}}. 
	\end{align*}
	This estimate gives rise to 
	\begin{align*}
		II \leq C(R) T^{\frac{1}{r}} \norm{\vphi_1- \vphi_2}_{C^0(W^{1, q})}. 
	\end{align*}
	Now it is again easy to combine \eqref{est:F_2_III} and \eqref{est:F_2_IV},  
	\begin{align*}
		 &\norm{\F_T^2 (\vphi_1, \bm{u}_1, \theta_1)  - \F_T^2 (\vphi_2, \bm{u}_2, \theta_2)}_{Y^2_T} 
		 \leq C(R, T) \Big( \norm{\vphi_1 - \vphi_2}_{\new{Z}^1_T} 
		 +  \norm{\bm{u}_1 - \bm{u}_2}_{\new{Z}^2_T}   
		 + \norm{\theta_1 - \theta_2}_{\new{Z}^3_T}
		 \Big)
		 	\numberthis \label{est:F_2}
	\end{align*}
	where $C(R, T) \rightarrow 0$ as $T \rightarrow 0$. 
		\end{subequations}
		\par 
		\medskip 
		\noindent
		\begin{subequations}
	\textit{Ad $\F_T^3$:} Lastly, we have 
	\begin{align*}
		&\norm{\F_T^3 (\vphi_1, \bm{u}_1, \theta_1) - \F_T^3 (\vphi_2, \bm{u}_2, \theta_2)}_{Y^3_T}\\ 
		& \quad \leq \norm{ - \nabla \cdot \big( \kappa(\vphi_0) M(\vphi_0)  \nabla  \big(\theta_1 - \theta_2 \big)  \big) 
			+ \nabla \cdot  \Big( \kappa(\vphi_1) M(\vphi_1)  \nabla \theta_1 -  \kappa(\vphi_2)  M(\vphi_2)  \nabla  \theta_2 \Big) }_{Y^3_T}\\ 
		& \quad  \quad + \norm{ \nabla \cdot  \Big( \kappa(\vphi_1) M'(\vphi_1) \nabla \vphi_1\,  \theta_1  -  \kappa(\vphi_2) M'(\vphi_2) \nabla \vphi_2 \, \theta_2 \Big) }_{Y^3_T} \\ 
		& \quad \quad + \norm{\nabla \cdot \Big(
			\kappa(\vphi_1) \nabla \big(M(\vphi_1) \alpha(\vphi_1) \nabla \cdot \bm{u}_1 \big)  
			- \kappa(\vphi_2) \nabla \big(M(\vphi_2) \alpha(\vphi_2) \nabla \cdot \bm{u}_2 \big)
			\Big)  }_{Y^3_T} \\ 
		& \quad \quad + \norm{S_f(\vphi_1, \bm{u}_1,  \theta_1) - S_f(\vphi_2, \bm{u}_2,  \theta_2)}_{Y^3_T} \\  
		& \quad \eqqcolon I + II + III +IV . 
	\end{align*}
	\end{subequations}
	\\ 
	\begin{subequations}
	\textit{Ad $\F_T^3.I$:} Similar to above, we set $\hat{\kappa}(\vphi_i) \coloneqq  \kappa(\vphi_i) M(\vphi_i)$ and compute with the help of \eqref{mult:-1q} and Theorem~\ref{thm:compositon}
	\begin{align*}
		&\norm{- \nabla \cdot \big(\hat{ \kappa}(\vphi_0)  \nabla  \big(\theta_1 - \theta_2 \big)  \big) 
			+ \nabla \cdot  \Big( \hat{\kappa}(\vphi_1)   \nabla \theta_1 -  \hat{\kappa}(\vphi_2)    \nabla  \theta_2 \Big) }_{W^{-2, q}_N} \\ 
		& \quad \leq C \Big( 
		\norm{  \big( \hat{  \kappa}(\vphi_0) -  \tilde{\kappa}(\vphi_1) \big) \nabla  \big(\theta_2 - \theta_1 \big)  \big)  }_{W^{-1, q}_N} 
		+ \norm{ \big( \tilde{  \kappa}(\vphi_1) -  \tilde{\kappa}(\vphi_2) \big)  \nabla \theta_2}_{W^{-1, q}_N} \Big) \\ 
		& \quad \leq C(R) \Big( \norm{ \vphi_0 - \vphi_1 }_{W^{1, q}}  \norm{ \theta_2 - \theta_1   }_{L^q} 
		+  \norm{ \vphi_1 - \vphi_2 }_{W^{1, q}}  \norm{ \theta_2}_{L^q} 
		\Big). 
	\end{align*}
	As above, see \eqref{est:F_1_I}, we exploit \eqref{embedding:hoelder} to deduce 
	\begin{align*}
		I \leq C(R) T^{\bar{\vartheta}}  
		\Big( \norm{\vphi_1}_{C^{ \bar{\vartheta}}}(W^{1, q}) \norm{\theta_1 - \theta_2}_{L^{3r}(L^q)} 
		+   \norm{\vphi_1 - \vphi_2}_{C^{ {\vartheta}}(W^{1, q})}  \norm{ \theta_2}_{L^{3r}(L^q)}  \Big). 
		\numberthis \label{est:F_3_I}
	\end{align*}
	\\ 
	\textit{Ad $\F_T^3.II$:} Using the notation $\tilde{\kappa}(\vphi_i) \coloneqq  \kappa(\vphi_i) M'(\vphi_i)$ along with the multiplication results \eqref{mult:-1q} and \eqref{mult:triple}, it follows that 
	\begin{align*}
		&\norm{  \nabla \cdot  \Big( \tilde{\kappa}(\vphi_1) \nabla \vphi_1\,  \theta_1  -  \tilde{\kappa}(\vphi_2) \nabla \vphi_2 \, \theta_2 \Big)  }_{W^{-2, q}_N} \\ 
		& \ \  \leq C \norm{   \tilde{\kappa}(\vphi_1) \nabla \vphi_1\,  \theta_1  -  \tilde{\kappa}(\vphi_2) \nabla \vphi_2 \, \theta_2 }_{W^{-1, q}} \\ 
		& \ \ 
		\leq C  \Big( \norm{\tilde{\kappa}(\vphi_1) \big( \nabla \vphi_1\,  \theta_1  -   \nabla \vphi_2\, \theta_2  \big) }_{W^{-1, q}}
		+ \norm{ \big( \tilde{\kappa}(\vphi_1) - \tilde{\kappa}(\vphi_2) \big) \nabla \vphi_2 \theta_2 }_{W^{-1, q}} 
		\Big) \\ 
		& \ \ \leq \begin{aligned}[t]
			C \norm{\tilde{\kappa} (\vphi_1)}_{W^{1, q}} \norm{\nabla \vphi_1}_{W^{2\delta, q}} \norm{ \theta_1 - \theta_2}_{W^{-\delta, q}}
			&+ C  \norm{\tilde{\kappa} (\vphi_1)}_{W^{1, q}}  \norm{\theta_2}_{W^{-\delta, q}} \norm{\nabla \big( \vphi_1 - \vphi_2 \big)}_{W^{2\delta, q}}\\ 
			& + C \norm{ \tilde{\kappa} (\vphi_1) - \tilde{\kappa} (\vphi_2)}_{W^{1, q}} \norm{ \nabla \vphi_2 }_{W^{2\delta, q}} \norm{\theta_2}_{W^{-\delta, q}}
		\end{aligned} \\ 
		& \ \ \leq C(R) \Big( \norm{\vphi_1 - \vphi_2 }_{W^{1 + 2 \delta, q}} + \norm{ \theta_1 - \theta_2}_{W^{- \delta, q}} \Big). 
	\end{align*}
	Note that we also applied Theorem \ref{thm:compositon} to $ \norm{ \tilde{\kappa} (\vphi_1) - \tilde{\kappa} (\vphi_2)}_{W^{1, q}}$. Along with the embedding results \eqref{embedding_X^1}, \eqref{embedding_X^3}, this leads to 
	\begin{align*}
		II \leq C(R) T^{\frac{1}{3r}} \Big(  \norm{\vphi_1 - \vphi_2 }_{C^0(W^{1 + 2 \delta, q})} 
		+ \norm{ \theta_1 - \theta_2}_{C^0(W^{- \delta, q})}
		 \Big). 
		 \numberthis \label{est:F_3_III}
	\end{align*}
	\\ 
	\textit{Ad $\F_T^3.III$:} We observe that exploiting the multiplication result \eqref{mult:-1q} instead of \eqref{mult:-1p} leads to the analogous estimate as in Lemma~\ref{lemma:general_estimate}, but in the space $W^{-2, q}_N(\Omega)$. Applying this result, we deduce 
	\begin{align*}
		 &\norm{\nabla \cdot \Big(
			\kappa(\vphi_1) \nabla \big(M(\vphi_1) \alpha(\vphi_1) \nabla \cdot \bm{u}_1 \big)  
			- \kappa(\vphi_2) \nabla \big(M(\vphi_2) \alpha(\vphi_2) \nabla \cdot \bm{u}_2 \big)
			\Big)  }_{W^{-2, q}_N} \\ 
		& \quad \leq  C (R)  \Big(  \norm{ M(\vphi_1) \alpha(\vphi_1) \nabla \cdot \bm{u}_1 - M(\vphi_2) \alpha(\vphi_2) \nabla \cdot \bm{u}_2}_{L^q}  + \norm{M(\vphi_2) \alpha(\vphi_2) \nabla \cdot \bm{u}_2}_{L^q} \norm{\vphi_1 - \vphi_2}_{W^{1 ,q}} \Big) \\ 
		& \quad \leq  \begin{aligned}[t]
			C (R)  \Big(  
			\norm{  M(\vphi_1) \alpha(\vphi_1)  \nabla \cdot (\bm{u}_1 - \bm{u}_2)}_{L^q} &+ \norm{ \big( M(\vphi_1) \alpha(\vphi_1) -  M(\vphi_2) \alpha(\vphi_2)\big) \nabla \cdot \bm{u}_2}_{L^q} \\ 
			&+ \norm{M(\vphi_2)}_{L^\infty} \norm{\bm{u}_2}_{\bm{W}^{1, q}} \norm{\vphi_1 - \vphi_2}_{W^{1 ,q}}
			\Big)
		\end{aligned}\\ 
		& \quad \leq C (R)  \Big(   
			 \norm{\vphi_1 - \vphi_2}_{W^{1 ,q}} 
			 + \norm{\bm{u}_1 - \bm{u}_2}_{\X} 
		\Big)
	\end{align*}
	and arrive at  
	\begin{align*}
		III \leq C(R) T^{\frac{1}{3r}} \Big( \norm{\vphi_1 - \vphi_2}_{C^0(\X)} 
		+ \norm{\bm{u}_1 - \bm{u}_2}_{C^0(\bm{W}^{1, q})}
		  \Big). 
		  	\numberthis \label{est:F_3_IV}
	\end{align*}
	\\ 
	\textit{Ad $\F_T^3.IV$:} Relying on assumption \ref{A:source_terms}, we conclude exactly as in \eqref{est:F_1_VII} that 
	\begin{align*}
			IV &\leq \ \norm{S_f(\vphi_1,  \bm{u}_1, \theta_1) - S_f(\vphi_2,  \bm{u}_2, \theta_2)}_{L^r (W^{-2, p}_N) }  
			\leq C(R, T) \norm{ (\vphi_1 - \vphi_2, \bm{u}_1 - \bm{u}_2, \theta_1 - \theta_2) }_{\new{Z}_T}, 
			\numberthis \label{est:F_3_V}
	\end{align*}
	where $C(R, T)$ tends to zero as $T \rightarrow 0$. 
	\\ 
	In summary, \eqref{est:F_3_I}-\eqref{est:F_3_V} lead to 
		\begin{align*}
		&\norm{\F_T^3 (\vphi_1, \bm{u}_1, \theta_1)  - \F_T^3 (\vphi_2, \bm{u}_2, \theta_3)}_{Y^3_T} 
		\leq C(R, T) \Big( \norm{\vphi_1 - \vphi_2}_{\new{Z}^1_T} 
		+  \norm{\bm{u}_1 - \bm{u}_2}_{\new{Z}^2_T}   
		+ \norm{\theta_1 - \theta_2}_{\new{Z}^3_T}
		\Big)
		\numberthis \label{est:F_3}
	\end{align*}
	for some constant $C(R, T) > 0$ such that $C(R, T) \rightarrow 0$ as $T \rightarrow 0$, which concludes the proof. 
	\end{subequations}

\end{proof}

\medspace

\fakepart{}

	\section*{Acknowledgments} 
	The second author is supported by the Graduiertenkolleg 2339 IntComSin of the Deutsche Forschungsgemeinschaft (DFG, German Research Foundation) – Project-ID 321821685. The support is gratefully acknowledged. Moreover, we would also like to thank the anonymous referee for helpful comments leading to an improvement of this paper.

	\addtocontents{toc}{\SkipTocEntry}
	\subsection*{Conflict of interests and data availability statement} There is no conflict of interests.
	There is no associated data to the manuscript.

	\bibliographystyle{siam}
	{\small
	\bibliography{refs}
	}
	\appendix

\section{$\textbf{W}^{2,p}$-regularity for elliptic systems} \label{sec:w^2^p}

The aim of this section is to provide a proof of the regularity result for elliptic systems with mixed boundary conditions as claimed in Theorem \ref{thm:higher_elliptic}. While there are many classical references, e.g., \cite{MR125307}, \cite{giaquintaMatrinazzi}, \cite{gilbarg2001elliptic}, \cite{Valent_Boundary}, most of these do not treat the problem in its full generality, restricting the exposition to the scalar case, the interior domain, a single type of boundary conditions or simply assume more regular coefficients than required.\\ 
We will use a different approach, beginning with a regularity result for elliptic systems in non-divergence form derived from parabolic maximal $L^r$-regularity theory. Thereafter, we proceed with an iterative argument that allows us to transfer the implications to systems in divergence form. \par 
\medskip 
\noindent
To simplify the notation, we will, for now, neglect the dependence of $\C$ on $\vphi$ and simply consider $\C = \C(\bm{x})$.   
First of all, let us recall the standard symmetry assumptions in linear elasticity  
\begin{equation*}
	\C_{ijkl} = \C_{jikl}= \C_{ijlk} = \C_{klij}, 
\end{equation*}
which imply 
\begin{align*}
	[\C \, \E(\bm{u}) ]_{ij} = [\C \,  \nabla \bm{u}]_{ij}. \numberthis \label{eq:symmetry_C}
\end{align*}
Therefore, it suffices to consider the elliptic mixed boundary problem 
	\begin{align*}
	\left\{ 
	\begin{alignedat}{2}
		- \nabla \cdot \big( \C \nabla\bm{v}   \big)  &= \bm{f} & & \quad \textnormal{in } \Omega,\\
		\C \nabla\bm{v} \bm{n} &= \bm{g} & &\quad \textnormal{on }  \Gamma_N, \\ 
			\bm{v}  &= \bm{0} & &\quad \textnormal{on }   \Gamma_D. 
	\end{alignedat} 
	\right. 
\end{align*}
However, we will first need to establish the corresponding results for the system in non-divergence form, i.e., 
\begin{equation*}
	-  \textnormal{tr}( \C \nabla^2 \bm{u})  - (\nabla \cdot \C )\, \nabla \bm{u} = \bm{f}  \quad \textnormal{in } \Omega,  \\ 
\end{equation*}
where the first multiplication is defined by 
\begin{equation*}
	\textnormal{tr} (\C \, \nabla^2 \bm{u})_i = \C_{ijkl} \, \partial_j \partial_l \bm{u}^k 
	= \sum_{j, l = 1}^n \sum_{k = 1}^n \C_{ijkl} \, \partial_j \partial_l \bm{u}^k. 
\end{equation*}
For fixed $j, l$, the sum over $k$ is the matrix-vector product of $\C_{ijkl}$ and the second-order derivative $\partial_j \partial_l \bm{u}$
\begin{equation*}
	\C_{ijkl}\,  \partial_j \partial_l \bm{u} =  \sum_{k = 1}^n \C_{ijkl}\, \partial_j \partial_l \bm{u}^k. 
\end{equation*}
Let $\alpha \in \N_{>0}^n$ be a multindex with $\abs{\alpha} = \sum \alpha_i$ and $D^\alpha  = i^{\abs{\alpha}} \partial_1^{\alpha_i}\cdots \partial_{n}^{\alpha_n}$, then we can rewrite the second-order derivatives as 
\begin{equation*}
	- \textnormal{tr} (\C \, \nabla^2 \bm{u} )= \sum_{j, l = 1}^n \C_{ijkl} \, i^2 \partial_j \partial_l \bm{u} = \sum_{\abs{\alpha} = 2} \C_{ijkl} \, D^{\alpha} \bm{u},  
\end{equation*}
where $j, l$ are the non-zero indices of $\alpha$, where $j = l$ is explicitly allowed. 
Hence, the principal symbol of this differential operator is given by
\begin{equation*}
	\mathcal{C}_{\#}(\bm{x}, \xi) \coloneqq \sum_{\abs{\alpha} = 2} \C_{ijkl} \, \xi^{\alpha} 
	= \sum_{j, l = 1}^n \C_{ijkl} \, \xi_j \xi_l \in \R^{n \times n}, \quad \textnormal{for all } \xi \in \R^n, 
\end{equation*}
and we want to show that $\mathcal{C}_{\#}$ is parameter elliptic with angle of ellipticity $\phi_{\mathcal{C}_{\# }} < \frac{\pi}{2} $, see \eqref{def:parameter_elliptic}. In order to prove this, let us recall that we call $\C_{\#}(\bm{x}, \cdot)$ \textit{strongly elliptic}, also known as \textit{Legendre} condition \cite[Def. 6.1.3]{PrussSimonett2016}, if 
\begin{equation*}
	\textnormal{Re}  \langle \mathcal{C}_{\#}(\bm{x}, \xi) \eta, \eta  \rangle_{\C^n} \geq c \norm{\xi}^2 \norm{\eta}^2 
	\quad \textnormal{for all } \quad \xi \in \R^n, \eta \in \C^n. 
\end{equation*}
By assumption \ref{A:W} we know that $\C$ is uniformly positive definite on symmetric matrices and the symmetry of $\C$ further entails that for any $\xi, \eta \in \R^n$ it holds that
\begin{align*}
	\langle \sum_{j, l = 1}^n
	&\C_{ijkl} \, \zeta_j \zeta_l \, \eta , \eta \rangle_{\R^n}
{\geq c \big( \norm{\eta}^2 \norm{\xi}^2 + (\eta \cdot \xi)^2 \big)}
	    \geq c \norm{\eta}^2 \norm{\xi}^2. 
\end{align*}
Hence, we compute for any $\eta = \eta^r + i \eta^i \in \C^n$ 
\begin{align*}
		\textnormal{Re}  \langle \mathcal{C}_{\#}(\bm{x}, \xi) ( \eta^r + i \eta^i),  \eta^r + i \eta^i   \rangle_{\C^n}
		&= \langle\sum_{j, l = 1}^n \C_{ijkl} \, \zeta_j \zeta_l\,  \eta^r , \eta^r \rangle_{\R^n}  +
		\langle \sum_{j, l = 1}^n \C_{ijkl} \, \zeta_j \zeta_l \, \eta^i, \eta^i \rangle_{\R^n} \\ 
		& \geq c \norm{\xi}^2 ( \norm{\eta^r}^2 + \norm{\eta^i}^2 ) = c \norm{\xi}^2 \norm{\eta}^2 
\end{align*}
and see that $\mathcal{C}_{\#}(\bm{x} , \cdot)$ is uniformly strongly elliptic. Thus, \cite[Rem.~6.1.4(ii)]{PrussSimonett2016} yields that $\mathcal{C}_{\#}$ is parameter elliptic with angle of ellipticity $\phi_{\mathcal{C}_{\# }} < \frac{\pi}{2}$. We also refer to \cite[Sec. 6.2.5]{PrussSimonett2016} for further details and note that the operator $\C(\vphi) \nabla^2 \bm{u}$ together with mixed boundary conditions satisfy the Lopatinskii-Shapiro conditions, cf. \cite[Prop. 6.2.13, Rem.~(ii)]{PrussSimonett2016}. \\ 
Thus, we can invoke \cite[Thm. 6.3.2]{PrussSimonett2016}, deducing that there exits some $\omega > 0$ such that the operator $-\C \, \nabla^2 \bm{u} - (\nabla \cdot \C )\, \nabla \bm{u}+ \omega \bm{u}$ combined with inhomogeneous mixed boundary conditions, and a sufficiently regular $\C$, exhibits maximal $L^r$-regularity on the whole real line. We emphasize that $\omega$ only depends on the ellipticity estimates on $\C$ from above. Hence \cite[Prop. 3.5.2 (ii)]{PrussSimonett2016} entails the following theorem.

\begin{theorem}\label{thm:w^2_p_non_div}
	Let $\Omega \subset \R^n$ be a domain with $C^{2}$-boundary and suppose $\Gamma_D  \subset \partial \Omega$ is relatively closed with $\mathcal{H}^{n-1} (\Gamma_D) > 0$. Define $\Gamma_N = \partial \Omega \setminus \Gamma_D$ and assume $\overline{\Gamma_D} \cap \overline{\Gamma_N} = \emptyset$. Moreover, suppose $1 < p \leq q < \infty$, $q > n$ and that $\C \in \bm{W}^{1, q}({\Omega}, \R^{4n})$ is uniformly strongly elliptic. Then, for every $\bm{f} \in \bm{L}^p(\Omega)$ and $\bm{g} \in \bm{W}^{1- \frac{1}{p}, p}(\Gamma_N)$ there exists a unique solution $\bm{u}$ to the elliptic boundary value problem 
	\begin{align*}
		- \textnormal{tr}( \C(\bm{x})    \nabla^2 \bm{u})  - (\nabla \cdot \C(\bm{x})) \nabla \bm{u} + \omega \bm{u} &= \bm{f} \quad \textnormal{in } \Omega, \\ 
		\C(\bm{x})   \nabla \bm{u} \, \bm{n}&= \bm{g} \quad \textnormal{on }  \Gamma_N, \\ 
		\bm{u} & = \bm{0} \quad \textnormal{on }  \Gamma_D.  
	\end{align*}
	Furthermore, there exists a constant $C > 0$, which is independent of $\bm{u}$ but may depend on the norm of $\C$, its ellipticity parameter, and the modulus of continuity, such that 
	\begin{equation*}
		\norm{\bm{u}}_{\bm{W}^{2, p}} \leq C \big(  \norm{\bm{f}}_{\bm{L}^p} +  \norm{\bm{g}}_{\bm{W}^{1- \frac{1}{p}, p}(\Gamma_N)}  \big). 
	\end{equation*}
\end{theorem}

\begin{theorem}\label{thm:elastic_non_symmetric}
	Suppose that $\C \in \bm{C}^2(\R, \R^{4n})$ is uniformly strongly elliptic, $\phi \in W^{1, q }(\Omega)$ for some $q > n \geq 2$, and that $\Omega$ has the same properties as in Theorem \ref{thm:w^2_p_non_div}.
	Moreover, assume that the data satisfies $\bm{f} \in \bm{L}^{p}(\Omega)$ and ${ \nolinebreak\bm{g} \in \bm{W}^{1- \frac{1}{p}, p} (\Gamma_N)}$ for some $p \in [2, q]$. Then, the elliptic problem 
	\begin{align}\label{eq:div_elasticity}
		\begin{alignedat}{1}
					- \nabla \cdot  \big( \C(\phi)   \nabla \bm{u} \big) &= \bm{f} \quad \textnormal{in } \Omega, \\ 
			\C(\phi)  \nabla \bm{u}\,  \bm{n}&= \bm{g} \quad \textnormal{on }  \Gamma_N, \\ 
			\bm{u} & = \bm{0} \quad \textnormal{on }  \Gamma_D, 
		\end{alignedat}
	\end{align}
	has a unique solution $\bm{u} \in \bm{W}^{2, p}(\Omega)$ and there exists a constant $C > 0$, which is independent of $\bm{u}$ but may depend on the norm $\norm{\C(\vphi)}_{W^{1, q}}$ and its ellipticity parameter, such that 
	\begin{equation*}
		\norm{\bm{u}}_{\bm{W}^{2, p}} \leq C \big(  \norm{\bm{f}}_{\bm{L}^p} +  \norm{\bm{g}}_{\bm{W}^{1- \frac{1}{p}, p} (\Gamma_N)}   \big). 
	\end{equation*}
\end{theorem}

\begin{proof}
	The idea of this proof is to iteratively apply Theorem \ref{thm:w^2_p_non_div} until we obtain the desired regularity.\\
	First of all, we observe, that due to the uniform strong ellipticty and Korn's inequality, we can invoke the Lax-Milagram theorem and deduce that there exists a unique weak solution $\bm{u} \in \bm{W}^{1, 2}_{\Gamma_N}(\Omega)$ to the system \eqref{eq:div_elasticity}, which further satisfies 
	\begin{equation}\label{eq:w^2_non_div_1}
		\norm{\bm{u}}_{\bm{W}^{1, 2}_{\Gamma_D}} 
		\leq C \big( \norm{\bm{f}}_{\bm{L}^2} + \norm{\bm{g}}_{\bm{W}^{-\frac{1}{2}, 2}} \big)
	\end{equation}
	for some $C > 0$. \\ 
	Moreover, we can exploit Theorem \ref{thm:compositon} to find $\C(\vphi) \in \bm{W}^{1, q}(\Omega)$, allowing us to consider the system
	\begin{alignat*}{2}
		- \textnormal{tr} (\C(\vphi)    \nabla^2 \bm{v} ) - (\nabla \cdot \C(\vphi)) \nabla \bm{v} + \omega \bm{v} &= \bm{f} + \omega \bm{u} \quad &&\textnormal{in } \Omega, \\ 
		\C(\bm{x})   \nabla \bm{v} \, \bm{n}&= \bm{g} \quad &&\textnormal{on }  \Gamma_N, \\ 
		\bm{v} & = \bm{0} \quad &&\textnormal{on }  \Gamma_D. 
	\end{alignat*}
	We take a closer look at the right-hand side $\bm{f} + \omega \bm{u}$ and distinguish between two cases: \par 
	\medskip 
	\underline{Case 1:} In this case, we suppose that $\bm{W}^{1, 2}(\Omega) \hookrightarrow \bm{L}^p(\Omega)$, entailing that $\bm{f} + \omega \bm{u} \in \bm{L}^p(\Omega)$. Hence, Theorem \ref{thm:w^2_p_non_div} yields the existence of a unique solution $\bm{v}$ along with the estimate 
	\begin{align*}
		\norm{\bm{v}}_{\bm{W}^{2, p}} \leq 
		C \big(  \norm{\bm{f} + \omega \bm{u}}_{\bm{L}^p} + \norm{\bm{g}}_{\bm{W}^{1 -\frac{1}{p}, p}} \big)
		&\leq C \big(  \norm{\bm{f} }_{\bm{L}^p} + \norm{\bm{u}}_{\bm{W}^{1, 2}_{\Gamma_D}}+ \norm{\bm{g}}_{\bm{W}^{1 -\frac{1}{p}, p}} \big)\\ 
		& \leq C \big(  \norm{\bm{f} }_{\bm{L}^p} + \norm{\bm{g}}_{\bm{W}^{1 -\frac{1}{p}, p}} \big), 
	\end{align*}
	where we used \eqref{eq:w^2_non_div_1} together with $p \geq 2$ for the last inequality. Moreover, it also follows that $\bm{v}$ is another weak solution to \eqref{eq:div_elasticity}, which entails by virtue of uniqueness that $\bm{v}$ and $\bm{u}$ coincide, concluding the proof in this case. \par 
	\medskip 
	\underline{Case 2:} If however $\bm{W}^{1, 2}(\Omega) \not\hookrightarrow \bm{L}^p(\Omega)$, which can only happen if $n \geq 3$, then we obtain by the same argument as above $\bm{u} \in \bm{W}^{2, p^*}(\Omega)$, where $p^* $ is the Sobolev exponent $\frac{2n}{n-2} > 2$, along with the estimate 
	\begin{equation*}
		\norm{\bm{u}}_{\bm{W}^{2, p^*}_{\Gamma_D}} 
		\leq C \big( \norm{\bm{f}}_{\bm{L}^p} + \norm{\bm{g}}_{\bm{W}^{1- \frac{1}{p}, p}} \big). 
	\end{equation*}
	Again there are two cases to consider: \par
	\underline{\textit{Case 2.1}} Similar to the first case, the assertion directly follows after one more iteration if it holds that $\bm{W}^{2, p^*}(\Omega) \hookrightarrow L^p(\Omega)$. We also observe that this condition is always satisfied if $2p^* \geq n$.\par 
	\underline{\textit{Case 2.2}} Conversely if $\bm{W}^{2, p^*}(\Omega) \not\hookrightarrow L^p(\Omega)$, and therefore that $2p^* < n$, we obtain by another application of Theorem \ref{thm:w^2_p_non_div} that $\bm{u} \in \bm{W}^{2, \tilde{p^*}}(\Omega)$, where $\tilde{p^*} = \frac{np^*}{n - 2p^*}$, along with the corresponding estimate. \\ 
	Repeating this argument increases the regularity iteratively and it only remains to show that a finite number of iterations suffices to attain the desired result. To this end, we recall that 
	\begin{equation*}
		\frac{1}{\tilde{p^*}} = \frac{1}{p^*} - \frac{2}{n}, 
	\end{equation*}
	which implies that $\frac{1}{\tilde{p^*}}$ increases by $\frac{2}{n}$ in each iteration. Hence, a finite amount of repetitions suffice to either obtain the assertion or the inequality $\frac{1}{2\tilde{p^*}} < \frac{1}{n}$, allowing us to proceed as in \textit{Case 2.1} and thus concluding this case and the proof.
\end{proof}

\end{document}